\documentclass[10pt,leqno]{article}
\usepackage[latin1]{inputenc}
\title{Triangulation et cohomologie étale sur une courbe analytique}
\author{\sc Antoine Ducros \footnote{L'auteur est membre des réseaux européens{ \em Arithmetic Algebraic Geometry} et {\em Real Algebraic and Analytic Geometry}.} \\IRMAR Universit\'e de Rennes 1 Campus de
  Beaulieu\\
35042 Rennes CEDEX
FRANCE}
\date{}
\usepackage{mathrsfs}
\usepackage{hyperref} 
\newfont{\cyr}{wncyr10}
\def \cha {\text{{\cyr X}}}
\usepackage[francais]{babel}
\usepackage{amssymb}
\date{}
\input{xypic}

\newcommand{\ra}[1]{\widehat{#1^{rad}}}

\newcommand{\rqt}[1] {\mbox{\rm \bf R}^{q}_{{\bf S}}{#1}}
\newcommand{\rqtc}[1] {\mbox{\rm \bf R}^{q}_{{\bf S},c}{#1} }

\renewcommand{\phi}{\varphi}
\newcommand{\cx}{\got{c}(X)}
\newcommand{\gmod}{^{gmod}}

\newcommand{\got}[1]{{\mathfrak #1}}

\renewcommand{\Bbb}{\mathbb}    
\renewcommand{\cal}{\mathscr}

\newcommand{\gm}{{\Bbb G}_{m}}
\newcommand{\gmu} {{\Bbb G}_{m,1}\an}

\newcommand{\G}{\mbox{\rm G}}
\renewcommand{\H}{\mbox{\rm  H}}
\newcommand{\RR}{{\Bbb R}}
\newcommand{\KK}{{\Bbb K}}
\newcommand{\ZZ}{{\Bbb Z}}
\newcommand{\EE}{{\Bbb E}}
\newcommand{\FF}{{\Bbb F}}
\newcommand{\DD}{{\Bbb D}}
\newcommand{\PP}{{\Bbb P}}

\newcommand{\CC}{{\Bbb C}}

\newcommand{\QQ}{{\Bbb Q}}

\renewcommand{\epsilon}{\varepsilon}

\newcommand{\an}{^{an}}
\newcommand{\zero}{^{\mbox{\tiny o}}}

\newcommand{\zn} {\widehat{\ZZ}_{\neq p}}
\newcommand{\znt} {\widehat{\ZZ}_{\neq p}(1)}

\newcommand{\itb}{\item[$\bullet$]}

\newcommand{\rqq}[1]{\mbox{\rm R}^{q}\pi_{*}#1}
\newcommand{\red}{\widetilde}
\newcommand{\ka}{\widehat{k^{a}}}
\newcommand{\cf}{{\em cf}}
\newcommand{\hres}{{\cal H}}

\newcommand{\deux}[1]{\refstepcounter{cpt}\label{#1}\noindent {\bf (\thesection.\thecpt)}\hspace{.1cm}}
\newcommand{\trois}[1]{\refstepcounter{cptbis}\label{#1}\noindent {\bf
    (\thesection.\thecpt .\thecptbis)}\hspace{.1cm}}
\newcounter{cpt}
\newcounter{cptbis}

\begin{document}
\maketitle
 
\tableofcontents

\setcounter{section}{-1}

\section*{Introduction}

Ce texte est dévolu à l'étude systématique de certains groupes de cohomologie sur une courbe analytique au sens de Berkovich (\cite{brk1},\cite{brk2}), et à leur comparaison avec leurs analogues schématiques lorsqu'on travaille avec l'analytification d'une courbe {\em algébrique}. Il étend, et redémontre par des méthodes moins {\em ad hoc}, les résultats établis par l'auteur dans \cite{duc}. Ceux-ci concernaient le $\H^{3}$ étale d'une courbe $p$-adique à coefficients dans $\mu_{n}^{\otimes 2}$ (pour $n$ premier à $p$) ; dans le présent article on traite le cas d'un corps ultramétrique complet quelconque (sauf au chapitre~\ref{DUAL} où l'on fait l'hypothèse que sa cohomologie possède une bonne dualité, ce qui est vrai par exemple aussi bien pour $\QQ_{p}$ que pour $\CC((t))$ ) et des groupes de cohomologie étale de tous degrés, à coefficients dans n'importe quel module galoisien fini de torsion première à la caractéristique résiduelle. 

\bigskip
\noindent
{\bf Précisions sur les objets étudiés dans l'article.} Désignons par $Z$ ou bien un espace de Berkovich, ou bien un schéma. Le site étale sur $Z$ sera encore noté $Z$ ; autrement dit un {\em faisceau} sur $Z$ sera un faisceau {\em étale}, les groupes $\H^{q}(Z,.)$ seront les groupes de cohomologie étale, {\em etc.} L'espace topologique sous-jacent sera quant à lui désigné par $|Z|$. Soit $\pi$ le morphisme canonique $Z\to |Z|$. Pour tout faisceau $\cal F$ sur $Z$ et tout entier $q$ on notera ${\rm R}^{q}\pi_{*}{\cal F}$ le faisceau {\em sur $|Z|$} associé au préfaisceau $U\mapsto \H^{q}(U,{\cal F})$. On peut aussi décrire ${\rm R}^{q}\pi_{*}$ comme le $q$-ième foncteur dérivé de $\pi_{*}$ ; les  ${\rm R}^{q}\pi_{*}$ mesurent en un sens la \og différence\fg~entre $Z$ et $|Z|$.  On dispose pour tout $q$ et tout $\cal F$ d'une application naturelle $\H^{q}(Z,{\cal F})\to \H^{0}(|Z|,\mbox{R}^{q}\pi_{*}{\cal F})$. Son noyau est formé des classes de cohomologie étale localement triviales sur {\em l'espace topologique $|Z|$} ; elle est surjective lorsque $Z$ est une courbe. 

\bigskip
Soit $\cal X$ une courbe {\em algébrique} sur un corps ultramétrique complet $k$, soit ${\cal X}\an$ son analytifiée à la Berkovich et soit $\cal F$ un faisceau sur $\cal X$. Fixons un entier $q$. Dans le diagramme commutatif $$\diagram \H^{q}({\cal X}\an,{\cal F})\rto& \H^{0}(|{\cal X}\an|,\mbox{R}^{q}\pi_{*}{\cal F})\\ \H^{q}({\cal X},{\cal F})\rto\uto & \H^{0}(|{\cal X}|,\mbox{R}^{q}\pi_{*}{\cal F})\uto \enddiagram$$ les flèches horizontales sont surjectives. Lorsque $\cal F$ est localement constant et de torsion première à la caractéristique résiduelle la flèche verticale de gauche est un isomorphisme d'après Berkovich (\cite{brk2}, cor. 6.3.11) ; celle de droite est donc dans ce cas toujours surjective et l'un des buts de cet article est de donner des conditions suffisantes pour qu'elle soit injective. Son noyau peut être vu comme l'ensemble des classes de $\H^{q}({\cal X},{\cal F})$ qui sont localement triviales sur l'espace topologique $|{\cal X}\an|$ modulo celles qui le sont déjà pour la topologie de Zariski sur $\cal X$. 

\bigskip
\noindent
{\bf Les résultats de {\rm \cite{duc}}.} Ils sont les suivants : prenons pour $k$ un corps {\em local}, c'est-à-dire ou bien une extension finie de $\QQ_{p}$ pour un certain $p$, ou bien un corps de la forme $F((t))$ avec $F$ fini. Soit $Y$ une $k$-courbe analytique lisse et soit $\Delta$ son squelette. C'est un sous-ensemble fermé de $Y$ qui est localement isomorphe à un graphe fini (\cite{brk1}, chap. 4). Orientons-le arbitrairement. Soit $n$ un entier premier à la caractéristique résiduelle ; notons $\mbox{\rm Harm}(\Delta,\ZZ/n)$ le groupe des cochaînes harmoniques sur $\Delta$ à coefficients dans $\ZZ/n$. Le théorème 4.2 de \cite{duc} affirme alors l'existence d'un isomorphisme entre $\H^ {0}(|Y|,\mbox{\rm R}^{3}\pi_{*}\mu_{n}^{\otimes 2})$ et $\mbox{\rm Harm}(\Delta,\ZZ/n)$, isomorphisme qui est construit à l'aide d'une sorte d'évaluation ponctuelle des classes de cohomologie. 

\bigskip
Par ailleurs le théorème 5.2 de \cite{duc} stipule (toujours avec la même hypothèse sur $k$) que si ${\cal X}$ est une $k$-courbe algébrique lisse le morphisme $$ \H^ {0}(|{\cal X}|,\mbox{\rm R}^{3}\pi_{*}\mu_{n}^{\otimes 2})\to \H^ {0}(|{\cal X}\an|,\mbox{\rm R}^{3}\pi_{*}\mu_{n}^{\otimes 2})$$ est un isomorphisme (le groupe de droite pouvant être décrit par le résultat précédent) ; dans le cas où ${\cal X}$ est projective la combinaison de ces deux résultats redonne (par d'autres méthodes) un résultat de Kato.

\subsection*{La démarche suivie ici}

\noindent
{\bf Un rappel du cas algébrique.} Si  $\cal X$ est une courbe {\em algébrique} lisse sur un corps $k$ et si $\cal F$ est un faisceau localement constant sur $k$ et de cardinal premier à la caractéristique de $k$ les groupes $\H^{0}(|{\cal X}|,\mbox{R}^{q}\pi_{*}{\cal F})$ et $\H^{1}(|{\cal X}|,\mbox{R}^{q}\pi_{*}{\cal F})$ sont isomorphes aux groupes de cohomologie d'un complexe à deux termes dit de {\em Bloch-Ogus} (qui se généralise en toute dimension) et qui peut être construit formellement à partir d'une suite exacte de cohomologie à support. {\em À ce sujet, on pourra consulter {\rm \cite{jlct}}.}

\bigskip
\noindent
{\bf Les triangulations.} On cherche à construire un complexe analogue pour une courbe {\em analytique} sur un corps ultramétrique complet $k$. Pour ce faire on définit tout d'abord au chapitre~\ref{TRI} (voir la définition~\ref{TRI}.\ref{deftri}) la notion de {\em triangulation} d'une courbe analytique : {\em grosso modo} une triangulation est la donnée d'un ensemble fermé et discret de points (les sommets) dont les composantes connexes du complémentaires (les arêtes) sont suffisamment élémentaires ; dans le cas algébriquement clos on demande par exemple que ce soient des disques ou des couronnes. Grâce au théorème de réduction semi-stable toute courbe dont le lieu singulier (au sens du critère jacobien) est nulle part dense peut être triangulée (proposition~\ref{TRI}.\ref{triexistsing}). 

\bigskip
\noindent
{\bf Les complexes.} Si $X$ est une courbe analytique munie d'une triangulation $\bf S$ on peut associer à tout faisceau $\cal F$ sur $X$ et tout entier $q$ deux complexes à deux termes  $\rqt{\cal F}$ et $\rqtc{\cal F}$ (voir le~\ref{COH}.\ref{complexe}). Le théorème~\ref{COH}.\ref{cohocomplexe} assure alors que les groupes $\H^{0}(|X|,\mbox{R}^{q}\pi_{*}{\cal F})$ et $\H^{1}(|X|,\mbox{R}^{q}\pi_{*}{\cal F})$ (resp. $\H^{0}_{c}(|X|,\mbox{R}^{q}\pi_{*}{\cal F})$ et $\H^{1}_{c}(|X|,\mbox{R}^{q}\pi_{*}{\cal F})$) sont isomorphes aux groupes de cohomologie de $\rqt{\cal F}$ (resp. $\rqtc{\cal F}$) {\em lorsque $\cal F$ est un faisceau {\em raisonnable}(définition~{\rm \ref{COH}.\ref{defmodfasc}})}.  Le démonstration du théorème~\ref{COH}.\ref{cohocomplexe} repose sur l'existence de suites exactes de cohomologie à support (voir le~\ref{COH}.\ref{introcohsup}) et sur les propriétés cohomologiques des arêtes. 

\bigskip
\noindent
Si $\cal F$ est un module galoisien fini {\em sur le corps $k$} et si sa torsion est première à la caractéristique résiduelle alors $\cal F$ est raisonnable. Sous cette hypothèse on peut (partiellement) expliciter le morphisme entre la cohomologie de $\rqt{\cal F}$ et celle de  $\mbox{R}^{q}\pi_{*}{\cal F}$  ; ce sont les calculs de la fin du chapitre~\ref{COH}. Ceci permet d'arriver à un résultat (théorème~\ref{COMP}.\ref{theocomp}) qui généralise le théorème 5.2 de \cite{duc} : 

\bigskip
\noindent
{\bf Théorème~\ref{COMP}.\ref{theocomp}.} {\em Soit $k$ un corps ultramétrique complet et soit $\cal X$ une courbe algébrique lisse sur $k$. Soit $\cal F$ un module galoisien sur $k$ fini et annulé par un entier $n$ premier à la caractéristique résiduelle. 

\bigskip
\begin{itemize}
\item[$i)$]  Si $\cal X$ a potentiellement bonne réduction alors quel que soit $q$ les flèches $$\H^{0}(|{\cal X}|,{\rm R}^{q+1}\pi_{*}{\cal F})\to  \H^{0}(|{\cal X}\an|,{\rm R}^{q+1}\pi_{*}{\cal F})$$ et $\H^{1}(|{\cal X}|,{\rm R}^{q}\pi_{*}{\cal F})\to \H^{1}(|{\cal X}\an|,{\rm R}^{q}\pi_{*}{\cal F})$ sont des isomorphismes. 

\bigskip
\item[$ii)$] Soit $q$ un entier tel que pour toute extension finie séparable $L$ de $k$ toute classe de $\H^{q}(L,{\cal F})$ soit décomposable.
Alors les flèches $$\H^{0}(|{\cal X}|,{\rm R}^{q+1}\pi_{*}{\cal F})\to  \H^{0}(|{\cal X}\an|,{\rm R}^{q+1}\pi_{*}{\cal F})$$ et $\H^{1}(|{\cal X}|,{\rm R}^{q}\pi_{*}{\cal F})\to \H^{1}(|{\cal X}\an|,{\rm R}^{q}\pi_{*}{\cal F})$ sont des isomorphismes.

\end{itemize}
}

\bigskip
\noindent
{\bf Commentaires.} Fournissons quelques explications : on dit que $\cal X$ a potentiellement bonne réduction si après une extension finie de $k$ elle s'immerge dans une courbe projective et lisse à bonne réduction. Quant à la définition des classes décomposables, elle est relativement technique (définition~\ref{COMP}.\ref{hypoh}). Contentons-nous de dire que la condition donnée en $ii)$ est satisfaite dans de nombreux cas, discutés à la fin du chapitre~\ref{COMP}. Elle est notamment toujours vérifiée quel que soit $\cal F$ lorsque $k$ est local et $q$ égal à $2$ (c'est ce qui redonne le théorème 5.2 de \cite{duc}) ; elle l'est pour tout $\cal F$ et tout $q$ si $k$ est égal à $\CC((t))$ ; elle l'est quel que soit $k$ pour le faisceau $\mu_{n}^{\otimes q}$ lorsque $q$ est égal à 1 ou 2 par le théorème de Merkurjev et Suslin, et pour $q$ au moins égal à $3$ si les résultats annoncés par Rost et Voevodsky sont vrais.

\bigskip
\noindent
{\bf  Au c\oe ur de la démonstration : le principe GAGA pour les fibrés en droite.} Voici {\em grosso modo} comment est établi le théorème ci-dessus. Il s'agit de montrer que certaines classes de cohomologie étale sur $\cal X$, que l'on suppose localement triviales sur $|{\cal X}\an|$, le sont déjà localement pour la topologie de Zariski sur $\cal X$. Cette dernière propriété se teste aux point génériques de $\cal X$ (voir par exemple \cite{jlct}, prop. 2.1.2 et th. 2.2.1). On en déduit que si $L$ est une extension finie séparable de $k$ et si $h$ est une classe de cohomologie de $\cal F$ sur ${\cal X}_{L}$ localement triviale pour la topologie de Zariski alors sa corestriction à la cohomologie de $\cal X$ l'est également. 

\bigskip
Grâce aux calculs explicites de la fin du chapitre~\ref{COH} on montre que les classes étudiées sont sommes de termes de la forme $\mbox{Cor}_ {L/k} (\eta)$ où $L$ est une extension finie séparable de $k$ et où $\eta$ est un élément de la cohomologie de $\cal F$ sur ${\cal X}_{L}$ qui peut être de deux types ; on désigne par $\overline{\cal X}$ une compactification projective de $\cal X$ et par $\partial$ le cobord $\H^{1}(.,\gm)\to \H^{2}(.,\mu_{n})$. 

\bigskip
\begin{itemize}
\item[$1)$] Premier cas : $\eta=\partial(\lambda)_{|{\cal X}_{L}\an}\cup\beta$ où $\lambda$ appartient à $\H^{1}(\overline{\cal X}\an_{L},\gm)$. Comme $\overline{\cal X}$ est projective le principe GAGA assure que $\lambda$ est la classe d'un fibré en droites {\em algébrique}, donc localement trivial pour la topologie de Zariski sur $\overline{\cal X}_{L}$ ; ceci permet de conclure. 

\bigskip
\item[$2)$] Second cas : $\eta=\alpha\cup \beta$ où $\alpha$ appartient à la cohomologie {\em de $L$} à coefficients dans $\cal F$ et où $\beta$ provient de la cohomologie {\em topologique} de $|{\cal X}_{L}\an|$, et plus précisément de $\H^{1}(|{\cal X}_{L}\an|,\ZZ/n)$. Si ${\cal X}$ a potentiellement bonne réduction alors $|{\cal X}_{L}\an|$ est contractile ; dans ce cas $\beta$ est nulle et c'est terminé. Si $\alpha$ est décomposable un calcul direct montre que la classe $\eta$ est somme de termes de la forme $\mbox{Cor}_{M/L}(\omega)$ où $M$ est une extension finie séparable de $L$ et où $\omega$ appartient à la cohomologie de $\cal F$ sur ${\cal X}_{M}$ et s'écrit $\partial(\lambda)_{|{\cal X}_{M}\an}\cup \gamma$ avec $\lambda$ appartenant à $\H^{1}(\overline{\cal X}_{M}\an,\gm)$. On conclut comme pour le premier cas.

\end{itemize}

\bigskip
\noindent
{\bf Optimalité du théorème.} Si $\cal X$ est une {\em courbe de Tate} le théorème~\ref{TATE}.\ref{theoopt} assure que pour tout entier $q$ le noyau de $\H^{0}(|{\cal X}|,{\rm R}^{q+1}\pi_{*}{\cal F})\to  \H^{0}(|{\cal X}\an|,{\rm R}^{q+1}\pi_{*}{\cal F})$ est isomorphe au quotient de $\H^{q}(k,{\cal F})$ par le groupe des classes décomposables. 

\bigskip
\noindent
{\bf Remarque.} Supposons toujours que $\cal X$ est une courbe de Tate. La proposition~\ref{TATE}.\ref{deschzerotate} décrit alors explicitement le groupe $\H^{0}(|{\cal X}\an|,{\rm R}^{q}\pi_{*}{\cal F})$. Ainsi dans certains cas au moins les résultats généraux du chapitre~\ref{COH} permettent des calculs relativement effectifs.

\subsection*{Résultats de dualité} 

\bigskip
\noindent
{\bf Motivation.} On part de la remarque suivante : l'isomorphisme entre les groupes $\H^ {0}(|Y|,\mbox{\rm R}^{3}\pi_{*}\mu_{n}^{\otimes 2})$ et $\mbox{\rm Harm}(\Delta,\ZZ/n)$ établi par le théorème 4.2 de \cite{duc} pour une courbe analytique lisse $Y$ sur un corps local $k$ (pour $n$ premier à la caractéristique résiduelle) peut se réinterpéter comme suit : la dualité de Poincaré fournit pour des raisons formelles une application bilinéaire $$\H^{0}(|Y|,\mbox{\rm R}^{3}\pi_{*}\mu_{n}^{\otimes 2})\times\H^{1}_{c}(|Y|,\ZZ/n)\to \ZZ/n$$ et le résultat ci-dessus peut se traduire en disant que l'accouplement en question est une dualité parfaite. 

\bigskip
\noindent
{\bf Le formalisme.} Au chapitre~\ref{DUAL} on fixe un entier $n$ premier à la caractéristique résiduelle et on fait l'hypothèse que la cohomologie des faisceaux  en $\ZZ/n$-modules sur $k$ possède un module dualisant $\DD$ en un certain degré $d$. Les hypothèses sont détaillées au~\ref{DUAL}.\ref{contextedual} ; si $k$ est local elles sont vérifiées avec $d$ égal à $2$ en prenant pour $\DD$ le faisceau $\mu_{n}$, si $k$ est égal à $\CC((t))$ elles le sont avec $d$ égal à $1$ et en prenant pour $\DD$ le faisceau constant $\ZZ/n$. Soit $X$ une $k$-courbe lisse et soit  $\cal F$ un faisceau raisonnable sur $X$. Une conséquence formelle des hypothèses faites sur $k$ et de la dualité de Poincaré sur un corps algébriquement clos établie par Berkovich au paragraphe 7.3 de \cite{brk2} est l'existence pour tout $q$  de deux accouplements $$\H^{0}(|X|,\mbox{R}^{q}\pi_{*}{\cal F})\times \H^{1}_{c}(|X|,\mbox{R}^{d+1-q}\pi_{*}\breve{\cal F})\to \ZZ/n$$ et $$\H^{0}_{c}(|X|,\mbox{R}^{q}\pi_{*}{\cal F})\times \H^{1}(|X|,\mbox{R}^{d+1-q}\pi_{*}\breve{\cal F})\to \ZZ/n$$ où les groupes à support compact sont discrets et ceux sans support profinis  (corollaire~\ref{DUAL}.\ref{conclucohoprofin}), et où la notation $\breve{\cal F}$ désigne le faisceau ${\cal F}^{\vee}\otimes\DD\otimes\mu_{n}$. 

\bigskip
\noindent
{\bf Les deux théorèmes de dualité.} Les théorèmes~\ref{DUAL}.\ref{maintheodua1} et~~\ref{DUAL}.\ref{maintheodua2} assurent alors que pour chacun des deux accouplements le noyau à droite est trivial tandis que le noyau à gauche est le sous-groupe des classes gratte-ciel, c'est-à-dire {\em à support dans un ensemble discret. }

\bigskip
Leurs démonstrations reposent sur la dualité de Poincaré, sur la description des groupes en jeu à l'aide des complexes introduits au chapitre~\ref{COH} et sur les propriétés cohomologiques des arêtes d'une triangulation. 

\bigskip
\noindent
{\bf Les classes gratte-ciel.} L'étude de ces classes fait l'objet du chapitre~\ref{GRAT} (sans hypothèse de dualité sur $k$) ; on y démontre que pour une courbe lisse munie d'une triangulation $\bf S$ leur support est forcément contenu dans le sous-ensemble ${\bf S}_{(2)}$ des points de type $(2)$ de $\bf S$ (proposition~\ref{GRAT}.\ref{ouvtri}). Si $P$ est un point de ${\bf S}_{(2)}$ le corps $\red{\hres(P)}$ est de degré de transcendance $1$ sur $\red{k}$. {\em Si $\cal F$ est non ramifié}, c'est-à-dire provient d'un faisceau sur le corps résiduel $\red{k}$ de $k$, alors la fibre en tout point $P$ de ${\bf S}_{(2)}$ du groupe des classes gratte-ciel  est somme directe de groupes de la forme $\cha^{q-j}(\red{\hres(P)},{\cal F}(-j))$ où $j$ est un entier positif. C'est la proposition~\ref{GRAT}.\ref{gptype2} ; la notation $\cha$ fait référence à l'ensemble des $\red{k}$-places de $\red{\hres(P)}$ (elle est précisément définie au~\ref{GRAT}.\ref{ouvuxp}) et ${\cal F}(-j)$ désigne le tordu à la Tate de $\cal F$. La présence de classes gratte-ciel non triviales dans le cas où $\cal F$ est non ramifié est donc reliée à l'existence de contre-exemples au principe de Hasse pour certains groupes de cohomologie sur les corps de fonctions $\red{\hres(P)}$ où $P$ appartient à ${\bf S}_{(2)}$. 

\subsection*{Exemples (voir le chapitre~\ref{RE})}

\noindent
{\bf Le théorème 4.2 de \cite{duc}}. Si $k$ est local la nullité de $\cha^{2}(.,\mu_{n})$ sur le corps des fonctions d'une courbe sur un corps fini permet de montrer que le sous-groupe de $\H^{0}(|X|,\mbox{\rm R}^{3}\pi_{*}\mu_{n}^{\otimes 2})$  formé des classes gratte-ciel est trivial ; l'accouplement défini ci-dessus est donc une dualité parfaite, et l'on retrouve ainsi le théorème 4.2 de \cite{duc}. 

\bigskip
\noindent
{\bf La dualité de Lichtenbaum}. Supposons toujours $k$ local. Alors la nullité des groupes $\cha^{2}(.,\mu_{n})$ et $\cha^{1}(.,\ZZ/n)$ sur le corps des fonctions d'une courbe sur un corps fini  permet de  montrer que le sous-groupe de $\H^{0}(|X|,\mbox{\rm R}^{2}\pi_{*}\mu_{n})$  formé des classes gratte-ciel est trivial. En conséquence on dispose d'une dualité parfaite entre $\H^{0}(|X|,\mbox{\rm R}^{2}\pi_{*}\mu_{n})$ et $\H^{1}_{c}(|X|,\mbox{\rm R}^{1}\pi_{*}\mu_{n})$. Ceci couplé au théorème de comparaison~\ref{COMP}.\ref{theocomp} permet dans le cas d'une $k$-courbe algébrique projective et lisse $\cal X$ de retrouver la dualité de Lichtenbaum (\cite{lic}) entre $_{n}\mbox{\rm Br}\;{\cal X}$  et $(\mbox{\rm Pic}\;{\cal X})/n$. 

\bigskip
\noindent
{\bf Dualité parfaite sur les ouverts de courbes de Mumford}. Pour tout corps $F$,  pour  tout module galoisien fini $\cal G$ sur $F$ de torsion première à la caractéristique de $F$ et pour tout entier $j$ le groupe $\cha^{j}(F(t),{\cal G})$ est nul (\cf. \cite{jps}, \S 4 p. 122) . Dès lors si pour tout $P$ appartenant à ${\bf S}_{(2)}$ le corps $\red{\hres(P)}$ est transcendant pur sur une extension finie de $\red{k}$ l'accouplement construit ci-dessus est une dualité parfaite dès que $\cal F$ est non ramifié (et sans autre hypothèse sur $k$ que l'existence du module dualisant $\DD$). {\em On est par exemple dans cette situation si $X$ est un ouvert de l'analytifiée d'une courbe de Mumford.}

\subsection*{Remerciements} 

Cet article doit beaucoup à un séjour d'un mois à Rehovot, en mai 2004, dans le cadre accueillant et propice à la réflexion de l'institut Weizmann. Je tiens à remercier cet organisme pour son hospitalité et à faire part de ma gratitude à Vladimir Berkovich, à l'origine de l'invitation, avec qui j'ai eu de longues et fructueuses discussions.

\setcounter{section}{0}

\subsection*{Quelques conventions relatives à la cohomologie} 

\bigskip
\begin{itemize}

\itb Si $Z$ est un espace de Berkovich ou bien un schéma on notera encore $Z$ le site étale correspondant et $|Z|$ l'espace topologique sous-jacent. 

\bigskip
\itb Quels que soient les sites en jeu tous les résultats concernant des groupes $\H^{q}$ ou des faisceaux $\mbox{R}^{q}\pi_{*}$ et présentés comme valables pour tout $q$ doivent être considérés comme vrais pour tout entier {\em relatif} $q$ avec la convention que les $\H^{q}$ et $\mbox{R}^{q}\pi_{*}$ sont toujours nuls dès que $q$ est strictement négatif. Cette remarque présente un intérêt lorsqu'un même énoncé (comme l'existence d'une suite exacte ou d'un isomorphisme) met en jeu de la cohomologie en degrés $q-1$, $q$, $q+1$, {\em etc.}, auquel cas l'information apportée par la considération des degrés négatifs est non vide. 

\bigskip
\itb En ce qui concerne la cohomologie étale (analytique ou schématique) si $f$ (resp. $n$) est une fonction (resp. un entier)  inversible sur le site annelé considéré et si le contexte indique clairement que l'on travaille avec des faisceaux en $\ZZ/n$-modules la notation $(f)$ désignera l'image de $f$ dans $\H^{1}(.,\mu_{n})$ par le cobord de la suite de Kummer, c'est-à-dire encore la classe du $\mu_{n}$-torseur $$U\mapsto \{g\in\gm(U)\;\mbox{tq}\;g^{n}=f\}.$$

\bigskip
\itb Pour tout entier $n$ et tout entier {\em négatif} $i$ le faisceau $\mu_{n}^{\otimes i}$ (sur un site annelé donné) sera {\em par définition} le faisceau ${\cal Hom}(\mu_{n}^{\otimes (-i)},\ZZ/n)$. Si $\cal F$ est un faisceau en $\ZZ/n$-modules le faisceau ${\cal F}\otimes\mu_{n}^{\otimes i}$ sera noté ${\cal F}(i)$ pour tout entier relatif $i$. 
\end{itemize} 

\section{La cohomologie de certains corps henséliens}\label{HENS}

\subsection*{Une décomposition de certains groupes de cohomologie} 

\bigskip
{\em Dans tout l'article le corps résiduel d'un corps valué $L$ sera systématiquement noté $\red{L}$.}

\bigskip
\deux{gval} Soit $K$ un corps hensélien pour une valuation $|\;|$ ; la notation est donc multiplicative  {\em et on ne suppose pas que le groupe des valeurs se plonge dans $\RR^{*}_{+}$.} Soit $K^{s}$ une clôture séparable de $K$. On note $p$ l'exposant caractéristique de $\red{K}$ et $\red{K}^{s}$ la fermeture séparable de $\red{K}$ dans $\red{K^{s}}$. Soit $G$ (resp. $\Pi$) le groupe de Galois de $\red{K}^{s}$ sur $\red{K}$ (resp. le quotient du groupe de Galois de $K^{s}$ sur $K$ par son sous-groupe de ramification). On dispose d'une suite exacte $$1\to \mbox{Hom}(|(K^{s})^{*}|/|K^{*}|,\red{K^{s}}^{*})\to \Pi\to G\to 1.$$ Le lemme de Zorn assure l'existence d'une extension modérément ramifiée $M$ de $K$ incluse dans $K^{s}$, telle que $\red{M}$ soit égal à $\red{K}$, et maximale pour cette propriété. Le choix de $M$ fournit une section de la suite exacte ci-dessus ; le groupe $|(K^{s})^{*}|/|M^{*}|$ est de torsion $p$-primaire.

\bigskip
\deux{isoqlzl} On choisit un isomorphisme entre $|(K^{s})^{*}|/|K^{*}|$, qui est divisible et de torsion, et $\bigoplus (\QQ_ {l}/\ZZ_{l})^{(I_{l})}$, où $I_{l}$ est pour tout nombre premier $l$ un ensemble d'indices ; notons que $|M^{*}|/|K^{*}|$ s'identifie par ce biais à $\bigoplus \limits_{l\neq p} (\QQ_ {l}/\ZZ_{l})^{(I_{l})}$. Quant au $G$-module $(\red{K^{s}}^{*})_{tors}$ il est {\em canoniquement} isomorphe à $\bigoplus\limits_{l \neq p} (\QQ_{l}/\ZZ_{l})(1)$. En conséquence les choix faits définissent un isomorphisme $\Pi\simeq \left(\prod\limits_{l \neq p} \ZZ_{l}(1)^{I_{l}}\right)\rtimes G$. On note $I$ la réunion disjointe $\coprod \limits_{l\neq p}I_{l}$.  

\bigskip
Si $\cal F$ est un faisceau étale provenant d'un $\Pi$-module de torsion première à $p$ alors pour tout entier $q$ le groupe $\H^{q}(K,{\cal F})$ s'identifie à $\H^{q}(\Pi, {\cal F})$ puisque le sous-groupe de ramification de $\mbox{Gal}(K^{s}/K)$ en est un pro-$p$-sous-groupe distingué.

\bigskip
\deux{h1cocycle} {\em On fixe un nombre premier $l$ différent de $p$}. Soit $n$ un entier. On identifie le sous-groupe de $l^{n}$-torsion de $|(K^{s})|^{*}/|K^{*}|$ à $(\ZZ/l^{n})^{(I_{l})}$. Soit $(a_{i})_{i\in I_{l}}$ une famille d'éléments de $K^{*}$ telle que  $(|a_{i}|^{1/l^{n}})_{i}$ forme une base du $\ZZ/l^{n}$-module libre $(\ZZ/l^{n})^{(I_{l})}$. 

\bigskip
\deux{basecan} {\bf Remarque.} Notons que si $(|a_{i}|^{1/l^{n}})_{i}$ est la base {\em canonique} de $(\ZZ/l^{n})^{(I_{l})}$ alors la classe $(a_{i})$ de $\H^{1}(K,\mu_{l^{n}})$ est donnée pour tout $i$ par le cocycle qui envoie un élément $(u,g)$ de $\left(\prod\limits_{l \neq p} \ZZ_{l}(1)^{I_{l}}\right)\rtimes G$ sur l'image du $i$-ième facteur de $u$ par la projection $$\ZZ_{l}(1)\to (\ZZ_{l}(1)/l^{n})\simeq \mu_{l^{n}}.$$

\bigskip
\deux{definicup} Munissons $I_{l}$ d'un ordre arbitraire. Pour toute partie finie $J$ de $I_{l}$ on note $(a_{J})$ l'élément du groupe $\H^{|J|}(K,\ZZ/l^{n}(|J|))$ égal au cup-produit des $(a_{i})$  où $i$ parcourt $J$ et où les $a_{i}$ sont rangés dans le sens des $i$ croissants. 

\bigskip
\deux{decompocup} {\bf Lemme.} {\em Soit $\cal F$ un faisceau étale sur $K$ donné par un $G$-module fini annulé par $l^{n}$. Soit $q$ un entier. Tout élément du groupe $\H^{q}(K,{\cal F})$ s'écrit alors de manière unique comme une somme $$\sum(a_{J})\cup h_{J}$$ où $J$ parcourt l'ensemble des parties de $I_{l}$ de cardinal majoré par $q$, où $h_{J}$ appartient pour tout $J$ au groupe $\H^{q-|J|}(\red{K},{\cal F}(-|J|))$ et où les $h_{J}$ sont presque tous nuls.} 

\bigskip
{\em Démonstration.} Par multilinéarité du cup-produit il suffit de le montrer dans le cas où $(|a_{i}|^{1/l^{n}})_{i}$ est la base canonique de $(\ZZ/l^{n})^{(I_{l})}$. Munissons $\ZZ^{(I)}$ de l'ordre lexicographique induit par un ordre quelconque sur $I$ et le corps $\red{K}(T_{i})_{i\in I}$ (où les $T_{i}$ sont des indéterminées) de la valuation qui envoie pour tout $i$ l'indéterminée $T_{i}$ sur l'élément de $\ZZ^{(I)}$ valant $1$ en $i$ et $0$ ailleurs. Soit $\red{K}(T_{i})_{i}^{h}$ le hensélisé de $\red{K}(T_{i})_{i}$ pour ladite valuation. Soit $\KK$ une extension algébrique de $\red{K}(T_{i})_{i}^{h}$ obtenue par adjonction pour tout $i$ un système compatible de racines de $T_{i}$ d'ordre premier à $\mu$, où $\mu$ est le nombre premier tel que $i$ appartienne à $I_{\mu}$. 

\setcounter{cptbis}{0}
\bigskip
\trois{vraicorpspart} {\em Le lemme est vrai pour le corps $\KK$, avec $a_{i}$ égal à $T_{i}$ pour tout $i$ dans $I_{l}$.} Pour le voir, on se ramène aussitôt, par un argument de limite inductive, au cas d'un nombre fini d'indéterminées. Ensuite une récurrence immédiate autorise à supposer qu'il n'y en a qu'une, que l'on note $T$. L'ensemble $I$ est alors un singleton et est égal à $I_{\mu}$ pour un nombre premier $\mu$ fixé (qui est différent de $p$). Le groupe $|\KK^{*}|$ est dans ces conditions isomorphe à $\ZZ_{(\mu)}$ et le groupe d'inertie à $\ZZ_{\mu}$. Si $\mu$ est différent de $l$ alors $I_{l}$ est vide et le groupe d'inertie étant un pro-$\mu$-groupe on dispose d'un isomorphisme $\H^{q}(\KK,{\cal F})\simeq \H^{q}(\red{K},{\cal F})$ ; le lemme est donc vrai dans ce cas. 

\bigskip
Si $\mu$ est égal à $l$ alors utilise là encore un argument de limite inductive pour modifier à nouveau $\KK$ et faire l'hypothèse que c'est le hensélisé de $\red{K}(T^{1/e})$ pour un certain $e$ premier à $l$. C'est en particulier un corps hensélien pour une valuation discrète dont $T^{1/e}$ est une uniformisante. Il est dès lors classique et bien connu que toute classe de $\H^{q}(\KK,{\cal F})$ a une unique écriture sous la forme $ \eta_{0}+(T^{1/e})\cup \eta_{1}$ où $\eta_{0}$ appartient à $\H^{q}(\red{K},{\cal F})$ et $\eta_{1}$ à $\H^{q-1}(\red{K},{\cal F}(-1))$. Comme $e$ est inversible modulo $l^{n}$ on a également existence et unicité de l'écriture sous la forme $h_{0}+(T)\cup h_{1}$ où  $h_{0}$ appartient à $\H^{q}(\red{K},{\cal F})$ et $h_{1}$ à $\H^{q-1}(\red{K},{\cal F}(-1))$ et là encore le lemme est vrai.

\bigskip
\trois{conclugalois} Le lemme est donc établi lorsque $K$ est égal à $\KK$. Or le~\ref{HENS}.\ref{basecan} permet d'en traduire l'énoncé, indépendamment du corps $K$, uniquement en termes de cohomologie du groupe $\left(\prod\limits_{l \neq p} \ZZ_{l}(1)^{I_{l}}\right)\rtimes G$. La validité du résultat souhaité sur le corps $\KK$ l'implique donc dans le cas général.~$\Box$

\subsection*{Composition avec une valuation discrète}

\deux{compvaldisc} On conserve les notations du~\ref{HENS}.\ref{gval} ; on suppose de plus que $\red{K}$ est muni d'une valuation {\em discrète} $v$ dont on note $\kappa$ le corps résiduel. Soit $\bf K$ le hensélisé de $K$ pour la composée de $|\;|$ et $v$. Le corps résiduel de $\bf K$ est égal à $\kappa$. Ce qu'on notera $p$ sera désormais l'exposant caractéristique {\em du corps $\kappa$}. Soit $l$ un nombre premier différent de $p$, soit $n$ un entier et soit $\cal F$ un $\kappa$-faisceau étale donné par un module galoisien fini annulé par $l^{n}$. Considérons une famille $(a_{i})$ de $K^{*}$ sujette à la même condition que celle du~\ref{HENS}.\ref{h1cocycle}. On ordonne $I_{l}$ et on définit $(a_{J})$ pour toute partie finie $J$ de $I_{l}$ comme au~\ref{HENS}.\ref{definicup}. Soit $\tau$ un élément de $K\zero$ tel que $\red{\tau}$ soit une uniformisante de $v$. 

\bigskip
\deux{compcohomo} Soit $q$ un entier et soit $h$ un élément de $\H^{q}(K,{\cal F})$. Le lemme~\ref{HENS}.\ref{decompocup} affirme que $h$ a une unique écriture sous la forme $\sum (a_{J})\cup h_{J}$ où $h_{J}$ appartient pour tout $J$ à $\H^{q-|J|}(\red{K},{\cal F}(-|J|))$. 

\bigskip
Soit $\red{K}^{h}$ le hensélisé de $\red{K}$ pour la valuation $v$. Fixons $J$. L'image de la classe $h_{J}$ dans $\H^{q-|J|}(\red{K}^{h},{\cal F}(-|J|))$ a (toujours d'après le lemme~\ref{HENS}.\ref{decompocup}) une unique écriture sous la forme $h_{J,0}+(\red{\tau})\cup h_{J,1}$ où $h_{J,0}$ appartient à $\H^{q-|J|}(\kappa,{\cal F}(-|J|))$ et  $h_{J,1}$ à $\H^{q-|J|-1}(\kappa,{\cal F}(-|J|-1))$. 

\bigskip
\deux{hasseloc} On en déduit (en décomposant la valuation de $\bf K$ à l'aide de celle de $\red{K}^{h}$ et en appliquant encore le lemme~\ref{HENS}.\ref{decompocup}), que $h_{|{\bf K}}$ s'écrit $$\sum\limits_{J} (a_{J})\cup h_{J,0}+(a_{J})\cup (\tau)\cup h_{J,1}.$$ Toujours d'après le lemme~\ref{HENS}.\ref{decompocup} la classe $h_{|{\bf K}}$ est nulle si et seulement si toutes les $h_{J,0}$ et $h_{J,1}$ sont nulles, donc si et seulement si la restriction de chacune des $h_{J}$ à $\red{K}^{h}$ est triviale.

\subsection*{Extension de rang rationnel relatif égal à $1$} 

\deux{fixenotation} {\bf Contexte et notations.} Soit $F$ un corps hensélien pour une valuation $|\;|$. On fixe une clôture séparable $F^{s}$ de $F$, on note $\Gamma$ le groupe ordonné $|(F^{s})^{*}|$ et $G$ le groupe de Galois de $F^{s}$ sur $F$. Soit $K$ un corps valué hensélien prolongeant $F$ dans lequel ce dernier est algébriquement clos. Soit $L$ le corps $K\otimes_{F}F^{s}$. On suppose  que $\red{K}$ est fini (et donc nécessairement purement inséparable) sur $\red{F}$ et que le groupe $|L^{*}|/\Gamma$ est libre de rang 1. Le quotient de $|K^{*}|/|F^{*}|$ par sa torsion s'injecte dans $|L^{*}|/\Gamma$ et s'identifie donc (le conoyau de ce plongement étant de torsion) à $(|L^{*}|/\Gamma)^{m}$ pour un certain entier $m$ strictement positif  que l'on appellera la {\em pathologie} de l'extension $K/F$. On dira que l'on a {\em orienté} $K$ si l'on choisit un générateur (parmi les deux possibles) de $|K^{*}|/|F^{*}|$ modulo sa torsion. Il revient au même de se donner un générateur de $|L^{*}|/\Gamma$. 

\bigskip
\deux{exemplepatho} {\bf Exemple de pathologie non triviale.} Soit $k$ un corps complet pour une valuation discrète et de caractéristique nulle. Supposons que le corps résiduel $\red{k}$ est de caractéristique $2$ et non parfait. Munissons $k$ de la valeur absolue compatible avec sa valuation et telle que $|2|$ soit égale à $1/2$. Soit $a$ appartenant à $k\zero$ tel que $\red{a}$ ne soit pas un carré. Soit $r$ un réel strictement compris entre $1/2$ et $1$ et n'appartenant pas à $\sqrt{|k^{*}|}$.  Soit $E$ le corps $k\{T/r\}(\sqrt{T-a})$. On a alors $\red{E}=\red{k}(\sqrt{\red{a}})$ et $|E^{*}|=|k^{*}|\oplus r^{\ZZ}$. L'élément $\sqrt{T-a}-\sqrt{a}$ de $E\otimes_{k}k(\sqrt{a})$ est de valeur absolue égale à $\sqrt{r}$ ; la pathologie de l'extension $E/k$ est en conséquence non triviale. 

\bigskip
\deux{cohomodereram} On revient aux notations du~\ref{HENS}.\ref{fixenotation}. Soit $L^{s}$ une clôture séparable de $L$ ; c'est également une clôture séparable de $K$. Comme le corps résiduel $\red{L}$ est algébriquement clos le groupe $\mbox{Gal}\;(L^{s}/L)$ est son propre sous-groupe d'inertie. Soit $W$ son groupe de ramification et soit $\pi$ le quotient de $\mbox{Gal}\;(L^{s}/L)$ par $W$. Il est canoniquement isomorphe à $\mbox{Hom}(|(L^{s})^{*}|/|L^{*}|,\red{L^{s}}^{*})$. 

\bigskip
Supposons que l'on a orienté $K$. Soit $\xi$ appartenant à $K^{*}$ dont l'image dans $|K^{*}|/|F^{*}|$  modulo sa torsion est le générateur choisi. On dira qu'un tel $\xi$ {\em représente} l'orientation choisie. Le groupe $\Gamma$ étant divisible il existe $\tau$ appartenant à $L^{*}$ tel que $|\tau|^{m}$ soit égal à $|\xi|$ ; l'image de $|\tau|$ dans $|L^{*}|/\Gamma$ en est un générateur. Comme $\Gamma$ est divisible le groupe $|L^{*}|$ est égal à $|\tau|^{\ZZ}\oplus \Gamma$, et $|(L^{s})^{*})|/|L^{*}|$ est en conséquence isomorphe à $\QQ/\ZZ$. Le groupe de ramification modérée est donc isomorphe au produit des $\ZZ_{l}$ pour $l$ premier à $p$, produit que l'on notera simplement $\zn$.  Pour tout entier $n$ premier à $p$ l'unique quotient d'ordre $n$ de $\zn$ correspond à l'extension obtenue en adjoignant à $L$ n'importe quelle racine $n$-ième de n'importe quel élément de valuation égale à $|\tau|$. Comme $\zn$ est abélien l'isomorphisme entre ce dernier et le groupe de ramification modérée est indépendant du choix de la clôture séparable $L^{s}$, {\em mais pas du choix de l'orientation : si on modifie cette dernière l'isomorphisme en question est changé en son opposé}. 

\bigskip
\deux{cohoenbas} Le quotient $\Pi$ du groupe de Galois de $L^{s}$ sur $K$ par $W$ s'insère donc indépendamment du choix de $L^{s}$ dans une suite exacte courte $$1\to \pi \to \Pi \to G \to 1.$$ On dira qu'une $K$-algèbre étale est {\em géométriquement modérément ramifiée} si elle peut être définie par un $\Pi$-ensemble. Une $K$-algèbre étale est géométriquement modérément ramifiée si et seulement si elle devient isomorphe, après tensorisation par $F^{s}$, à un produit d'extensions modérément ramifiées de $L$. 

\bigskip
Comme $\pi$ est abélien la suite exacte ci-dessus induit une vraie action de $G$ sur $\pi$. Un choix d'orientation sur $K$ fournit d'après ce qui précède un isomorphisme entre $\pi$ et $\zn$, et donc une action de $G$ sur $\zn$. Un calcul immédiat montre que $\zn$, vu comme $G$-groupe abélien profini {\em via} cette action, est isomorphe à $\znt$. 

\bigskip
\deux{suitespec} Soit $M$ un $\Pi$-module fini de cardinal premier à $p$. Il définit en particulier un faisceau étale sur $K$, noté encore $M$. Comme le groupe de ramification est un pro-$p$-groupe distingué dans $\mbox{Gal}(L^{s}/K)$ la flèche naturelle  $\H^{q}(\Pi,M)\to \H^{q}(K,M)$ est un isomorphisme. 

\bigskip
{\em On fixe une orientation sur $K$ et on s'en donne un représentant $\xi$. On rappelle que $m$ désigne la pathologie de $K$ et on choisit un élément $\tau$ de $(K^{s})^{*}$ vérifiant l'égalité $|\tau|^{m}=|\xi|$.}

\bigskip
\deux{pathopprim} {\bf Proposition.} {\em La pathologie $m$ de $K$ est une puissance de $p$.}

\bigskip
{\em Démonstration.} La flèche naturelle de $\mbox{Gal} \;(L/K)$ dans $\mbox{Gal} \;(F^{s}/F)$ est un isomorphisme. Comme  $\red{K}$ est une extension finie purement inséparable de $\red{F}$ le sous-groupe d'inertie de $\mbox{Gal} \;(L/K)$ est le même que celui de $\mbox{Gal} \;(F^{s}/F)$, et les sous-groupes de ramifications sont également les mêmes (le groupe de ramification étant simplement l'unique $p$-sous-groupe de Sylow du groupe d'inertie). Soit $E$ l'extension modérément ramifiée maximale de $F$ dans $F^{s}$. D'après ce qu'on vient de voir $K\otimes_{F}E$ est l'extension modérément ramifiée maximale de $K$ dans $L$. Le morphisme $\mbox{Gal}\;((K\otimes_{F}E)/K)\to \mbox{Gal}\;(E/F)$ est un isomorphisme. Or le premier de ces groupes s'identifie à $\mbox{Hom}(|(K\otimes_{F}E)^{*}|/|K^{*}|,(\red{K\otimes_{F}E})^{*})$, et même plus précisément à $\mbox{Hom}(|(K\otimes_{F}E)^{*}|/|K^{*}|,\red{E}^{*})$ puisque $\red{K\otimes_{F}E}$ est une extension purement inséparable de $\red{E}$, d'où il découle que la torsion de  $(\red{K\otimes_{F}E}^{*})$ est incluse dans $\red{E}^{*}$. Le second de ces groupes est canoniquement isomorphe à $\mbox{Hom}(|E^{*}|/|F^{*}|,\red{E}^{*})$. 

Soit $s$ le plus grand entier strictement positif tel que $|\xi|^{1/s}$ appartienne à $|(K\otimes_{F}E)^{*}|$. Comme $K\otimes_{F}E$ est modérément ramifié sur $K$ l'entier $s$ est premier à $p$. Comme $F^{s}/E$ est une $p$-extension l'entier $m$ est de la forme $p^{r}s$ pour un certain $r$. Fixons une racine primitive $s$-ième de l'unité dans $\red{E}$ et notons $\phi$ le morphisme de $|\tau|^{p^{r}\ZZ}\oplus\Gamma$ dans $\red{E}^{*}$ qui est nul sur $\Gamma$ et envoie $|\tau|^{p^{r}}$ sur la racine en question. Alors $\phi$ induit un morphisme de $|(K\otimes_{F}E)^{*}|/|K^{*}|$ vers $\red{E}^{*}$ dont la restriction à $|E^{*}|/|F^{*}|$ est triviale. Le morphisme lui-même est trivial en vertu de ce qui précède, ce qui montre que $s$ est égal à 1 et donc que $m$ est égal à $p^{r}$.~$\Box$

\bigskip
\deux{residucup} {\bf Proposition.} {\em Soit $n$ un entier premier à $p$ et soit $M$ un $G$-module fini annulé par $n$. Soit $q$ un entier. L'application $(h,\eta)\mapsto h+(\xi)\cup \eta$ induit un isomorphisme $$\H^{q}(F,M)\oplus\H^{q-1}(F,M(-1))\simeq \H^{q}(K,M).$$}

\bigskip
{\em Démonstration.} Il suffit de le démontrer dans le cas où $n$ est de la forme $l^{d}$ où $l$ est un nombre premier différent de $p$ et $d$ un entier. Plaçons-nous sous cette hypothèse et soit $(a_{i})_{i }$ une famille d'éléments de $F^{*}$ telle que $(|a_{i}|)_{i}$ forme une $\ZZ/l^{d}$-base de la $l^{d}$-torsion de $|(F^{s})^{*}|/|F^{*}|$. La pathologie étant une puissance de $p$ est en particulier première à $l$, et de ce fait la réunion de $(|a_{i}|)_{i}$ et de $|\xi|$ forme une $\ZZ/l^{d}$-base de $|(K^{s})^{*}|/|K^{*}|$. On peut dès lors appliquer la proposition~\ref{HENS}.\ref{decompocup} : si l'on munit l'ensemble $I$ des indices d'un ordre arbitraire alors toute classe de $\H^{q}(K,M)$ a une unique décomposition sous la forme $$h=\sum (a_{J})\cup h_{J,0}+\sum (\xi)\cup (a_{J})\cup h_{J,1}$$ où $J$ parcourt l'ensemble des parties finies de $I$, et où pour $J$ fixée les classes $h_{J,0}$ et $h_{J,1}$ appartiennent respectivement aux groupes $\H^{q-|J|}(\red{K},M(-|J|))$ et $\H^{q-|J|-1}(\red{K}, M(-|J|-1))$ ; la notation $(a_{J})$ a quant à elle été introduite au~\ref{HENS}.\ref{definicup}. D'autre part, toujours par la proposition~\ref{HENS}.\ref{decompocup}, on peut écrire pour tout entier $\delta$ toute classe $\lambda$ de $\H^{\delta}(F,M)$ d'une unique manière sous la forme $\sum (a_{J})\cup \lambda_{J}$ où $\lambda_{J}$ est pour toute partie finie $J$ de $I$ un élément de $\H^{\delta-|J|}(\red{F},M(-|J|))$. Mais $\red{F}\hookrightarrow \red{K}$ est finie et purement inséparable ; en conséquence les flèches de restriction de la cohomologie de $\red{F}$ vers celle de $\red{K}$ sont des isomorphismes. Le résultat cherché s'en déduit aussitôt.~$\Box$

\section{Pseudo-couronnes et pseudo-disques}\label{PSEU}

\setcounter{cpt}{0}

\subsection*{Définitions générales et premières propriétés} 

\deux{intro} {\bf Notations.} Si $L$ est un corps ultramétrique complet on notera $L^{rad}$ sa clôture radicielle et $\ra{L}$ le complété de cette dernière. On se donne {\bf pour toute la suite} un corps $k$ complet pour une valeur absolue ultramétrique non triviale $|\; .\;|$. On désignera par $p$ l'exposant caractéristique du corps résiduel $\red{k}$. On fixe une clôture algébrique $k^{a}$ de $k$ et l'on note $\ka$ son complété. On désigne par $k^{s}$ la fermeture séparable de $k$ dans $k^{a}$. La notion d'espace $k$-analytique sera à prendre au sens de Berkovich (\cite{brk1},\cite{brk2}). Si $X$ est un espace $k$-analytique dont les domaines affinoïdes sont réduits on notera $\cx$ l'anneau des sections globales du faisceau des fonctions constantes sur $X$. Lorsque $X$ est connexe et non vide $\cx$ est un corps de dimension finie sur $k$, séparable si $X$ est géométriquement réduit. Une {\em $k$-courbe analytique} (ou simplement une {\em $k$-courbe} s'il n'y a pas d'ambiguïté) sera un espace $k$-analytique séparé, paracompact et purement de dimension 1. 

\bigskip
\deux{courdisque} On appellera {\em $k$-disque} (resp. {\em $k$-couronne}) (ou simplement disque ou couronne s'il n'y a pas d'ambiguïté sur le corps de base) tout espace $k$-analytique isomorphe à un domaine de la droite affine défini par une inégalité de la forme $|t|<R$ (resp. $r<|t|<R$) où $R$ est strictement positif et peut éventuellement être pris égal à $+\infty$, où $r$ est positif ou nul et où $r$ est strictement inférieur à $R$. Une fonction induisant un tel isomorphisme sera appelé une {\em fonction coordonnée}. 

\bigskip
Si $X$ est une couronne et si $Y$ est un ouvert de $X$ on dira que $Y$ est une {\em sous-couronne} de $X$ si c'est une couronne et si la restriction à $Y$ de toute fonction coordonnée de $X$ en est une fonction coordonnée (il suffit que ce soit vrai pour {\em une} fonction coordonnée). Si $X$ est un disque et si $Y$ est un ouvert de $X$ on dira que $Y$ est une {\em sous-couronne} de $X$ si c'est une couronne et si la restriction à $Y$ de toute fonction coordonnée de $X$ en est une fonction coordonnée ayant même borne supérieure en norme sur $X$ et sur $Y$ (il suffit que ce soit vrai pour {\em une} fonction coordonnée). Les notions de couronne, de disque, de sous-couronne, de fonction coordonnée...sont trivialement stables par changement de corps de base. 

\setcounter{cptbis}{0}

\bigskip
\trois{gmcouronne} {\bf  Remarque.} On considère donc ici la droite affine comme un disque, et un disque ouvert épointé ou le groupe multiplicatif comme des couronnes ; ces abus de langage sont justifiés par le fait que la théorie des faisceaux étales "raisonnables" sur la droite affine (resp. le disque ouvert épointé et le groupe multiplicatif) est la même que sur un disque ouvert (resp. une couronne) comme on le verra ci-dessous, et c'est cette théorie qui nous intéresse ici. 

\bigskip
\trois{squelette} Soit $X$ une $k$-couronne. L'espace topologique sous-jacent est un arbre réel avec deux bouts. Si $x$ et $y$ sont deux points de $X$ on dit que $x\leq y$ si pour toute fonction $f$ analytique sur $X$ on a l'inégalité $|f(x)|\leq |f(y)|$. L'ensemble $S(X)$ des points maximaux pour cette relation d'ordre est l'unique intervalle réel ouvert qui joint les deux bouts de $X$ et est appelé le {\em squelette} de $X$. Pour tout $x$ appartenant à $X$ l'ensemble des majorants de $x$ possède un plus grand élément noté $x_{1}$ qui est donc sur $S(X)$. Pour tout $t$ compris entre $0$ et $1$ l'ensemble des points $y$ tels que pour toute fonction $f$ sur $X$ l'on ait l'encadrement $$|f(x)\leq |f(y)|\leq (1-t)|f(x)|+t|f(x_{1})|$$ possède un plus grand élément $x_{t}$ (la notation est compatible avec la précédente lorsque $t$ vaut 1, et $x_{0}$ est égal à $x$). L'application $(x,t)\mapsto x_{t}$ définit une rétraction de $X$ sur $S(X)$. Par sa construction même cette rétraction commute à tout automorphisme de $X$ même s'il n'induit pas l'identité sur $k$. 

\bigskip
Identifions $X$ à un domaine de la droite affine défini par une condition de la forme $r<|t|<R$. Le squelette $S(X)$ correspond alors à l'ensemble $\{\eta_{s}\}_{r<s<R}$ où pour tout réel positif $s$ on note $\eta_{s}$ le point défini par la semi-norme $\sum a_{i}t_{i}\mapsto \max |a_{i}|s^{i}$.

\bigskip
Notons enfin qu'il n'y a, modulo les fonctions de norme constante, que deux fonctions coordonnées sur $X$. En choisir une équivaut à orienter $S(X)$ (dans le sens qui la rend croissante en norme). Si l'on fait un tel choix on dira que l'on a {\em orienté $X$} et on dira de toute fonction coordonnée croissante en norme sur $S(X)$ qu'elle {\em représente} l'orientation fixée. 

\bigskip
\deux{defpseudo} Soit $X$ un bon espace $k$-analytique connexe, non vide et réduit. On dira que $X$ est un {\em pseudo-disque} s'il existe une extension finie $L$ de $\ra{k}$ telle que $X\times_{k}L$ soit une somme disjointe de $L$-disques. Il revient au même de demander que $\cx$ soit séparable sur $k$ et que $X\times_{\cx}L$ soit un $L$-disque pour une certaine extension finie $L$ de $\ra{\cx}$.  On dira que $X$ est une {\em pseudo-couronne} si topologiquement $X$ est un arbre réel à deux bouts et s'il existe une extension finie $L$ de $\ra{k}$ telle que $X\times_{k}L$ soit une somme disjointe de $L$-couronnes. Cette dernière condition équivaut à demander que $\cx$ soit séparable sur $k$ et que $X\times_{\cx}L$ soit une $L$-couronne pour une certaine extension finie $L$ de $\ra{\cx}$. Si $X$ est un pseudo-disque (resp. une pseudo-couronne) et si $L$ est une extension de $\cx$ on dira que $L$ {\em déploie} $X$ si $X\times_{\cx}L$ est un disque (resp. une couronne). 

\bigskip
Soit $X$ un pseudo-disque (resp. une pseudo-couronne). On appellera {\em sous-pseudo-couronne} de $X$ tout ouvert $Y$ de $X$ tel que $Y\times_{\cx}L$ soit une sous-couronne de $X\times_{\cx}L$ pour une (et donc pour toute) extension finie $L$ de $\cx$ déployant $X$. Notons que dans ce cas $\got{c}(Y)$ est égal à $\cx$. 

\bigskip
Si $X$ est un pseudo-disque (resp. une pseudo-couronne) et que $Z$ en est une sous-pseudo-couronne alors pour tout corps valué complet $L$ au-dessus de $\cx$ l'espace $X\times_{\cx}L$ est un $L$-pseudo-disque (resp.  une $L$-pseudo-couronne) et $Z\times_{\cx}L$ en est une sous-pseudo-couronne. 

\bigskip
\deux{purinsep} {\bf Exemple de pseudo-disque géométriquement connexe non trivial.} Soit $X$ la fibre générique d'un $k\zero$-schéma formel lisse $\got{X}$ purement de dimension 1 et soit $\bf x$ un point fermé de ${\got X}_{s}$ dont le corps résiduel est purement inséparable et de degré au moins égal à $p$ sur $\red{k}$. Soit $Z$ l'image réciproque de $\bf x$ dans $X$ par l'application de réduction. Alors $Z$ ne possède aucun point $k$-rationnel (puisque l'image d'un tel point sur ${\got X}_{s}$ est nécessairement $\red{k}$-rationnelle). Il existe une extension finie séparable $L$ de $k$ telle que $\red{L}$ s'identifie au corps résiduel de $\bf x$. L'image réciproque de $\bf x$ sur ${\got X}_{s}\times_{\red{k}}\red{L}$ est un singleton formé d'un point $\red{L}$-rationnel ; on en déduit que $Z\times_{k}L$ est un $L$-disque. 

\bigskip
\deux{pseudocour} Soit $X$ une pseudo-couronne. On notera $S(X)$ l'unique intervalle ouvert qui joint les deux bouts de $X$ et on l'appelle le {\em squelette} de $X$. Une {\em orientation} sur $X$ sera une orientation sur $S(X)$. Soit $L$ une extension finie galoisienne de $\ra{\cx}$ déployant $X$ ; notons $Y$ la couronne $X\times_{\cx}L$. Désignons par $G$ le groupe de Galois de $L$ sur $\ra{\cx}$. En vertu des remarques faites en~\ref{PSEU}.\ref{courdisque}.\ref{squelette} l'action de $G$ sur $Y$ commute à la rétraction de $Y$ sur $S(Y)$ et stabilise en particulier ce dernier. Soit $t$ une fonction coordonnée sur $Y$ et soit $g$ appartenant à $G$. L'image de $t$ par $g$ est de la forme $\lambda t^{\varepsilon}(1+u)$ où $\lambda$ appartient à $L^{*}$, où $\epsilon$ vaut $1$ ou $-1$ et où $u$ est une fonction $L$-analytique de valeur absolue partout strictement inférieure à 1 sur $Y$. Si $\epsilon$ valait $-1$ alors $g$ échangerait les deux bouts de $Y$ et donc $X$ n'aurait qu'un seul bout, ce qui est exclu par définition d'une pseudo-couronne. Un calcul immédiat montre alors que $g\mapsto \lambda$ est un cocycle, donc par le théorème de Hilbert 90 est de la forme $g\mapsto g(\mu)/\mu$ pour un certain élément $\mu$ fixé de $L^{*}$. En particulier $\lambda$ est toujours de valeur absolue égale à 1, d'où il découle que $G$ agit trivialement sur l'espace topologique $S(Y)$. La projection $Y\to X$ induit en conséquence un homéomorphisme entre $S(Y)$ et $S(X)$. Pour tout point $x$ de $S(X)$ le corps $\cx$ est algébriquement clos dans $\hres(x)$ et si $y$ désigne son unique antécédent sur $Y$ alors $\hres(y)$ s'identifie à l'extension $L\otimes_{\cx}\hres(x)$.  

\bigskip
Soit $L$ un corps valué complet au-dessus de $\cx$. De ce qui précède on déduit que $X\times_{\cx}L\to X$ induit un homéomorphisme entre $S(X\times_{\cx}L)$ et $S(X)$. En particulier toute orientation de $X$ en définit une sur $X\times_{\cx}L$.

\bigskip
\deux{pseudodisc} Soit $X$ un pseudo-disque. Soit $L$ une extension finie galoisienne de $\ra{\cx}$ telle que $X\times_{\cx}L$ soit un disque, que l'on notera $Y$. Désignons par $G$ le groupe de Galois de $L$ sur $k$ et identifions $Y$ (par le choix d'une fonction coordonnée) à un domaine de la droite affine sur $L$ défini par une inégalité de la forme $|t|<R$. Soit $g$ un élément de $G$. On peut écrire $g(t)$ sous la forme $\alpha +\lambda t +t^{2}\phi(t)$ où $\alpha$ est un scalaire de valeur absolue strictement inférieure à $R$ et où $\phi$ est une fonction analytique telle que $|t^{2}\phi(t)|<|\lambda||t|$ en tout point de $Y$. On montre comme ci-dessus que $|\lambda|$ est égale à 1. On en déduit qu'il existe une sous-couronne $Z$ de $Y$ stabilisée par $G$. D'après le~\ref{PSEU}.\ref{pseudocour} le groupe $G$ agit trivialement sur l'espace topologique $S(Z)$. L'image de $Z$ sur $X$ en est une sous-pseudo-couronne dont  l'image de $S(Z)$ est le squelette. L'un des deux bouts de cet intervalle est l'unique bout de $X$. 

\subsection*{La cohomologie des pseudo-disques et pseudo-couronnes} 

\deux{geomod}  Soit $X$ un espace $k$-analytique. Un revêtement fini étale $Y \to X$ sera dit {\em géométriquement modéré} si toute composante connexe de $Y\times_{k}\ka$ s'identifie à un quotient d'un revêtement fini galoisien d'ordre premier à $p$ d'une composante connexe de $X\times_{k}\ka$. Des exemples de ce type de revêtement sont fournis par les revêtements galoisiens d'ordre premier à $p$ ou, si les domaines affinoïdes de $X$ sont réduits, par les revêtements provenant d'une $\cx$-algèbre étale. Si $X$ est connexe alors tout point géométrique $x$ de $X$ fournit un groupe profini $\pi_{1}\gmod(X,x)$ et une équivalence naturelle entre la catégorie des revêtements étales géométriquement modérés de $X$ et celle des $\pi_{1}\gmod(X,x)$-ensembles discrets. On dira que $\pi_{1}\gmod(X,x)$ est le {\em groupe fondamental géométriquement modéré} de l'espace pointé $(X,x)$.

\bigskip
Si $X$ est géométriquement connexe et si $x$ est un point géométrique qui se factorise par $X\times_{k}\ka$ on dispose d'une suite exacte $$1\to \pi_{1}\gmod (X\times_{k}\ka,x) \to\pi_{1}\gmod(X,x)\to \mbox{Gal}\;(k^{a}/k)\to 1.$$ 

\bigskip
\deux{geomodka} Soit $Y$ un espace $\ka$-analytique et soit $y$ un point géométrique de $Y$. Berkovich a établi les deux résultats suivants : 

\bigskip
\begin{itemize}
\itb Si $Y$ est un disque alors $\pi_{1}\gmod(Y,y)$ est trivial. 
\itb Si $Y$ est une couronne alors la donnée d'une fonction coordonnée sur $Y$ induit un isomorphisme entre $\pi_{1}{\gmod}(Y,y)$  et $\zn$. 
\end{itemize}

\bigskip
Plaçons-nous dans le cas de la couronne et donnons-nous une fonction coordonnée $t$ sur $Y$. Pour tout entier $n$ premier à $p$ l'unique quotient de $\zn$ d'ordre $n$ correspond au revêtement fourni par l'extraction d'une racine $n$-ième de $t$.  On en déduit que l'isomorphisme $\pi_{1}{\gmod}(Y,y)\simeq \zn$ ne dépend que de la classe de $t$ modulo les fonctions de norme constante ; il est donc bien déterminé par le choix d'une orientation sur $Y$. Si l'on change l'orientation en son opposée (ce qui revient à remplacer $t$ par $t^{-1}$) l'isomorphisme est composé avec la multiplication par $(-1)$. Notons enfin que comme $\zn$ est abélien l'on dispose d'un isomorphisme {\em canonique} $\pi_{1}\gmod(Y,y)\simeq \pi_{1}\gmod(Y,z)$ pour tout couple $(y,z)$ de points géométriques de $Y$. 

\bigskip
\deux{souscouronne} {\bf Remarque.} Soit $Y$ une $\ka$-couronne, soit $Z$ une sous-couronne de $Y$ et soit $z$ un point géométrique de $Z$. La flèche naturelle $\pi_{1}\gmod(Z,z)\to \pi_{1}\gmod(Y,z)$ est alors un isomorphisme. 

\bigskip
\deux{indep} Soit $Y$ un espace $\ka$-analytique qui est un disque ou bien une couronne. Le groupe $\pi_{1}\gmod(Y,y)$ est indépendant (à isomorphisme unique près) du point géométrique $y$ de $Y$. On le note simplement $\pi$. Soit $\cal F$ un faisceau sur $Y$ provenant d'un $\pi$-module fini { \em de cardinal premier à $p$} que l'on notera encore $\cal F$. 

\bigskip
\deux{cohogalois} {\bf Lemme.} {\em Avec les hypothèses et notations ci-dessus la flèche naturelle $\H^{i}(\pi,{\cal F})\to \H^{i}(Y,{\cal F})$ est un isomorphisme pour tout entier $i$. En particulier le groupes $\H^{i}(Y,{\cal F})$ est fini pour tout $i$, et nul si $Y$ est un disque et si $i$ est strictement positif.} 

\bigskip
{\em Démonstration.}  C'est clair si $i$ est égal à $0$. Supposons que $i$ soit égal à $1$ et soit $h$ une classe appartenant à $\H^{1}(Y,{\cal F})$. Comme $\cal F$ est géométriquement modéré il existe un revêtement géométriquement modéré $Z\to Y$ tel que ${\cal F}_{|Z}$ soit isomorphe au faisceau constant associé à un certain groupe fini $G$ d'ordre premier à $p$. La classe $h_{|Z}$ est celle d'un revêtement $T\to Z$ de groupe $G$, qui est donc géométriquement modéré. Or $T$ trivialise $h$, ce qui permet de conclure. 

\bigskip
Si $i$ est au moins égal à 3 le lemme est vrai puisque les deux groupes sont nuls pour des raisons de dimension cohomologique. Dans le cas où $i$ est égal à 2 et où $\cal F$ est constant c'est vrai aussi puisque les deux groupes sont nuls, le premier encore pour des raisons de dimension cohomologique et le second d'après Berkovich. Dans le cas général il existe un revêtement fini galoisien $Z\to Y$ connexe et d'ordre premier à $p$ tel que la restriction de $\cal F$ à $Z$ soit constante. La comparaison des suites spectrales de Hochschild-Serre relatives à la cohomologie des groupes et à la cohomologie étale permet de conclure.~$\Box$ 

\bigskip
\deux{pigen} Soit $X$ une courbe $k$-analytique qui est ou bien un pseudo-disque ou bien une pseudo-couronne. Fixons un $k$-plongement $\cx\hookrightarrow k^{a}$ et notons $G$ le groupe $\mbox{Gal}\;(k^{a}/\cx)$. Soit $x$ un point géométrique de $X$ se factorisant par $X\times_{\cx}\ka$. Le groupe $\pi_{1}\gmod(X\times_{\cx}\ka,x)$ ne dépend pas (à isomorphisme unique près) du choix de $x$ ; notons-le $\pi$. On dispose d'une suite exacte $$1\to \pi \to \pi_{1}\gmod(X,x)\to G\to 1.$$ 
Notons que $\pi_{1}\gmod(X,x) $ est lui aussi indépendant du choix de $x$ à un isomorphisme unique près ({\em pourvu que $x$ soit astreint à se factoriser par le plongement fixé $\cx\hookrightarrow k^{a}$}). 
Dans le cas où $X$ est un pseudo-disque $\pi$ est trivial et l'exactitude de la suite signifie simplement que $\pi_{1}\gmod(X,x)\to G$ est un isomorphisme. Dans le cas où $X$ est une pseudo-couronne orientons-la (on oriente du même coup $X\times_{\cx}\ka$). Ce choix fournit un isomorphisme $\pi\simeq \zn$ et la suite exacte ci-dessus induit une action de $G$ sur $\pi$ (une vraie action puisque $\pi$ est abélien) et donc sur $\zn$. 
Un calcul simple montre que ce dernier est alors isomorphe, comme groupe profini muni d'une action de $G$, à $\znt$. 

\bigskip
\deux{cohogen} Conservons les notations $\pi$ et $G$ introduites ci-dessus ; notons $\Pi$ le groupe $\pi_{1}\gmod(X,x) $. Soit $\cal F$ un faisceau sur $X$ provenant d'un $\Pi$-module fini de cardinal premier à $p$ que l'on note encore $\cal F$. Si $X$ est un pseudo-disque alors $\Pi$ s'identifie à $G$ et $\cal F$ peut être vu comme un $\cx$-faisceau étale.

\bigskip
\deux{etalegalois} {\bf Lemme.} {\em Avec les hypothèses et notations ci-dessus la flèche naturelle $\H^{i}(\Pi,{\cal F})\to \H^{i}(X,{\cal F})$ est un isomorphisme pour tout entier $i$. En particulier si $X$ est un pseudo-disque (auquel cas $\cal F$ est un $\cx$-faisceau étale) alors $\H^{i}(\cx,{\cal F})\to \H^{i}(X,{\cal F})$ est un isomorphisme pour tout entier $i$.} 

\bigskip
{\em Démonstration.} On dispose de deux suites spectrales $$\H^{p}(G,\H^{q}(\pi,M))\Rightarrow \H^{p+q}(\Pi,M)$$ et $$\H^{p}(G,\H^{q}(X\times_{\cx}\ka,{\cal F})\;)\Rightarrow \H^{p+q}(X,{\cal F}).$$ Pour tout entier $q$ la flèche $\H^{q}(\pi,M)\to \H^{q}(X\times_{\cx}\ka,{\cal F})$  est un isomorphisme en vertu du lemme~\ref{PSEU}.\ref{cohogalois}. Le résultat cherché vient aussitôt.~$\Box$ 

\bigskip
\deux{remsouspseu} {\bf Remarque.} Si $X$ est une pseudo-couronne et si $Z$ en est une sous-pseudo-couronne il résulte de ce qui précède que pour tout point géométrique $z$ de $Z$ le morphisme naturel $\pi_{1}\gmod(Z,z)\to \pi_{1}\gmod(X,z)$ est un isomorphisme. Pour tout faisceau $\cal F$ sur $X$ provenant d'un $\pi_{1}\gmod(X,z)$-module fini de torsion première à $p$ et pour tout entier $q$ la flèche naturelle $\H^{q}(X,{\cal F})\to \H^{q}(Z,{\cal F})$ est un isomorphisme.

\section{Les systèmes de composantes et les valuations associées}\label{GERM} 

\setcounter{cpt}{0}
\setcounter{cptbis}{0}
\deux{introql} On va rappeler et étendre certains résultats de \cite{duc} (prop. 1.11, prop. 1.23....). Ils y sont énoncés dans le cas d'une courbe lisse, ils vont l'être ci-dessous {\em mutatis mutandis} dans le cas d'une courbe quasi-lisse, c'est-à-dire dont tout point a un voisinage isomorphe à un domaine affinoïde d'une courbe lisse. Le passage de lisse à quasi-lisse ne sera pas détaillé ; il repose sur le fait qu'un domaine affinoïde est toujours réunion de domaines affinoïdes rationnels (\cite{duc2}, lemme 2.4) et sur l'étude locale de la variation des normes de fonctions sur une courbe lisse (\cite{duc},  (1.19), (1.20) et prop. 1.21). 

\bigskip
\deux{introgerme} Soit $X$ une courbe $k$-analytique quasi-lisse et soit $P$ un point de $X$. On appellera {\em système de composantes adhérent à $P$} tout élément de $\lim\limits_{\leftarrow} \pi_{0}(V-\{P\})$ où $V$ parcourt l'ensemble des voisinages ouverts de $P$. Si ${\cal E}$ est un tel système et si $V$ est un ouvert contenant $P$ on notera ${\cal E}_{V}$ la composante connexe de $V-\{P\}$ définie par $\cal E$. 

\bigskip
\deux{normeconstante} {\bf Lemme.} {\em Soit $P$ un point de $X$ et soit $f$ une fonction méromorphe au voisinage de $P$. Soit ${\cal E}$ un système de composantes adhérent à  $P$. Il existe un voisinage ouvert $V$ de $P$ tel que $f$ soit analytique sur ${\cal E}_{V}$ et tel que la valeur absolue de $|f|$ soit ou bien partout strictement inférieure à $1$, ou bien partout égale à $1$, ou bien partout strictement supérieure à $1$ sur ${\cal E}_{V}$. Pour tous les systèmes de composantes adhérents à $P$ à l'exception d'un nombre fini d'entre eux c'est le second terme de l'alternative qui est vrai.} 

\bigskip
{\em Démonstration.} C'est établi dans le cas lisse dans  \cite{duc} ( proposition 1.11), lorsque $f$ est analytique et inversible au voisinage de $P$. Cela s'étend sans problèmes au cas quasi-lisse. Si l'anneau ${\cal O}_{X,P}$ est un corps la fonction $f$ est ou bien identiquement nulle ou bien définie et inversible au voisinage de $P$, ce qui permet de conclure d'après ce qui précède. Sinon  $P$ est un point rigide et si $f$ a un zéro ou un pôle en $P$ alors $|f|$ tend vers zéro ou vers l'infini à l'approche de $P$. Le résultat souhaité s'en déduit aussitôt.~$\Box$ 

\bigskip
\deux{vallocal} Si $P$ est un point de $X$ notons ${\cal M}_{P}$ la fibre en $P$ du faisceau des fonctions méromorphes sur $X$. Notons que si $X$ n'est pas rigide alors ${\cal M}_{P}$ est égal à ${\cal O}_{X,P}$, et possède une valeur absolue canonique induite par l'évaluation en $P$. Soit ${\cal E}$ un système de composantes adhérent à $P$. Le lemme ci-dessus permet d'associer à ${\cal E}$ une valuation sur ${\cal M}_{P}$ que l'on notera $v_{\cal E}$. Si $P$ est un point de type $(1)$ ou $(4)$ il n'y a qu'un seul système de composantes ${\cal E}$ adhérent à $P$ ; si de plus $P$ est rigide alors $v_{\cal E}$ est composée de la valuation discrète de ${\cal M}_{P}$ avec la valeur absolue de $\hres(P)$, sinon toute inégalité entre $|f|$ et un scalaire valable en $P$ l'est encore au voisinage de $P$ et donc $v_{\cal E}$ est simplement la valeur absolue canonique de ${\cal M}_{P}$. Si $P$ est un point de type $(3)$ il y a au plus deux systèmes de composantes adhérents à $P$ (exactement 2 si $P$ n'est pas situé sur le bord de la courbe). Soit ${\cal E}$ un tel système ; alors $v_{\cal E}$ est simplement la valeur absolue canonique de ${\cal M}_{P}$ (\cite{duc}, (1.17)). Soit $P$ un point de type $(2)$ et soit ${\cal C}$ l'unique courbe projective, normale, et géométriquement intègre sur la fermeture algébrique de $\red{k}$ dans $\red{\hres(P)}$ dont le corps des fonctions est $\red{\hres(P)}$. L'ensemble des systèmes de composantes adhérents à $P$ est en bijection avec l'ensemble des points fermés d'un ouvert non vide de $\cal C$ (égal à $\cal C$ si $P$ est un point intérieur de $X$). Si ${\cal E}$ est un système de composantes la valuation $v_{\cal E}$ est la composée de la valeur absolue canonique de ${\cal M}_{P}$ et de la valuation discrète de $\red{\hres(P)}$ associée au point fermé de $\cal C$ correspondant (\cite{duc}, 1.18). 

\bigskip
Notons que d'après ce qui précède si $P$ n'est pas de type $(3)$ alors $v_{\cal E}$ détermine $\cal E$. Pour tout point $P$ de $X$ et tout système de composantes $\cal E$ adhérent à $P$ on désigne par ${\cal O}({\cal E})$ l'anneau $\lim\limits_{\to} {\cal O}({\cal E}_{V})$ où $V$ parcourt l'ensemble des voisinages ouverts de $P$. 

\bigskip
\deux{henslocal} {\bf Proposition.} {\em Soit $P$ un point de $X$ et soit ${\cal E}$ un système de composantes adhérent à $P$. Soit $f$ un élément de ${\cal O}({\cal E})$ entier sur ${\cal M}_{P}$. Alors les conclusions du lemme~{\rm \ref{GERM}.\ref{normeconstante}} sont encore valables pour $f$ et permettent de munir la fermeture intégrale ${\cal O}({\cal E})_{\rm alg}$ de ${\cal M}_{P}$ dans ${\cal O}({\cal E})$ d'une valuation prolongeant $v_{\cal E}$ qui identifie le corps ${\cal O}({\cal E})_{\rm alg}$ au hensélisé de ${\cal M}_{P}$ pour $v_{\cal E}$.}

\bigskip
{\em Démonstration.} On procède en plusieurs étapes. 

\setcounter{cptbis}{0}

\bigskip
\trois{oealgsep} Commençons par montrer que ${\cal O}({\cal E})_{\rm alg}$ est séparable sur ${\cal M}_{P}$. On suppose que $p$ est un nombre premier et que $k$ lui-même est de caractéristique $p$. Soit $f$ appartenant à ${\cal M}_{P}$ s'écrivant $g^{p}$ pour un certain $g$ appartenant à ${\cal O}({\cal E})_{\rm alg}$. Le but est de prouver que $g$ appartient à ${\cal M}_{P}$. On peut toujours supposer (quitte à multiplier $f$ par une puissance convenable d'une uniformisante dans le cas où $P$ est rigide) que $f$ appartient à ${\cal O}_{X,P}$, et même à ${\cal O}_{X}(X)$ en remplaçant $X$ par un voisinage convenable de $P$. Soit $Q$ l'unique antécédent de $P$ sur le revêtement fini et plat $Y$ de $X$ obtenu par l'extraction d'une racine $p$-ième de $f$. L'anneau local ${\cal O}_{Y,Q}$ s'identifie à ${\cal O}_{X,P}[T]/(T^{p}-f)$. Il existe un voisinage ouvert $V$ de $P$ tel que $f$ soit une puissance $p$-ième sur $V-\{P\}$. En conséquence si $W$ désigne l'image réciproque de $V$ sur $Y$ l'espace $W-\{Q\}$ n'est réduit en aucun de ses points. Le lieu des points en lequel un bon espace de Berkovich est réduit est un ouvert. Dès lors $Y$ n'est pas réduit en $Q$ et donc ${\cal O}_{X,P}[T]/(T^{p}-f)$ n'est pas réduit. Or cet anneau s'injecte dans ${\cal M}_{P}[T]/(T^{p}-f)$ puisque $T^{p}-f$ est unitaire. On en déduit que ${\cal M}_{P}[T]/(T^{p}-f)$ n'est pas réduit, et donc que $f$ est la puissance $p$-ième d'un élément de ${\cal M}_{P}$. Autrement dit $g$ appartient à ${\cal M}_ {P}$. 

\bigskip
\trois{oealgrig} Si $P$ est un point rigide $\cal E$ est nécessairement l'unique système de composantes adhérent à $P$ et la valuation correspondante est hensélienne ; il suffit dès lors de montrer que ${\cal M}_{P}$ est intégralement clos dans ${\cal O}({\cal E})$. Soit $\tau$ une uniformisante de ${\cal O}_{X,P}$. Si $\hres(P)$ est séparable sur $k$ alors $P$ possède par quasi-lissité de $X$ un voisinage dans $X$ isomorphe à un $\hres(P)$-disque. Tout élément de ${\cal O}({\cal E})$ annulant un polynôme à coeffcients dans ${\cal M}_{P}$ est majoré par $|\tau|^{m}$ pour un certain entier relatif $m$ au voisinage épointé de $P$. La forme explicite des fonctions au voisinage épointé du centre d'un disque permet aussitôt d'en déduire que $f$ est méromorphe. Dans le cas général il existe une extension finie purement inséparable $L$ de $k$ telle que si $Q$ désigne l'unique antécédent de $P$ sur $X\times_{k}L$ alors $\hres(Q)$ soit séparable sur $L$. On en déduit que $f$ est méromorphe sur $X_{L}$ au voisinage de $Q$ puis que $f^{p^{n}}$ appartient à ${\cal M}_{P}$ pour un certain entier $n$. Il découle du~\ref{GERM}.\ref{henslocal}.\ref{oealgsep} ci-dessus que $f$ appartient à ${\cal M}_{P}$.

\bigskip
\trois{oealgtype3} Si $P$ est de type $(3)$ la valuation $v_{\cal E}$ est la valeur absolue canonique sur ${\cal M}_{P}$, pour laquelle ce dernier est hensélien. Soit $\cal R$ un polynôme unitaire irréductible et séparable à coefficients dans ${\cal M}_{P}$ annulant une fonction $f$ de ${\cal O}({\cal E})$. Quitte à restreindre $X$ on peut supposer que les coefficients de $\cal R$ sont définis sur $X$, et que $f$ est définie sur ${\cal E}_{X}$. Soit $Y$ le revêtement fini et plat de $X$ donné par l'adjonction d'une racine de $\cal R$. Alors $P$ n'a qu'un antécédent $Q$ sur $Y$, et il y a un et un seul système de composantes sur $Y$ au-dessus de ${\cal E}$ (\cite{duc}, (1.17)). Comme $P$ annule une fonction analytique sur ${\cal E}_{X}$ il est nécessairement de degré 1 : dans le cas contraire le polynôme quotient de $\cal R$ par $T-f$ sur ${\cal O}({\cal E}_{X})$ définirait au moins un second système de composantes au-dessus de $\cal E$. Ceci montre que ${\cal M}_{P}$ est algébriquement clos dans ${\cal O}({\cal E})$ et permet de conclure. 

\bigskip
\trois{henselcour} Supposons maintenant que $P$ n'est pas rigide ni de type $(3)$. Dans ce cas ${\cal O}_{X,P}$ est égal au corps ${\cal M}_{P}$ et un système de composantes adhérent en $P$ est déterminé par la valuation correspondante. Fixons un système projectif $(X_{i},P_{i})_{i}$ d'espaces analytiques pointés au-dessus de $(X,P)$ tel que $X_{i}$ soit pour tout $i$ un revêtement étale galoisien d'un voisinage ouvert de $P$ dans $X$ dont la fibre en $P$ est (ensemblistement) réduite à $\{P_{i}\}$ et tel que la limite inductive des $\hres(P_{i})$ constitue une clôture séparable de $\hres(P)$. La limite inductive ${\cal M}$ des ${\cal M}_{P_{i}}$ est alors une clôture séparable de ${\cal M}_{P}$. 

\bigskip
{\em Dans ce  qui suit on utilise implicitement, à plusieurs reprises, le fait qu'un morphisme fini et plat est ouvert et fermé.} Choisissons une famille compatible $({\cal E}_{i})_{i}$ où ${\cal E}_{i}$ est pour tout $i$ un système de composantes adhérent à $P_{i}$ et situé au-dessus de ${\cal E}$. Soit $f$ appartenant à ${\cal O}({\cal E)}_{\rm alg}$. Il annule un polynôme unitaire séparable $\cal R$ à coefficients dans ${\cal M}_{P}$. Il existe $i$ et un voisinage ouvert $W$ de $P_{i}$ dans $X_{i}$, que l'on peut supposer irréductible et dont on notera $V$ l'image dans $X$, tels que la fonction $f$ soit définie sur ${\cal E}_{V}$ et tels que l'on puisse écrire ${\cal R}=\prod(T-f_{j})$ où les $f_ {j}$ sont des fonctions analytique sur $W$. Comme ${\cal E}_{i,W}$ est connexe et réduit (par lissité), et comme ${\cal R}(f)$ est nul on a l'égalité $f=f_{j}$ pour l'un des $j$ sur ${\cal E}_{i,W}$, et donc sur $W$ puisque ${\cal E}_{i,W}$ en est un ouvert non vide et que $W$ est réduit et irréductible. On en déduit que les conclusions du lemme~\ref{GERM}.\ref{normeconstante} (que l'on applique à $f_{j}$ sur l'espace $X_{i}$)  sont satisfaites par $f$. 

\bigskip
Soit $v$ la valuation sur $\cal M$ définie par la famille des $v_{{\cal E}_{i}}$ ; elle prolonge $v_{\cal E}$ et induit la valuation souhaitée sur ${\cal O}({\cal E}_{\rm alg}$. Le groupe de décomposition $G$ de $v$ dans ${\cal M}$ est exactement le sous-groupe de $\mbox{Gal}({\cal M}/{\cal M}_{P})$ qui stabilise $({\cal E}_{i})_{i}$. Un argument de descente galoisienne montre aussitôt qu'un élément $f$ de ${\cal M}$ ( qui provient nécessairement de l'un des ${\cal M}_{P_{i}}$) est invariant sous $G$ si et seulement il existe un voisinage ouvert $W$ de $P_{i}$ dans $X_{i}$ tel que la restriction de $f$ à ${\cal E}_{i,W}$ provienne de ${\cal O}_{X}({\cal E}_{V})$ où $V$ désigne l'image de $W$ sur $X$.  Ceci achève la démonstration.~$\Box$

\bigskip
\deux{systcompel} Soit $P$ un point de $X$. Un système de composantes ${\cal E}$ adhérent à $P$ sera dit {\em élémentaire} si l'on est dans l'une des trois situations suivantes :

\bigskip
\begin{itemize} 

\itb Le point $P$ est $k$-rationnel. 

\bigskip
\itb Le point $P$ est de type $(2)$ et il existe une $k\zero$-courbe intègre, propre et lisse $\cal X$, un voisinage $V$ de $P$ dans $X$ et un isomorphisme entre $V$ et un domaine $k$-analytique de ${\cal X}_{\eta}\an$ {\em via} lequel $P$ correspond à l'unique antécédent du point générique de ${\cal X}_{s}$ et ${\cal E}$ au système de composantes donné par $\pi^{-1}(\{x\})$ pour un certain $\red{k}$-point $x$ de ${\cal X}_{s}$, où $\pi$ est la flèche de réduction.

\bigskip
\itb Le point $P$ est de type $(3)$ et possède un voisinage dans $X$ isomorphe à une couronne {\em compacte} dont $P$ est élément du squelette.

\end{itemize}
\bigskip
\deux{systpseuco} Soit  $P$ un point de $X$ ou bien rigide, ou bien de type $(2)$, ou bien de type $(3)$. Soit ${\cal E}$ un système de composantes adhérent à $P$. Il existe alors un système cofinal $\cal S$ de voisinages $V$ de $P$ dans $X$ possédant la propriété suivante :  {\em pour tout $V$ appartenant à $\cal S$ l'espace ${\cal E}_{V}$ est une pseudo-couronne déployée par une extension finie de $k$ et si $V$ et $W$ sont deux éléments de $\cal S$ tels que $W\subset V$ alors ${\cal E}_{W}$ est une sous-pseudo-couronne de ${\cal E}_{V}$.}

\bigskip
Par ailleurs le corps ${\cal O}({\cal E})_{\rm alg}$ satisfait, en tant qu'extension du corps valué hensélien $k$, les hypothèses du paragraphe~\ref{HENS}.\ref{fixenotation} ; on le déduit des faits rappelés au~\ref{GERM}.\ref{vallocal}. Notons $\got{c}({\cal E})$ la fermeture algébrique de $k$ dans ${\cal O}({\cal E})$. Pour $V$ appartenant à $\cal S$ la flèche naturelle $\got{c}({\cal E}_{V})\to \got{c}({\cal E})$ est un isomorphisme. Fixons un plongement de $\got{c}({\cal E})$ dans $k^{s}$ et une clôture séparable ${\cal O}({\cal E})_{\rm alg}^{s}$ de ${\cal O}({\cal E})_{\rm alg}\otimes_{\got{c}({\cal E})}k^{s}$. Notons $G$ le groupe $\mbox{Gal}(k^{s}/\got{c}({\cal E}))$. 

\bigskip
\deux{inductetale} Pour tout $V$ appartenant à ${\cal S}$ désignons par $\Pi_{V}$ le groupe fondamental géométriquement modérés de ${\cal E}_{V}$ (qui est défini sans ambiguïté une fois choisi un plongement $\got{c}({\cal E})\hookrightarrow k^{a}$). Il résulte de la remarque \ref{PSEU}.\ref{remsouspseu} que si $V$ et $W$ sont deux éléments de $\cal S$ tels que $W\subset V$ alors la flèche $\Pi_{W}\to \Pi_{V}$ est un isomorphisme. Notons $\Pi_{\cal E}$ la limite projective des $\Pi_{V}$. Si $\cal F$ est un $\Pi_{\cal E}$-module on peut le voir comme un faisceau étale sur ${\cal E}_{V}$ pour tout $V$ appartenant à $\cal S$ ; s'il est fini et de torsion première à $p$ les flèches de transition $\H^{q}({\cal E}_{V},{\cal F})\to \H^{q}({\cal E}_{W},{\cal F})$ sont des isomorphismes pour tout couple $(V,W)$ d'éléments de $\cal S$ tels que $W\subset V$. La limite inductive des  $\H^{q}({\cal E}_{V},{\cal F})$ sera dans ce cas notée $\H^{q}({\cal E},{\cal F})$. 

\bigskip
\deux{foncteuretale} Soit $L$ une ${\cal O}({\cal E})_{\rm alg}$-algèbre finie étale. Elle définit, {\em via} la tensorisation avec ${\cal O}({\cal E})$, un revêtement fini étale de ${\cal E}_{V}$ pour $V$ suffisamment petit dans $\cal S$. 

\bigskip
\deux{equivetale} {\bf Lemme.} {\em On garde les hypothèses et notations introduites ci-dessus. Supposons que $L$ soit géométriquement modérément ramifiée et soit $V$ tel que $L$ soit définie sur ${\cal O}({\cal E}_{V})$. 

\bigskip

$i)$ Le revêtement de ${\cal E}_{V}$ induit par $L$ est géométriquement modéré et se prolonge en conséquence d'une unique manière en un revêtement de ${\cal E}_{W}$ pour tout ouvert $W$ de $\cal S$ contenant $V$. 

\bigskip
$ii)$ Cette construction induit un isomorphisme entre $\Pi_{\cal E}$ et  le quotient du groupe  de Galois de ${\cal O}({\cal E})_{\rm alg}^{s}$ sur ${\cal O}({\cal E})_{\rm alg}$ par son groupe de ramification. 

\bigskip
$iii)$ Soit $\cal F$ un $\Pi_{\cal E}$-module fini de torsion première à $p$. On peut le voir d'après le $ii)$ comme un faisceau étale sur ${\cal O}({\cal E})_{\rm alg}$ aussi bien que comme un faisceau étale sur chacun des ${\cal E}_{V}$. L'on dispose alors d'un système d'isomorphismes naturels $\H^{q}({\cal O}({\cal E})_{\rm alg},{\cal F })\simeq \H^{q}({\cal E},{\cal F})$. }

\bigskip
{\em Démonstration.} Si le système ${\cal E}$ est élémentaire les deux premières assertions sont établies grâce aux descriptions explicites des deux groupes concernés et des revêtements qu'ils décrivent (\cf. ~\ref{HENS}.\ref{cohomodereram} et~\ref{PSEU}.\ref{geomodka} {\em supra}) ; la troisième se déduit du fait que les groupes de cohomologie en jeu sont en fait les groupes de cohomologie du $\Pi_{\cal E}$-module discret $\cal F$ d'après le~\ref{HENS}.\ref{suitespec} et le~\ref{PSEU}.\ref{cohogalois}.

\bigskip
Dans le cas général il existe une extension finie normale $L$ de $k$, un point $Q$ sur $X\times_{k}L$ situé au-dessus de $P$, et un système de composantes ${\cal E}'$ adhérent à $Q$ qui est situé au-dessus de $\cal E$ et est élémentaire. On peut remplacer $k$ par sa fermeture radicielle dans $L$ (ce qui ne modifie pas le topos étale) et donc supposer $L$ galoisienne sur $k$. Soit $H$ le groupe de Galois de $\hres(Q)$ sur $\hres(P)$ et soit $H'$ le sous-groupe de $H$ qui fixe ${\cal E}'$. Pour tout voisinage $W$ suffisamment petit de $Q$ l'espace ${\cal E}'_{W}$ est galoisien de groupe $H'$ sur ${\cal E}_{V}$, où $V$ désigne l'image de $W$ sur $X$. Le revêtement ${\cal E}'_{W}\to {\cal E}_{V}$ est géométriquement modéré puisqu'il est induit par une extension de $k$. Il est  par ailleurs immédiat que ${\cal O}({\cal E})_{\rm alg}$ s'identifie au sous-corps de ${\cal O}({\cal E}')_{\rm alg}$ formé des invariants sous l'action de $H'$. En conséquence ${\cal E}'_{\rm alg}$ est de manière naturelle une extension galoisienne de ${\cal O}({\cal E})_{\rm alg}$ de groupe $H'$. On en déduit aussitôt, à l'aide du cas élémentaire, les assertions $i)$ et $ii)$. On utilise le~\ref{HENS}.\ref{suitespec} et le~\ref{PSEU}.\ref{etalegalois} pour établir la validité de $iii)$.~$\Box$

\bigskip
\deux{orient} {\bf Remarque.} Toute orientation de {\em l'une} des pseudo-couronnes ${\cal E}_{V}$ pour $V$ dans $\cal S$ en induit une sur {\em chacune} des ${\cal E}_{V}$. On dira alors que l'on a {\em orienté} $\cal E$. Il revient au même d'orienter ${\cal O}({\cal E}_{\rm alg})$ ; le lien entre les deux notions consiste à associer à un représentant $\xi$ de l'orientation choisie sur ${\cal O}({\cal E}_{\rm alg})$, représentant défini sur ${\cal E}_{V}$ pour un certain $V$ élément de $\cal S$, l'orientation de ${\cal E}_{V}$ qui rend $\xi$ croissante en norme sur le squelette. Un tel choix étant fait l'isomorphisme évoqué au point $ii)$ du lemme ci-dessus est un isomorphisme d'extensions de $G$ par $\znt$. 

\bigskip
\deux{constcour} Orientons ${\cal O}({\cal E})_{\rm alg}$ et donc $\cal E$. Donnons-nous un élément $\xi$ de ${\cal O}({\cal E})_{\rm alg}^{*}$ représentant l'orientation choisie. {\em Comme cette propriété ne dépend que de la valuation de $\xi$ ce dernier peut être choisi dans ${\cal M}_{P}$}. Soit $V$ un ouvert élément de $\cal S$ tel que $\xi$ soit définie sur ${\cal E}_{V}$. Soit $n$ un entier premier à $p$. Désignons par $\eta_{n}$ la classe de $\xi$ dans le groupe $\H^{1}({\cal O}({\cal E})_{\rm alg},\mu_{n})\simeq \H^{1}({\cal E},\mu_{n})$. Pour tout $W$ appartenant à $\cal S$ on note encore $\eta_{n}$ l'élément correspondant de $\H^{1}({\cal E}_{W},\mu_{n})$. 

\bigskip
Soit $\cal F$ un $G$-module fini annulé par $n$ et soit $W$ appartenant à $\cal S$. Soit $q$ un entier. En vertu de l'assertion $iii)$ du lemme~\ref{GERM}.\ref{equivetale} et de la proposition~\ref{HENS}.\ref{residucup} l'application $(\lambda,h)\mapsto \lambda +\eta_{n}\cup h$ induit un isomorphisme $$\H^{q}(\got{c}({\cal E}_{W}),{\cal F})\oplus \H^{q-1}(\got {c}({\cal E}_{W},{\cal F}(-1))\simeq \H^{q}({\cal E}_{W},{\cal F}).$$

\bigskip
\deux{courskel} On garde les hypothèses et notations du {\rm \ref{GERM}.\ref{systpseuco}, {\em en supposant de plus que $X$ est une pseudo-couronne, que $P$ est situé sur $S(X)$ et que $\cal E$ est le système de composantes adhérent à $P$ défini par l'une des deux composantes connexes de $X-\{P\}$}. On peut choisir le système cofinal ${\cal S}$ de sorte que pour tout $V$ appartenant à $\cal S$ l'espace ${\cal E}_{V}$ soit une sous-pseudo couronne de $X$. Soit $\cal F$ un $\Pi$-module fini annulé par un entier $n$ premier à $p$. Pour tout ouvert $V$ appartenant à $\cal S$ et pour tout entier $q$ la flèche $\H^{q}(X,{\cal F})\to \H^{q}({\cal E}_{V},{\cal F})$ est un isomorphisme. Ces deux groupes s'identifient d'après le lemme~\ref{GERM}.\ref{equivetale} à $\H^{q}({\cal O}({\cal E})_{\rm alg},{\cal F})$. Ces remarques ont, en vertu de ce qui précède, les conséquences suivantes : 

\setcounter{cptbis}{0}
\bigskip
\trois{locglob} {\bf  Proposition.} {\em Soit $X$ une pseudo-couronne. Fixonx un plongement de $\cx$ dans $k^{a}$ et notons $G$ le groupe de Galois correspondant. Soit $\Pi$ le groupe fondamental géométriquement modéré de $X$ et soit $\cal F$ un faisceau étale sur $X$ donné par un $\Pi$-module fini de cardinal premier à $p$. Soit $P$ un point de $S(X)$. Donnons-nous un entier $q$. Si $h$ est un élément de $\H^{q}(X,{\cal F})$ nul en $P$ alors $h$ est triviale.}

\bigskip
{\em Démonstration.} Comme $h$ est nulle en $P$ elle est nulle sur ${\cal E}_{V}$ pour un certain $V$, et donc est égale à zéro d'après ce qui précède.~$\Box$ 

\bigskip
\trois{exactescinde} {\bf Proposition.} {\em Les notations $X$ et $G$ sont celles de la proposition ci-dessus. Il existe une classe $(\eta_{n})_{n}$ appartenant à $\H^{1}(X,\znt)$ telle que pour tout $G$-module fini $\cal F$ de cardinal premier à $p$ et tout entier $q$ l'application  $(\lambda,h)\mapsto \lambda +\eta_{n}\cup h$ induise un isomorphisme $$\H^{q}(\got{c}(X),{\cal F})\oplus \H^{q-1}(\got {c}(X),{\cal F}(-1))\simeq \H^{q}(X,{\cal F}).\;\Box$$}

\bigskip
\deux{loczeropseudisc} Soit $X$ un pseudo-disque. Fixons un plongement de $\cx$ dans $k^{a}$ et notons $G$ le groupe de Galois correspondant. Soit ${\cal F}$ un $G$-module fini de cardinal premier à $p$, soit $q$ un entier et soit $h$ une classe de $\H^{p}(X,{\cal F})$ qui s'identifie à $H^{p}(\cx,{\cal F})$ d'après le lemme~\ref{PSEU}.\ref{etalegalois}. Donnons-nous un point $P$ situé sur le squelette d'une sous-pseudo-couronne $Z$ de $X$. Supposons que $h(P)$ est nulle. Alors la restriction de $h$ à $Z$ est nulle d'après la proposition~\ref{GERM}.\ref{courskel}.\ref{locglob} ; comme $\got{c}(Z)$ est égal à $\cx$ et comme $\H^{q}(\got{c}(Z),{\cal F})\to \H^{q}(Z,{\cal F})$ est injective d'après la proposition~\ref{GERM}.\ref{courskel}.\ref{exactescinde} ci-dessus la classe $h$ est elle-même nulle. 

\subsection*{À propos des images directes supérieures}

\bigskip
\deux{rpqcourbes} Soit $X$ une courbe $k$-analytique quelconque et $\cal F$ un faisceau en groupes abéliens sur $X$. La dimension topologique de $X$ étant inférieure ou égale à 1 les suites spectrales $$\H^{p}(|X|,\rqq{\cal F})\Rightarrow \H^{p+q}(X,{\cal F})\;\mbox{et}\;\H^{p}_{c}(|X|,\rqq{\cal F})\Rightarrow \H^{p+q}_{c}(X,{\cal F})$$ induisent pour tout entier $q$ deux suites exactes $$0\to\H^{1}(|X|,\mbox{\rm R}^{q-1}\pi_{*}{\cal F})\to \H^{q}(X,{\cal F})\to \H^{0}(|X|,\rqq{\cal F})\to 0$$ et  $$0\to\H^{1}_{c}(|X|,\mbox{\rm R}^{q-1}\pi_{*}{\cal F})\to \H^{q}_{c}(X,{\cal F})\to \H^{0}_{c}(|X|,\rqq{\cal F})\to 0.$$

\bigskip
\deux{compcech} {\bf Remarque.} La flèche $\H^{1}(|X|,\mbox{\rm R}^{q-1}\pi_{*}{\cal F})\to \H^{q}(X,{\cal F})$ peut se décrire comme suit. Soit $h$ une classe appartenant à $\H^{1}(|X|,\mbox{\rm R}^{q-1}\pi_{*}{\cal F})$. La dimension topologique de $X$ étant inférieure ou égale à $1$ cette classe provient d'un élément $\lambda$ de $\check{\H }^{1}({\cal U}/X,\underline{\H}^{q-1}({\cal F}))$ pour un certain recouvrement $\cal U$ de $X$ par des ouverts $U_{i}$ tels que $U_{i}\cap U_{j}\cap U_{r}$ soit vide pour tout triplet $(i,j,r)$ d'indices distincts, et où $\underline{\H}^{q-1}({\cal F})$ est le préfaisceau $U\mapsto \H^{q}(U,{\cal F})$. La suite spectrale $\check{\H}^{i}({\cal U}/X,\underline{\H}^{j}({\cal F})\Rightarrow \H^{i+j}(X,{\cal F})$ dégénère pour $i$ strictement plus grand que $1$ et fournit donc une flèche de $\check{\H }^{1}({\cal U}/X,\underline{\H}^{q-1}({\cal F}))$ vers $\H^{q}(X,{\cal F})$. L'image de $\lambda$ par la flèche en question est alors l'élément cherché : c'est immédiat lorsque $q$ vaut $0$ ou $1$, et on conclut par une récurrence facile en plongeant $\cal F$ dans un $\ZZ/n$-module injectif. {\em On dispose d'une description analogue de la flèche $\H^{1}_{c}(|X|,\mbox{\rm R}^{q-1}\pi_{*}{\cal F})\to \H^{q}_{c}(X,{\cal F})$.}

\bigskip
\deux{introaltercup} Donnons-nous maintenant deux faisceaux en groupes abéliens ${\cal F}$ et ${\cal G}$ sur $X$. Soient $q$ et $l$ deux entiers. On peut définir {\em via} le cup-produit une application naturelle {\small $$\H^{0}(|X|,\mbox{R}^{q}\pi_{*}{\cal F})\times \H^{1}_{c}(|X|,\mbox{R}^{l}\pi_{*}{\cal G})\to \H^{1}_{c}(|X|,\mbox{R}^{q+l}\pi_{*}{\cal F}\otimes {\cal G})\to  \H^{q+l+1}_{c}(X,{\cal F}\otimes {\cal G}).$$} Donnons-en une description alternative : soit $h$ appartenant à $\H^{0}(|X|,\mbox{R}^{q}\pi_{*}{\cal F})$. Elle possède un antécédent $\lambda$ dans $\H^{q}(X,{\cal F})$. Soit $\eta$ un élément du groupe  $\H^{1}_{c}(|X|,\mbox{R}^{l}\pi_{*}{\cal G})$ et soit $\mu$ son image dans $\H^{l+1}_{c}(X,{\cal G})$. Alors le cup-produit de $\lambda$ et $\mu$ ne dépend pas du choix de $\lambda$ et est précisément égal à l'image de $(h,\eta)$ par l'accouplement mentionné ci dessus ; on déduit ce fait de la remarque~\ref{GERM}.\ref{compcech} ci-dessus et de la description des cup-produits à partir des complexes de \v{C}ech (\cf. \cite{mil}, chap. V, rem. 1.19). {\em On dispose de résultats semblables concernant $\H^{0}_{c}(|X|,.)\times \H^{1}(|X|,.)$.}

\bigskip
\deux{aretetriv} {\bf Proposition.} {\em Soit $X$ une courbe $k$-analytique qui est ou bien un pseudo-disque, ou bien une pseudo-couronne. Soit $\Pi$ le groupe $\pi_{1}\gmod(X,x)$ où $x$ est un point géométrique se factorisant par un plongement $\cx \hookrightarrow \ka$ préalablement fixé, et soit $\cal F$ un $\Pi$-module fini de cardinal premier à $p$. Alors pour tout entier $q$ la flèche $\H^{q}(X,{\cal F})\to \H^{0}(|X|,\rqq{\cal F})$ est un isomorphisme et  le groupe $\H^{1}(|X|,\mbox{\rm R}^{q-1}\pi_{*}{\cal F})$ est trivial.} 

\bigskip
{\em Démonstration.} En vertu de la suite exacte mentionnée au~\ref{GERM}.\ref{rpqcourbes} il suffit de prouver que $\H^{q}(X,{\cal F})\to \H^{0}(|X|,\rqq{\cal F})$ est injective. Commençons par le cas où $X$ est une pseudo-couronne. Si $h$ est une classe de $\H^{q}(X,{\cal F})$ qui s'annule dans $\H^{0}(|X|,\rqq{\cal F})$ alors $h(P)$ est nulle pour tout point $P$ de $X$. C'est en particulier vrai en au moins un point de $S(X)$ et on déduit alors de la proposition~\ref{GERM}.\ref{courskel}.\ref{locglob} que $h$ est triviale. Supposons maintenant que $X$ est un pseudo-disque. Dans ce cas $\cal F$ est un $G$-module et $\H^{q}(X,{\cal F})$ s'identifie à $\H^{q}(\cx,{\cal F})$ d'après le lemme~\ref{PSEU}.\ref{etalegalois}. Comme  $X$ contient une sous-pseudo-couronne $Z$ d'après le~\ref{PSEU}.\ref{pseudodisc} le~\ref{GERM}.\ref{loczeropseudisc} entraîne la trivialité de $h$.~$\Box$ 

\section{Triangulation d'une courbe} \label{TRI}

\setcounter{cpt} {0}
\deux{deftri} {\bf Définition.} {\em Soit $X$ une $k$-courbe. On appelle {\em triangulation} de $X$ tout sous-ensemble fermé et discret $\bf S$ de $X$ formé de points ou bien rigides, ou bien de type $(2)$, ou bien de type $(3)$, qui rencontre toutes les composantes connexes de $X$ et dont le complémentaire est somme disjointe de pseudo-disques et de pseudo-couronnes. Les points de $\bf S$ seront dans ce cas appelés les {\em sommets} ; les composantes connexes de $X-{\bf S}$ seront les {\em arêtes}.}

\bigskip
\deux{extri} {\bf Exemples.} Si $X$ est la fibre générique d'un schéma formel pluristable $\got{X}$ alors l'ensemble $\bf S$ des images réciproques des points génériques des composantes irréductibles de $\got{X}_{s}$ constitue une triangulation. 

\bigskip
La donnée de n'importe quel $k$-point de $\PP^{1,an}_{k}$ constitue une triangulation avec une seule arête. 

\bigskip
Si $X$ est un pseudo-disque (resp. une pseudo-couronne) et si $P$ appartient au squelette d'une sous-pseudo-couronne de $X$ (resp. au squelette de $X$) alors $\{P\}$ est une triangulation de $X$. 

\bigskip
\deux{arsom} Soit $X$ une $k$-courbe et soit $\bf S$ une triangulation de $X$. Soit $Y$ une arête. Comme elle ne constitue pas une composante connexe de $X$ son adhérence rencontre $\bf S$. Les éléments du bord de $Y$ dans $X$ seront appelés les {\em sommets de $Y$.} Si $Y$ est un pseudo-disque c'est un arbre réel à un bout ; en conséquence $Y$ possède un et un seul sommet $S$ (rappelons que par définition une courbe est séparée) et $Y\cup \{S\}$ est compact et contractile. Si $Y$ est une pseudo-couronne c'est un arbre réel à deux bouts. On a dans ce cas trois possibilités : ou bien $Y$ a deux sommets distincts $S_{1}$ et $S_{2}$ chacun situé à l'un des deux bouts (auquel cas $Y\cup\{S_{1},S_{2}\}$ est compact et contractile) ou bien $Y$ a un sommet $S$ et $Y\cup\{S\}$ a encore un bout (et est contractile), ou bien $Y$ a un sommet $S$ et $Y\cup\{S\}$ est compact et se rétracte sur $S(Y)\cup\{S\}$ qui est alors un cercle. 

\bigskip
\deux{fincour} Soit $S$ un point de $\bf S$ et soit $V$ un voisinage $k$-affinoïde de $S$ dans $X$. Soit $U$ l'intérieur de $V$ dans $X$. Alors $V-U$ est formé d'un nombre fini de points. En particulier $U$ a un nombre fini de bouts. Soit $E$ l'ensemble des arêtes de $X$ aboutissant à $S$. Soit $Y$ appartenant à $E$. Un bout de $Y\cap U$ provient ou bien d'un bout de $Y$, ou bien d'un bout de $U$. Il existe donc un sous-ensemble fini $F$ de $E$ tel que pour toute arête $Y$ appartenant à $E-F$ les bouts de $Y\cap U$ soient exactement ceux de $Y$, ce qui implique que $Y$ est contenu dans $U$. Si $Y$ est une pseudo-couronne appartenant à $E-F$ et si $S(Y)\cup\{S\}$ n'est pas compact alors le bout de $S(Y)\cup\{S\}$ provient d'un bout de $U$. En conséquence il existe un sous-ensemble fini $F'$ de $E$ contenant $F$ tel que pour toute pseudo-couronne $Y$ appartenant à $E-F$ l'espace $S(Y)\cup\{S\}$ soit compact (et donc homéomorphe à un cercle sur lequel $Y\cup\{S\}$ se rétracte). Comme $U$ est homotopiquement équivalent à un arbre fini l'ensemble des arêtes appartenant à $E$ qui sont des pseudo-couronnes est fini. 

\bigskip
\deux{triexist} {\bf Proposition.} {\em Soit $X$ une $k$-courbe quasi-lisse. Soit $\bf T$ un sous-ensemble fermé et discret de $X$ formé de points de type $(2)$ ou $(3)$. Il existe une triangulation $\bf S$ de $X$ formée uniquement de points de type $(2)$ ou $(3)$ et contenant $\bf T$.}

\bigskip
{\em Démonstration.}  On procède en plusieurs étapes. {\em Le caractère ouvert et fermé des morphismes finis et plats sera implicitement utilisé tout au long de ce qui suit.} 
\setcounter{cptbis}{0}

\bigskip
\trois{triangloc} Soit $P$ un point de $X$. Comme $X$ est quasi-lisse il existe un voisinage de $P$ dans $X$, que l'on peut supposer compact, s'identifiant à un domaine $k$-analytique d'une courbe lisse $Y$. De la proposition 2.2.1 de \cite{brk5} (elle-même fondée sur le théorème de réduction semi-stable), du fait qu'un domaine $k$-analytique compact d'un espace $k$-affinoïde est réunion finie de domaines rationnels (\cf. \cite{duc2}, lemme 2.4) et de la dernière assertion du lemme~\ref{GERM}.\ref{normeconstante} on déduit que $P$ possède un voisinage $V$ dans $X$ tel que les composantes connexes de $V-\{P\}$ soient toutes des pseudo-disques à l'exception d'un nombre fini d'entre elles qui sont des pseudo-couronnes. Si $P$ est un point de type $(1)$ ou $(4)$ alors $X$ est lisse au voisinage de $P$ et dans ce cas on déduit directement de la proposition 2.2.1 de \cite{brk5} que $P$ possède dans $X$ un voisinage qui est un pseudo-disque. 

\bigskip
\trois{premierjet} En vertu du~\ref{TRI}.\ref{triexist}.\ref{triangloc} ci-dessus on peut, quitte à agrandir $\bf T$, supposer que tout point de $X-{\bf T}$ possède un voisinage qui est ou bien un peudo-disque, ou bien une pseudo-couronne,  que toute composante connexe de $X$ rencontre $\bf T$ et que toute composante connexe de $X-{\bf T}$ est contractile. 

\bigskip
Comme $\bf T$ possède au moins un point sur chaque composante connexe de $X$ le bord de toute composante connexe de $X-{\bf T}$ comprend au moins un point de $\bf T$ ; en particulier une telle composante a au moins un bout. Elle en a au plus deux : en effet chacun des points d'une telle composante possède un voisinage ouvert qui est ou bien un pseudo-disque, ou bien une pseudo-couronne, et qui donc a au plus deux bouts. 

\bigskip
\trois{compunbout} Soit $U$ une composante connexe de $X-{\bf T}$ ayant exactement un bout. Soit $Q$ un point extrémal de $U$. Si $Q$ n'est pas de type $(1)$ ou $(4)$ alors la valeur absolue de $k$ est triviale et $U$ est topologiquement réduit à un intervalle semi-ouvert ne contenant aucun point rigide. Dès lors aucun voisinage de $Q$ dans $U$ ne peut être un pseudo-disque, et aucun ne peut non plus être une pseudo-couronne puisqu'un voisinage connexe de $Q$ est alors homéomorphe à un intervalle semi-ouvert, et n'a en particulier qu'un bout. 

\bigskip
On en déduit que $Q$ est de type $(1)$ ou $(4)$ ; il possède donc un voisinage $V$ dans $U$ qui est un pseudo-disque. Soit $Z$ une sous-pseudo-couronne de $V$ et soit $z$ un point du squelette de $Z$. Alors toutes les composantes connexes de $U-\{z\}$ sont des pseudo-disques, à l'exception éventuelle de celle qui contient l'intervalle ouvert joignant $z$ au bout de $U$. 

\bigskip
\trois{compdeuxbouts} Soit $U$ une composante connexe de $X-{\bf T}$ ayant deux bouts et soit $I$ son squelette. Soit $x$ un point de $I$. Il possède un voisinage $V$ dans $U$ qui est ou bien un pseudo-disque, ou bien une pseudo-couronne ; ce voisinage ayant au moins deux bouts, comme on le voit en considérant son intersection avec $I$, c'est une pseudo-couronne. Son squelette n'est autre que $V\cap I$. Tout intervalle ouvert strictement inclus dans $V\cap I$ et contenant $x$ est encore le squelette d'une pseudo-couronne qui est un voisinage de $x$ dans $U$. Les propriétés topologiques élémentaires de l'intervalle $I$ entraînent alors l'existence d'un intervalle $J$ de $\ZZ$ et d'un ensemble discret $\{x_{i}\}_{i\in J}$ de points de $I$ tels que pour tout $i$ le point $x_{i}$ possède un voisinage dans $U$ qui est une pseudo-couronne de squelette $]x_{i-1},x_{i+1}[$. Ces voisinages recouvrent $U$ et toute composante connexe de $U-\{x_{i}\}_{i\in J}$ est un pseudo-disque ou une pseudo-couronne. 

\bigskip
Soit $P$ un point du bord de $U$ dans $X$. Alors $P$ appartient à $\bf T$ et donc est de type $(2)$ ou $(3)$. Il possède en conséquence un voisinage dont la trace sur $U$ est une pseudo-couronne de squelette inclus dans $I$. La famille $(x_{i})$ ci-dessus peut de ce fait être choisie de manière à être discrète dans l'adhérence {\em topologique} de $U$ dans $X$. 

\bigskip
\trois{conclutriang} Soit $\Omega$ l'ensemble des composantes connexes de $X-{\bf T}$ qui ne sont pas des pseudo-disques. Soit $U$ appartenant à $\Omega$. Il existe un sous-ensemble $E_{U}$ de $U$ formé de points de type $(2)$ ou $(3)$, fermé et discret dans l'adhérence topologique de $U$ dans $X$, et qui est  tel que les composantes connexes de $U-E_{U}$ soient  des pseudo-disques ou des pseudo-couronnes : si $U$ a deux bouts c'est fait au~\ref{TRI}.\ref{triexist}.\ref{compdeuxbouts} ci-dessus. Si $U$ n'en a qu'un, on choisit un point $z$ comme au~\ref{TRI}.\ref{triexist}.\ref{compunbout} puis on applique le~\ref{TRI}.\ref{triexist}.\ref{compdeuxbouts} à la composante de $U-\{z\}$ qui contient l'intervalle joignant $z$ au bout de $U$. 

\bigskip
Soit $P$ un point de $\bf T$. D'après le~\ref{TRI}.\ref{triexist}.\ref{triangloc} il possède un voisinage $V$ tel que toutes les composantes connexes de $V-\{P\}$ à l'exception d'un nombre fini soient des pseudo-disques, les autres étant des pseudo-couronnes. On en déduit que la réunion $\bf S$ de $\bf T$ et des $E_{U}$, où $U$ parcourt $\Omega$, est fermé et discret. Par construction $\bf S$ constitue une triangulation de $X$.~$\Box$ 

\bigskip
\deux{triexistsing} {\bf Proposition.} {\em Soit $X$ une $k$-courbe.

\bigskip
\begin{itemize}

\item[$i)$] Si $X$ possède une triangulation alors le lieu quasi-lisse de $X$ est dense. 

\bigskip
\item[$ii)$] Réciproquement supposons que le lieu  quasi-lisse de $X$ est dense. Alors $X$ possède une triangulation. Plus précisément, soit $\bf T$ un sous-ensemble fermé et discret de $X$ constitué de points ou bien rigides, ou bien de type $(2)$, ou bien de type $(3)$ ; on suppose que $\bf T$ contient tous les points non lisses de $X$. Alors il existe une triangulation de $X$ constituée de $\bf T$ et de points de type $(2)$ ou $(3)$ et telle que toute arête aboutissant à un point rigide de $\bf T$ soit une pseudo-couronne ayant un second sommet de type $(2)$ ou $(3)$.
\end{itemize}
}

\bigskip
{\em Démonstration.} Si $X$ possède une triangulation la réunion des arêtes est un ouvert dense et quasi-lisse de $X$, d'où l'assertion $i)$. Traitons maintenant la réciproque. Comme $X\times_{k}\ra{k}\to X$ induit un homéomorphisme entre les espaces topologiques sous-jacents on peut remplacer $k$ par $\ra{k}$ et donc le supposer {\em parfait}. Soit ${\bf T}_{\tiny \rm rig}$ l'ensemble des points rigides de $\bf T$. Soit $P$ un élément de  ${\bf T}_{\tiny \rm rig}$. Soit $V$ un voisinage ouvert de $P$ tel que $V-\{P\}$ ne rencontre pas $\bf T$ (en particulier $V-\{P\}$ est quasi-lisse). Soit $W$ la normalisation de $V$. Comme $k$ est parfait $W$ est quasi-lisse. 

\bigskip
Soit $r_{P}$ le nombre d'antécédents de $P$ sur $W$ et $P_{1},\ldots,P_{r_{P}}$ les antécédents en question. Fixons $i$. Comme $W$ est quasi-lisse et $k$ parfait il existe un voisinage $W_{i}$ de $P_{i}$ dans $W$ qui ne contient aucun autre antécédent de $P$ et est un $\hres(P_{i})$-disque compact ; fixons un voisinage $W'_{i}$ de $P_{i}$ dans $W_{i}$ qui s'identifie à un $\hres(P_{i})$-disque compact de rayon strictement inférieur. Notons $\eta_{i}$ (resp. $\eta'_{i}$) le point maximal de $W_{i}$ (resp. $W'_{i}$). Notons $U_{i}$ (resp. $U'_{i}$) le $\hres(P_{i})$-disque ouvert de centre $P_{i}$ et de même rayon que $W_{i}$ (resp. $W'_{i}$). Le singleton $\{\eta'_{i}\}$ constitue une triangulation de $U_{i}$. La composante connexe de $P_{i}$ dans $U_{i}-\{\eta'_{i}\}$ n'est autre que $U'_{i}$ qui est un $\hres(P_{i})$-disque ; comme $U'_{i}-\{P_{i}\}$ est une $\hres(P_{i})$-couronne $\{\eta'_{i}\}$ constitue également une triangulation de $U'_{i}-\{P_{i}\}$. 

\bigskip
Le morphisme $W\to V$ induit un isomorphisme entre $W_{i}-\{P_{i}\}$ et un domaine $k$-analytique $\Omega_{i}$ de $V-\{P\}$. Soit $\omega_{P,i}$ (resp. $\omega'_{P,i}$) l'image de $\eta_{i}$ (resp. $\eta'_{i}$) sur $\Omega_{i}$. Soit $Z_{P,i}$  l'image du domaine $U_{i}-\{P_{i}\}$ dans $\Omega_{i}$. Alors $\{\omega'_{P,i}\}$ constitue une triangulation de $Z_{P,i}$. Soit $Y$ le complémentaire dans $X$ de la réunion de ${\bf T}_{\tiny \rm rig}$ et des $Z_{P,i}$ où $P$ parcourt l'ensemble des points rigides de $\bf T$ et où $i$ varie, $P$ étant fixé, entre $1$ et $r_{P}$. C'est un domaine $k$-analytique fermé de $X$ qui est quasi-lisse puisque $\bf T$ contient l'ensemble des points non lisses de $X$. La réunion des $\omega_{P,i}$ et de ${\bf T}-{\bf T}_{\tiny \rm rig}$ constitue un sous-ensemble fermé et discret ${\bf E}$ de $Y$ dont tous les points sont de type $(2)$ ou $(3)$. Par la proposition~\ref{TRI}.\ref{triexist} il existe une triangulation $\bf F$ de $Y$ contenant $\bf E$ et ne comprenant que des points de type $(2)$ ou $(3)$. La réunion $\bf S$ de $\bf F$, de ${\bf T}_{\tiny \rm rig}$ et des $\omega'_{P,i}$ est une triangulation de $X$ satisfaisant aux conditions requises.~$\Box$ 

\section{Un analogue analytique du complexe de Bloch-Ogus}\label{COH}

\setcounter{cpt}{0}
\subsection*{Préliminaires sur la cohomologie à support compact} 

\bigskip
\deux{notlambda} {\bf Quelques notations.} Si $j$ est une immersion ouverte entre espaces $k$-analytiques paracompacts et séparés on notera $j_{c}$ le foncteur induit entre les groupes de cohomologie à support compact. Si $h$ est une classe de cohomologie à support compact (sur un espace convenable) on notera $h_{\bullet}$ son image dans la cohomologie sans support. Soit $X$ un espace $k$-analytique séparé et paracompact. Pour tout faisceau $\cal F$ sur $X$ et tout entier $q$ on dispose d'un morphisme (qui s'insère dans une suite exacte naturelle de cohomologie à support) $$\lim_{\stackrel{\to}{K}} \H^{q}(X-K,{\cal F})\to \H^{q+1}_{c}(X,{\cal F})$$ où $K$ parcourt l'ensemble des parties compactes de $X$. Soit $K$ un compact de $X$, soit $\xi$ une fonction inversible sur $X-K$ et soit $n$ un entier inversible dans $k$. L'adjonction d'une racine $n$-ième de $\xi$ définit une classe de $\H^{1}(X-K,\mu_{n})$ ; l'image de ladite classe dans $\H^{2}_{c}(X,\mu_{n})$ sera notée $\Lambda_{n,\xi}$. Par construction $\Lambda_{n,\xi}$  provient de $\H^{1}_{c}(X,\gm)$ {\em via} le cobord de la suite de Kummer. 

\bigskip
\deux{hrestz} Soit $X$ un espace $k$-analytique paracompact et séparé dont les domaines $k$-affinoïdes sont géométriquement réduits et soit $Y$ un ouvert connexe et non vide de $X$ ; notons $j$ l'immersion $Y\hookrightarrow X$. Sous ces hypothèses $\got{c}(Y)$ est une extension finie séparable de $k$. On la notera $L$ et pour tout espace $k$-analytique $V$ on désignera par $V_{L}$ l'espace $V\times_{k}L$. L'ouvert $Y$ s'identifie à une composante connexe $Z$ de $Y_{L}$ {\em via} une section de la première projection. Si $h$ est une classe de cohomologie sur $Y$ on notera $h_{|Z}$ l'élément de la cohomologie de $Y_{L}$ dont la restriction à $Z$ correspond à $h$ modulo l'isomorphisme $Z\simeq Y$ et qui est nul sur les autres composantes. On a alors $\mbox{Cor}_{L/k}(h_{|Z})=h$ où $\mbox{Cor}$ désigne la corestriction. 

\bigskip
\deux{cupclasseconst} Soit $\cal F$ (resp. $\cal G$) un faisceau en groupes abéliens sur $k$ (resp. sur $X$). Donnons-nous deux entiers $l$ et $q$, une classe $h$ appartenant à $\H^{l}(L,{\cal F})$ et une classe $\eta$ dans le groupe $\H^{q}_{c}(Y,{\cal G})$. On peut écrire $h \cup\eta=\mbox{Cor}_{L/k}(h_{|Z}\cup \eta_{|Z})$ et donc par fonctorialité $j_{c}(h \cup\eta)=\mbox{Cor}_{L/k}(j_{c}(h_{|Z}\cup \eta_{|Z}))$. Comme $L$ se plonge dans l'anneau des fonctions de $X_{L}$ la classe $h_{|Z}$ se prolonge en une classe de $\H^{l}(X_{L},{\cal F})$, qui est simplement la restriction (notée encore $h$) de la classe $h$ de $\H^{l}(L,{\cal F})$. On peut dès lors écrire $$j_{c}(h \cup\eta)=\mbox{Cor}_{L/k}( h\cup j_{c}(\eta_{|Z})) $$ et donc $$j_{c}(h \cup\eta)_{\bullet}=\mbox{Cor}_{L/k}( h\cup (j_{c}(\eta_{|Z}))_{\bullet}\;).$$ Les deux lemmes ci-dessous joueront un rôle crucial dans la démonstration du théorème de comparaison (\cf. th.~\ref{COMP}.\ref{theocomp} {\em infra}). 

\bigskip
\setcounter{cptbis}{0}
\trois{toptriv} {\bf Lemme.} {\em On garde les notations introduites ci-dessus. On suppose que $\cal G$ est égal au faisceau constant $\ZZ/n$ pour un certain entier $n$,  que la classe $\eta$ provient de $\H^{q}_{c}(|Y|,\ZZ/n)$ et que $\H^{q}(|X_{L}|,\ZZ/n)$ est trivial. Alors $j_{c}(h \cup\eta)_{\bullet}$ est nulle.} 

\bigskip
{\em Démonstration.}  Dans ce cas $\eta_{|Z}$ provient de $\H^{q}_{c}(| Y_{L}|,\ZZ/n)$ et $j_{c}(\eta_{|Z})_{\bullet}$ appartient en conséquence à l'image de $\H^{q}(|X_{L}|,\ZZ/n)\to \H^{q}(X_{L},\ZZ/n)$ ; elle est de ce fait nulle en vertu des hypothèses. Comme $j_{c}(h \cup\eta)_{\bullet}$ est égal à $\mbox{Cor}( h\cup (j_{c}(\eta_{|Z}))_{\bullet}\;)$ elle est triviale.~$\Box$ 

\bigskip
\trois{gagafibre}  {\bf  Lemme.} {\em On reprend les notations du~{\rm \ref{COH}.\ref{notlambda} et du~\ref{COH}.\ref{hrestz}}. Soit $n$ un entier inversible dans $k$ et soit $M$ une extension finie séparable de $L$.  Soit $\cal F$ un faisceau en groupes abéliens sur $k$ donné par un module galoisien fini et annulé par $n$. Faisons l'hypothèse qu'il existe une immersion ouverte de $X$ dans l'analytifié ${\cal X}\an$ d'une $k$-variété {\em algébrique,  projective et lisse} $\cal X$. Soit $l$ un entier et soit $\lambda$ une classe appartenant à $\H^{l+2}_{c}(Y,{\cal F}(1))$ de la forme $\mbox{\rm Cor}_{M/L}(h\cup \Lambda_{n,\xi})$ où $h$ appartient à $\H^{l}(M,{\cal F})$ et où $\xi$ est une fonction définie et inversible sur le complémentaire d'un compact de $Y\times_{L}M$. Alors il existe un recouvrement $(U_{i})$ de $\cal X$ par des ouverts de Zariski tel que la restriction de $j_{c}(\lambda)_{\bullet}$ à $U_{i}\an\cap X$ soit triviale pour tout $i$. }

\bigskip
{\em Démonstration.} Notons que comme ${ \cal X}$ est projective la cohomologie à support compact de  ${\cal X}\an$ coïncide avec la cohomologie sans support. La commutativité du diagramme $$\diagram \H^{l+2}_{c}(Y,.)\rto\dto&\H^{l+2}_{c}(X,.)\rto\dto&\H^{l+2}({\cal X}\an,.)\dto\\ \H^{l+2}(Y,.)&\H^{l+2}(X,.)\lto&\H^{l+2}({\cal X}\an,.)\lto \enddiagram$$ permet de supposer que $X$ est égal à ${\cal X}\an$. 

\bigskip
On sait d'après le~\ref{COH}.\ref{hrestz} que $\lambda$ est égale à $\mbox{Cor}_{L/k}(\lambda_{|Z})$ et $\lambda_{|Z}$ est elle-même la corestriction (de $M$ à $L$) de l'élément  $(h\cup \Lambda_{n,\xi})_{|Z\times_{L}M}$ de $\H^{l+2}_{c}(Y_{M},{\cal F}(1))$, qui est défini comme étant égal à $h\cup \Lambda_{n,\xi}$ sur $Z\times_{L}M$ et à zéro sur les autres composantes. Par transitivité on a donc $$\lambda=\mbox{Cor}_{M/k}(\;(h\cup \Lambda_{n,\xi})_{|Z\times_{L}M}\;).$$ La classe $j_{c}((h\cup \Lambda_{n,\xi})_{|Z\times_{L}M})$ de $\H^{l+2}_{c}({\cal X}\an_{M},{\cal F}(1))$ peut se réécrire $$h\cup j_{c}((\Lambda_{n,\xi})_{|Z\times_{L}M})$$ puisque $h$ provient de la cohomologie de $M$. D'autre part  $j_{c}((\Lambda_{n,\xi})_{|Z\times_{L}M})$ se déduit, {\em via} le cobord de la suite de Kummer, d'un élément de $\H^{1}({\cal X}_{M}\an,\gm)$. Or $\cal X$ étant projective ce dernier groupe s'identifie, par GAGA et par le théorème 90 de Hilbert, à $\H^{1}(|{\cal X}_{M}|,\gm)$. On en déduit que $j_{c}((\Lambda_{n,\xi})_{|Z\times_{L}M})$ est l'image d'une classe de cohomologie étale du {\em schéma} ${\cal X}_{M}$ qui est localement triviale pour la topologie de Zariski. En conséquence $h\cup j_{c}((\Lambda_{n,\xi})_{|Z\times_{L}M})$ provient elle-même d'une certaine $\mu$ appartenant à $\H^{l +2}({\cal X}_{M},{\cal F}(1)$ et localement triviale pour la topologie de Zariski.

\bigskip

La classe $j_{c}(\lambda)$ est égale à l'image de $\mbox{Cor}_{M/k}(\mu)$ dans $\H^{l+2}({\cal X}\an,{\cal F}(1))$. Comme $\mu$ est localement triviale sur ${\cal X}_{M}$ pour la topologie de Zariski  $\mbox{Cor}_{M/k}(\mu)$ est nulle aux points génériques de ${\cal X}$. La variété $\cal X$ étant lisse ceci implique (Gabber, \cf.~\cite{jlct}, prop. 2.1.2 et th. 2.2.1) que $\mbox{Cor}_{M/k}(\mu)$ est localement triviale pour la topologie de Zariski, ce qui permet de conclure.~$\Box$

\subsection*{Complexes associés à une triangulation} 

\bigskip
\deux{defmodfasc} {\bf Définition.} {\em On appellera faisceau {\em raisonnable} sur un espace $k$-analytique  tout faisceau étale {\em en groupes abéliens} dont la restriction à chaque composante connexe de l'espace en question provient d'un module fini et de cardinal premier à $p$ sur le groupe fondamental géométriquement modéré.} 

\bigskip
\deux{introcohsup} Soit $X$ un espace $k$-analytique séparé et paracompact et soit $\bf S$ une partie fermée et discrète de $X$. Notons $\underline{\bf S}$ l'ensemble des complémentaires de voisinages ouverts de $\bf S$. Alors $\underline{\bf S}$ est une {\em famille de supports} de l'espace $X-{\bf S}$ au sens de \cite{brk2}, \S 5.1. La famille $\underline{\bf S}$ est paracompactifiante. Soit $\cal F$ un faisceau sur $X$. La suite exacte $ii)$ du paragraphe 5.2.1 de \cite{brk2}, appliquée en prenant pour $\Phi$ l'ensemble de tous les fermés de $X$, se réécrit de la manière suivante en tenant compte de la proposition 4.3.4 de  {\em loc. cit.} : $$\ldots\to \H^{i}(X,{\cal F})\to \prod_{P\in {\bf S}} \H^{i}(\hres(P),{\cal F}_{P})\to \H^{i+1}_{\underline{\bf S}}(X-{\bf S},{\cal F})\to \H^{i+1}(X,{\cal F})\to\ldots $$

Si l'on prend maintenant pour $\Phi$ la famille des parties compactes de $X$ la suite exacte  $ii)$ du paragraphe 5.2.1 de \cite{brk2} devient $$\ldots\to \H^{i}_{c}(X,{\cal F})\to \bigoplus_{P\in {\bf S}} \H^{i}(\hres(P),{\cal F}_{P})\to \H^{i+1}_{c}(X-{\bf S},{\cal F})\to \H^{i+1}_{c}(X,{\cal F})\to\ldots $$

\bigskip
\deux{deuxsuites} Soit $X$ une $k$-courbe. Soit $\cal F$ un faisceau sur $X$. La dimension de $|X|$ étant inférieure ou égale à 1 le morphisme naturel  entre les suites spectrales $$\H^{p}(|X|,{\rm R}^{q}\pi_{*}{\cal F})\Rightarrow \H^{p+q}(X,{\cal F})\;\mbox{et}\;  \H^{p}_{c}(|X|,{\rm R}^{q}\pi_{*}{\cal F})\Rightarrow \H^{p+q}_{c}(X,{\cal F})$$ fournit pour tout entier $i$ positif ou nul un diagramme commutatif à lignes exactes : $$\diagram 0\rto& \H^{1}(|X|, {\rm  R}^{i-1}\pi_{*}{\cal F})\rto \dto & \H^{i}(X,{\cal F})\rto\dto &  \H^{0}(|X|,{\rm R}^{i}\pi_{*}{\cal F})\dto\rto& 0 \\ 0\rto& \H^{1}_{c}(|X|, {\rm R}^{i-1}\pi_{*}{\cal F})\rto &\H^{i}_{c}(X,{\cal F})\rto& \H^{0}_{c}(|X|,{\rm  R}^{i}\pi_{*}{\cal F})\rto& 0\enddiagram$$

\bigskip
\deux{nulsurs} {\bf Lemme.} {\em {\em Supposons que $\cal F$ est raisonnable.} Soit $\bf S$ une triangulation de $X$ et soit $i$ un entier. Le noyau de $\H^{i}(X,{\cal F})\to \prod\limits_{P\in {\bf S}} \H^{i}(\hres(P),{\cal F}_{P})$  est égal à celui de $\H^{i}(X,{\cal F})\to \H^{0}(|X|,{\rm  R}^{i}\pi_{*}{\cal F})$ ; le noyau de $$\H^{i}_{c}(X,{\cal F})\to \bigoplus\limits_{P\in {\bf S}} \H^{i}(\hres(P),{\cal F}_{P})$$ est égal à celui de $\H^{i}_{c}(X,{\cal F})\to \H^{0}_{c}(|X|,{\rm R}^{i}\pi_{*}{\cal F})$ .}

\bigskip
{\em Démonstration.} Soit $h$ une classe de $\H^{i}(X,{\cal F})$. Si l'image de $h$ dans le groupe $\H^{0}(|X|,{\rm R}^{i}\pi_{*}{\cal F})$ est  nulle alors $h$ est nulle en tout point de $X$, et en particulier en tout sommet de $\bf S$. Réciproquement, supposons que $h$ est nulle en tout sommet de $\bf S$. Soit $Y$ une arête de $\bf S$. Comme $\bf S$ est une triangulation $Y$ possède au moins un sommet $P$. Comme $h$ est nulle en $P$ elle est nulle au voisinage de $P$. Si $Y$ est une pseudo-couronne alors l'un des deux bouts de $S(Y)$ (au moins) aboutit à $P$ et $h$ est donc nulle en au moins un point de $S(Y)$. On en déduit grâce à la proposition~\ref{GERM}.\ref{courskel}.\ref{locglob} que $h_{|Y}$ est triviale. Si $Y$ est un pseudo-disque alors $Y$ possède une sous-peudo-couronne $Z$ dont le squelette aboutit à $P$. Il existe donc un point de $S(Z)$ en lequel $h$ est nulle. On déduit du~\ref{GERM}.\ref{loczeropseudisc} que $h_{|Y}$ est triviale. Si $x$ est un point de $X$ alors ou bien $x$ appartient à $\bf S$ auquel cas $h(x)$ est nulle par hypothèse, ou bien $x$ est sur une arête auquel cas $h(x)$ est nulle d'après ce qui précède. En conséquence $h$ appartient au noyau de $\H^{i}(X,{\cal F})\to \H^{0}(|X|,{\rm R}^{i}\pi_{*}{\cal F})$. 

\bigskip
Le résultat pour la cohomologie à supports compacts s'ensuit aussitôt compte-tenu de l'injectivité de $\H^{0}_{c}(|X|,{\rm R}^{i}\pi_{*}{\cal F})\to \H^{0}(|X|,{\rm  R}^{i}\pi_{*}{\cal F})$.~$\Box$ 

\bigskip
\deux{complexe}{\bf Deux complexes de type "Bloch-Ogus".}  On désigne toujours par $X$ une $k$-courbe, par $\bf S$ une triangulation de $X$ et par $\cal F$ un faisceau raisonnable sur $X$. Soit $q$ un entier positif. On note $\rqt{\cal F}$ (resp. $\rqtc{\cal F}$) le complexe à deux termes concentrés en degrés zéro et un $$\prod\limits_{P\in {\bf S}} \H^{q}(\hres(P),{\cal F}_{P})\to \H^{q+1}_{\underline{\bf S}}(X-{\bf S},{\cal F})$$ (resp. $\bigoplus\limits_{P\in {\bf S}} \H^{q}(\hres(P),{\cal F}_{P})\to \H^{q+1}_{c}(X-{\bf S},{\cal F})$). Du lemme~\ref{COH}.\ref{nulsurs} ainsi que des suites exactes du~\ref{COH}.\ref{introcohsup} et du~\ref{COH}.\ref{deuxsuites} on déduit le théorème suivant :  

\bigskip
\deux{cohocomplexe} {\bf Théorème.} {\em Soit $i$ appartenant à $\{0,1\}$. Le groupe $\H^{i}(\rqt{\cal F})$ s'identifie à $\H^{i}(|X|,{\rm R}^{q}\pi_{*}{\cal F})$ ; le groupe $\H^{i}(\rqtc{\cal F})$ s'identifie à $\H^{i}_{c}(|X|,{\rm R}^{q}\pi_{*}{\cal F})$.}~$\Box$

\subsection*{Quelques calculs explicites} 

\deux{calcexpl} On suppose toujours que $X$ est une $k$-courbe et que $\bf S$ en constitue une triangulation. Soit $\cal F$ un faisceau sur $k$ donné par un module galoisien fini annulé par un entier $n$ premier à $p$. Soit $q$ un entier. Soit $Y$ une arête de $\bf S$. On suppose que l'un au moins des sommets de $Y$ est un point quasi-lisse de $X$ ; on choisit un tel sommet et on le note $P$. On se propose de décrire le groupe $\H^{q+1}_{\underline{\bf S}_{|Y}}(Y,{\cal F})$. On distingue trois cas :

\setcounter{cptbis}{0}

\bigskip
\trois{ydisc} {\bf Le cas où $Y$ est un pseudo-disque.} Dans ce cas $P$ est l'unique sommet de $Y$ et $Y\cup\{P\}$ est compact. Dès lors la famille $\underline{\bf S}_{|Y}$ est exactement celle des compacts de $Y$. On dispose d'une suite exacte $$\ldots \to \H^{i}(Y,{\cal F})\to\lim_{\stackrel{\to}{K \tiny \mbox{compact}}} \H^{i}(Y-K,{\cal F}) \to \H^{i+1}_{c}(Y,{\cal F})\to \H^{i+1}(Y,{\cal F})\to\ldots$$ 

Soit $Z$ une sous-pseudo-couronne de $Y$ et soit $i$ un entier. Le lemme~\ref{PSEU}.\ref{etalegalois} assure que $\H^{i}(Y,{\cal F})$ s'identifie à $\H^{i}(\got {c}(Y),{\cal F})$ ; d'autre part $\got{c}(Z)$ est égal à $\got{c}(Y)$ et $\H^{i}(\got {c}(Z),{\cal F})\to \H^{i}(Z,{\cal F})$ est injective. Enfin si $T$ est une sous-pseudo-couronne de $Z$ alors $\H^{i}(Z,{\cal F})\to \H^{i}(T,{\cal F})$ est un isomorphisme d'après la remarque~\ref{PSEU}.\ref{remsouspseu}. 

\bigskip
Soit $\cal E$ le système de composantes adhérent à $P$ que définit $Y$. En reprenant les notations du~\ref{GERM}.\ref{systpseuco} et du~\ref{GERM}.\ref{inductetale} le corps $\got{c}({\cal E})$ est égal à $\got{c}(Y)$ ; la famille des sous-pseudo-couronnes de $Y$ est par ailleurs cofinale dans celle des complémentaires de parties compactes. On déduit dès lors de ce qui précède que $\H^{q+1}_{c}(Y,{\cal F})$ est naturellement isomorphe à $\H^{q}({\cal E},{\cal F})/\H^{q}(\got{c}({\cal E}),{\cal F})$.  D'après le~\ref{GERM}.\ref{constcour} il existe une fonction $\xi$ appartenant à ${\cal M}_{P}$ telle que $h\mapsto (\xi)\cup h$ induise un isomorphisme $$\H^{q-1}(\got{c}({\cal E}),{\cal F}(-1))\simeq \H^{q}({\cal E},{\cal F})/\H^{q}(\got{c}({\cal E}),{\cal F}).$$ 

\bigskip
\noindent
{\small {\bf Remarque.} La flèche naturelle $\H^{q}(\hres(P),{\cal F})\to \H^{q+1}_{\underline{\bf S}_{|Y}}(Y,{\cal F})$ déduite de la différentielle du complexe $\rqt{\cal F}$ est, modulo l'isomorphisme ci-dessus, induite par la restriction  $$\H^{q}(\hres(P),{\cal F})\to \H^{q}({\cal O}({\cal E})_{\rm alg},{\cal F})\simeq \H^{q}({\cal E},{\cal F}).$$}

\bigskip
Ce qui précède peut se traduire par l'existence d'un isomorphisme $$\H^{q-1}(\got{c}(Y),{\cal F}(-1))\simeq \H^{q+1}_{c}(Y,{\cal F})=\H^{q+1}_{\underline{\bf S}_{|Y}}(Y,{\cal F})$$ donné par $h\mapsto h\cup \Lambda_{n,\xi}$, la notation $\Lambda_{n,\xi}$ ayant été introduite au~\ref{COH}.\ref{notlambda}. 

\bigskip
\trois{ypseucofuite} {\bf Le cas où $Y$ est une pseudo-couronne dont l'adhérence n'est pas compacte.} Dans cette situation $P$ est encore le seul sommet de $Y$ et $Y\cup \{P\}$ a exactement un bout : l'extrémité ouverte de $S(Y)\cup\{P\}$. La famille des  sous-pseudo-couronnes de $Y$ qui aboutissent à $P$ est cofinale dans celle des complémentaires de fermés de $\underline{\bf S}_{|Y}$. Or si $Z$ est une sous-pseudo-couronne de $Y$ la flèche $\H^{i}(Y,{\cal F})\to \H^{i}(Z,{\cal F})$ est un isomorphisme pour tout $i$ (\cf. remarque~\ref{PSEU}.\ref{remsouspseu}). On en déduit que $\H^{q+1}_{\underline{S}_{|Y}}(Y,{\cal F})$ est nul. 

\bigskip
\trois{ypseucocomp} {\bf Le cas où $Y$ est une pseudo-couronne dont l'adhérence est compacte.} La famille $\underline{\bf S}_{|Y}$ est alors celle des parties compactes de $Y$. La famille des ouverts de $Y$ qui sont la réunion disjointe de deux sous-pseudo-couronnes ayant chacune un bout commun avec $Y$ est cofinale dans celle des complémentaires de parties de $\underline{\bf S}_{|Y}$. Pour toute sous-pseudo-couronne $Z$ de $Y$ et tout entier $i$ la flèche $\H^{i}(Y,{\cal F})\to \H^{i}(Z,{\cal F})$ est un isomorphisme d'après la remarque~\ref{PSEU}.\ref{remsouspseu}. De ces faits et de la suite exacte mentionnée au~\ref{COH}.\ref{calcexpl}.\ref{ydisc} découle l'existence d'une suite exacte $$0\to \H^{q}(Y,{\cal F}) \to \H^{q}(Y,{\cal F})\times \H^{q}(Y,{\cal F})\to  \H^{q+1}_{c}(Y,{\cal F}) \to 0$$ où la flèche de gauche est la diagonale. 

\bigskip
\noindent
{\small {\bf Remarque.} Soit $Q$ l'un des sommets de $Y$. Supposons d'abord qu'il n'adhère qu'à un bout de $Y$ et soit ${\cal E}$ le système de composantes correspondant. La flèche naturelle $\H^{q}(\hres(Q),{\cal F})\to \H^{q+1}_{\underline{\bf S}_{|Y}}(Y,{\cal F})$ déduite de la différentielle du complexe $\rqt{\cal F}$ est alors, modulo l'isomorphisme ci-dessus, donnée par la composée de $$\H^{q}(\hres(Q),{\cal F})\to \H^{q}({\cal O}({\cal E})_{\rm alg},{\cal F})\simeq \H^{q}({\cal E},{\cal F})\simeq \H^{q}(Y,{\cal F})$$ et de l'identification de   $\H^{q}(Y,{\cal F})$ avec le facteur de $\H^{q}(Y,{\cal F})\times \H^{q}(Y,{\cal F})$ qui correspond au bout en question. Si $Q$ adhère aux deux bouts de $Y$ on note ${\cal E}$ et ${\cal E}'$ les deux systèmes de composantes concernés et la flèche étudiée s'identifie alors à 
$$\H^{q}(\hres(Q),{\cal F})\to  \H^{q}({\cal E},{\cal F})\times \H^{q}({\cal E}',{\cal F})\simeq \H^{q}(Y,{\cal F})\times \H^{q}(Y,{\cal F}).$$}

\bigskip
Le choix d'un des deux bouts de $Y$ (ou, si l'on préfère, d'une orientation) fournit en conséquence un isomorphisme $\H^{q+1}_{c}(Y,{\cal F})\simeq \H^{q}(Y,{\cal F})$ ; on fixe un tel choix de sorte que $P$ adhère au bout en question. 

\bigskip
Explicitons l'isomorphisme ci-dessus, ou plus précisément sa réciproque. Le cup-produit avec la section $\overline{1}$ de $\H^{0}(Y,\ZZ/n)$ est l'identité de $\H^{q}(Y,{\cal F})$ modulo l'isomorphisme ${\cal F}\simeq {\cal F}\otimes \ZZ/n$. Par fonctorialité la flèche $\H^{q}(Y,{\cal F})\to \H^{q+1}_{c}(Y,{\cal F})$ est donc donnée par le cup-produit avec l'image $\varpi_{n}$ de $\overline{1}$ dans $\H^{1}_{c}(Y,\ZZ/n)$. La classe $\varpi_{n}$ se décrit facilement : on se donne deux sous-pseudo-couronnes $Z$ et $Z'$ de $Y$ disjointes et ayant chacune un bout commun avec $Y$. On suppose que $Z$ est celle qui est située du côté que l'on a choisi. Le complémentaire de $Z\coprod Z'$ dans $Y$ est compact. On considère la section du $\ZZ/n$-torseur trivial sur $Z\coprod Z'$ qui vaut $\overline{1}$ sur $Z$ et $\overline{0}$ sur $Z'$ ; l'image de cette section par la flèche $\lim\limits_{\to} \H^{0}(Y-K,\ZZ/n)\to \H^{1}_{c}(Y,\ZZ/n)$ est égale à $\varpi_{n}$. Notons que $\varpi_{n}$ provient de la classe de $\H^{1}_{c}(|Y|,\ZZ/n)$ définie de manière analogue (laquelle est l'image réciproque par la rétraction de la classe fondamentale modulo $n$ de l'intervalle orienté $S(Y)$). 

\bigskip
Du~\ref{GERM}.\ref{constcour} et du~\ref{GERM}.\ref{courskel}.\ref{exactescinde} découle l'existence d'une sous-pseudo-couronne $Z$ de $Y$ aboutissant à $P$ (que l'on peut toujours supposer différente de $Y$) et d'une fonction inversible $\xi$ sur $Z$ telle que $(\lambda,h) \mapsto \lambda + h\cup (\xi)$ définisse un isomorphisme $\H^{q}(\got{c}(Y),{\cal F})\oplus \H^{q}(\got{c}(Y),{\cal F}(-1))\simeq \H^{q}(Y,{\cal F})$, où la classe $(\xi)$ est vue comme appartenant à $\H^{1}(Y,\mu_{n})$ par le biais de $\H^{q}(Z,{\cal F})\simeq \H^{q}(Y,{\cal F})$. 

\bigskip
Fixons une sous-pseudo-couronne $Z'$ de $Y$ disjointe de $Z$ et contenant l'autre bout de $Y$. Le cup-produit de $(\xi)$ (vue comme classe de $\H^{1}(Y,\mu_{n})$) avec $\varpi_{n}$ est par construction l'élément $\Lambda_{n,\vec{\xi}}$ de $\H^{2}_{c}(Y,\mu_{n})$ où $\vec{\xi}$ désigne la fonction inversible sur $Z\coprod Z'$ qui vaut $\xi$ sur $Z$ et $1$ sur $Z'$. 

\bigskip
En conclusion le groupe $\H^{q+1}_{\underline{\bf S}_{|Y}}(Y,{\cal F})$ est égal à $\H^{q+1}_{c}(Y,{\cal F})$ et chaque élément de ce dernier possède une écriture de la forme $h_{1}\cup \varpi_{n}+h_{2}\cup \Lambda_{n,\vec{\xi}}$ où $h_{1}$ (resp. $h_{2}$) appartient à $\H^{q}(\got{c}(Y),{\cal F})$ (resp. à $\H^{q-1}(\got{c}(Y),{\cal F}(-1)$)~). 

\section{Un théorème de comparaison}\label{COMP}  

\setcounter{cpt}{0} 

\deux{potbonne} Si $\cal X$est  une courbe {\em algébrique} sur $k$ et si $L$ est une extension de $k$ on notera ${\cal X}_{L}$ la courbe ${\cal X}\times_{k}L$. On dira que ${\cal X}$ a {\em potentiellement bonne réduction} s'il existe une extension finie $L$ de $k$ et une immersion ouverte de $ {\cal X}_{L}$ dans la fibre générique d'une $L\zero$-courbe relative propre et lisse. Il revient au même de demander que ${\cal X}$ soit lisse et que $|{\cal X}\an_{L}|$ soit simplement connexe pour toute extension $L$ de $k$.

\bigskip
\deux{diversdiag} Soit $\cal X$ une $k$-courbe algébrique et soit ${\cal F}$ un faisceau sur ${\cal X}$. On note encore $\cal F$ le faisceau correspondant sur ${\cal X}\an$. Pour tout $q$ on dispose d'un diagramme commutatif à lignes exactes {\small $$\diagram 0\rto &\H^{1}(|{\cal X}|,{\rm R}^{q}\pi_{*}{\cal F})\rto \dto &\H^{q+1}({\cal X},{\cal F})\rto \dto & \H^{0 }(|{\cal X}|,{\rm R}^{q+1}\pi_{*}{\cal F}) \rto \dto&0\\  0\rto &\H^{1}(|{\cal X}\an|,{\rm R}^{q}\pi_{*}{\cal F})\rto &\H^{q+1}({\cal X}\an,{\cal F})\rto  & \H^{0 }(|{\cal X}\an|,{\rm R}^{q+1}\pi_{*}{\cal F}) \rto & 0\enddiagram$$} dans lequel la flèche verticale du milieu est un isomorphisme d'après Berkovich dès que $\cal F$ est constructible de torsion première à $p$ (\cite{brk2}, cor. 7.5.4). {\em Plaçons-nous sous cette dernière hypothèse}. La flèche $\H^{0}(|{\cal X}|,{\rm R}^{q+1}\pi_{*}{\cal F})\to  \H^{0}(|{\cal X}\an|,{\rm R}^{q+1}\pi_{*}{\cal F})$ est alors surjective par un argument de chasse au diagramme élémentaire, et $\H^{1}(|{\cal X}|,{\rm R}^{q}\pi_{*}{\cal F})\to \H^{1}(|{\cal X}\an|,{\rm R}^{q}\pi_{*}{\cal F})$ est injective. Si l'un de ces deux morphismes est un isomorphisme, il en va de même de l'autre. On se propose de donner deux conditions suffisantes pour qu'il en soit ainsi. Afin d'en rendre l'énoncé plus digeste, introduisons une définition. 

\bigskip
\deux{hypoh} {\bf Définition.} {\em Soit $K$ un corps, soient $n$ et $q$ deux entiers et soit $\cal G$ un faisceau en $\ZZ/n$-modules sur $k$.  On dira qu'une classe appartenant à $\H^{q}(K,{\cal G})$ est {\em décomposable } si elle peut s'écrire comme une somme de termes de la forme $\mbox{\rm Cor}_{M/K}(\alpha \cup\beta)$ où :

\bigskip
\begin{itemize} 
\itb $M$ est une extension finie séparable de $k$ ;
\bigskip
\itb $\alpha$ appartient à $\H^{1}(M,\mu_{n})$ et $\beta$ à $\H^{q-1}(M,{\cal F}(-1))$.

\end{itemize}}

\bigskip
\deux{theocomp} {\bf Théorème.} {\em Soit $k$ un corps ultramétrique complet et soit $p$ l'exposant caractéristique du corps résiduel $\red{k}$. Soit $n$ un entier premier à $p$ et soit $\cal F$ un faisceau sur $k$ donné par un module galoisien fini de $n$-torsion. Soit $\cal X$ une courbe algébrique lisse sur $k$ et ${\cal X}\an$ l'espace analytique associé. 

\bigskip
\begin{itemize}
\item[$i)$]  Si $\cal X$ a potentiellement bonne réduction alors quel que soit $q$ les flèches $$\H^{0}(|{\cal X}|,{\rm R}^{q+1}\pi_{*}{\cal F})\to  \H^{0}(|{\cal X}\an|,{\rm R}^{q+1}\pi_{*}{\cal F})$$ et $\H^{1}(|{\cal X}|,{\rm R}^{q}\pi_{*}{\cal F})\to \H^{1}(|{\cal X}\an|,{\rm R}^{q}\pi_{*}{\cal F})$ sont des isomorphismes. 

\bigskip
\item[$ii)$] Soit $q$ un entier tel que pour toute extension finie séparable $L$ de $k$ toute classe de $\H^{q}(L,{\cal F})$ soit décomposable. Alors les flèches $$\H^{0}(|{\cal X}|,{\rm R}^{q+1}\pi_{*}{\cal F})\to  \H^{0}(|{\cal X}\an|,{\rm R}^{q+1}\pi_{*}{\cal F})$$ et $\H^{1}(|{\cal X}|,{\rm R}^{q}\pi_{*}{\cal F})\to \H^{1}(|{\cal X}\an|,{\rm R}^{q}\pi_{*}{\cal F})$ sont des isomorphismes.

\end{itemize}
}

\bigskip
{\em Démonstration.} Les flèches radicielles étant sans effet sur les topoi en jeu on peut supposer que ${\cal X}$ s'identifie à un ouvert dense d'une courbe projective lisse $\overline{\cal X}$. Soit $\bf E$ le complémentaire de  ${\cal X}\an$ dans $\overline{\cal X}\an$. La proposition~\ref{TRI}.\ref{triexistsing} assure l'existence d'une triangulation ${\bf T}$ de ${\cal X}\an$ qui est réunion de $\bf E$ et d'un ensemble de points de type $(2)$ ou $(3)$ (dont on peut supposer qu'il rencontre chaque composante connexe de ${\cal X}\an$) et qui est telle que toute arête aboutissant à un point de $\bf E$ soit une pseudo-couronne. Par compacité de $\overline{\cal X}\an$ l'ensemble $\bf T$ est fini. Son intersection avec ${\cal X}\an$ constitue une triangulation de cette dernière, que l'on note $\bf S$.

\bigskip
Pour chacune des assertions $i)$ et $ii)$ il suffit de montrer que si $h$ est une classe de $\H^{q+1}({\cal X}\an,{\cal F})$ localement triviale sur l'espace topologique $|{\cal X}\an|$ alors il existe un recouvrement $(U_{i})$ de ${\cal X}$ par des ouverts de Zariski tel que la restiction de $h$ à $U_{i}\an$ soit nulle pour tout $i$. Tout élément du noyau de $\H^{q+1}({\cal X}\an,{\cal F})\to \H^{0}(|{\cal X}\an|,{\rm R}^{q+1}\pi_{*}{\cal F})$ provient de $\H^{q+1}_{\underline{\bf S}}(X-{\bf S},{\cal F})$ d'après le lemme~\ref{COH}.\ref{nulsurs}. Le cardinal de $\bf S$ est fini, et tout voisinage d'un sommet contient toutes les arêtes y aboutissant à l'exception d'un nombre fini d'entre elles. En conséquence le groupe  $\H^{q+1}_{\underline{\bf S}}(X-{\bf S},{\cal F})$ s'identifie à $\bigoplus \limits_{Y} \H^{q+1}_{\underline{\bf S}_{|Y}}(Y,{\cal F})$ où $Y$ parcourt l'ensemble des arêtes. 

\bigskip
Il suffit dès lors de montrer que pour toute arête $Y$ et toute classe $\lambda$ appartenant à $\H^{q+1}_{\underline{\bf S}_{|Y}}(Y,{\cal F})$ l'image de $\lambda$ dans $\H^{q+1}({\cal X}\an,{\cal F})$ est trivialisée par l'analytification d'un recouvrement de Zariski de ${\cal X}$. On fixe donc $Y$ et $\lambda$ ;  l'on note $j$ l'immersion ouverte $Y\hookrightarrow {\cal X}\an$ et l'on distingue trois cas. On réutilise sans en rappeler le sens les notations des paragraphes~\ref{COH}.\ref{notlambda},~\ref{COH}.\ref{calcexpl}.\ref{ydisc}, \ref{COH}.\ref{calcexpl}.\ref{ypseucofuite} et~\ref{COH}.\ref{calcexpl}.\ref{ypseucocomp}. 

\bigskip
\setcounter{cptbis}{0}
\bigskip
\trois{comppdisc} Si $Y$ est un pseudo-disque alors d'après le~\ref{COH}.\ref{calcexpl}.\ref{ydisc} le groupe $\H^{q+1}_{\underline{\bf S}_{|Y}}(Y,{\cal F})$ est égal à $\H^{q+1}_{c}(Y,{\cal F})$, et $\lambda$ est de la forme $ h\cup \Lambda_{n,\xi}$. L'image de la classe $\lambda$ dans $\H^{q+1}({\cal X}\an,{\cal F})$ est simplement $j_{c}(\lambda)_{\bullet}$. Cette image satisfait les conclusions requises d'après le lemme~\ref{COH}.\ref{cupclasseconst}.\ref{gagafibre}. 

\bigskip
\trois{comppcour} Si $Y$ est une pseudo-couronne dont l'adhérence n'est pas compacte alors d'après le~\ref{COH}.\ref{calcexpl}.\ref{ypseucofuite} le groupe 
$\H^{q+1}_{\underline{\bf S}_{|Y}}(Y,{\cal F})$ est trivial et il n'y a donc rien à démontrer. 

\bigskip
\trois{comppcourcomp} Si $Y$ est une pseudo-couronne dont l'adhérence est compacte alors d'après le~\ref{COH}.\ref{calcexpl}.\ref{ypseucocomp} le groupe $\H^{q+1}_{\underline{\bf S}_{|Y}}(Y,{\cal F})$ est égal à  $\H^{q+1}_{c}(Y,{\cal F})$, et $\lambda$ est de la forme $h_{1}\cup \varpi_{n}+h_{2}\cup \Lambda_{n,\vec{\xi}}$. Là encore l'image de $\lambda$ dans $\H^{q+1}({\cal X}\an,{\cal F})$ est simplement $j_{c}(\lambda)_{\bullet}$. 

\bigskip
\noindent
La classe $j_{c}(h_{2}\cup \Lambda_{n,\vec{\xi}})_{\bullet}$ satisfait les conclusions requises d'après le lemme~\ref{COH}.\ref{cupclasseconst}.\ref{gagafibre}. 

\bigskip
{\bf Il reste à traiter le cas de $j_{c}(h_{1}\cup \varpi_{n})_{\bullet}$. C'est ici (et uniquement ici) qu'interviennent les hypothèses spécifiques aux assertions $i)$ et $ii)$.} 

\bigskip
\begin{itemize}
\item[]{\bf L'assertion $i)$.} La classe $\varpi_{n}$ appartient à $\H^{1}_{c}(Y,\ZZ/n)$ et provient de $\H^{1}_{c}(|Y|,\ZZ/n)$. Or $\H^{1}(|{\cal X}_{\got{c}(Y)}\an|,\ZZ/n)$ est trivial d'après l'hypothèse faite sur ${\cal X}$ et la classe $j_{c}(h_{1}\cup \varpi_{n})_{\bullet}$ est de ce fait nulle en vertu du  lemme~\ref{COH}.\ref{cupclasseconst}.\ref{toptriv}. 

\bigskip

\item[]{\bf L'assertion $ii)$.} Il suffit de traiter le cas où la classe $h_{1}$ est de la forme $\mbox{Cor}_{M/\got{c}(Y)}(h_{3}\cup (\alpha))$, où $M$ est une extension finie séparable de $\got{c}(Y)$, où $h_{3}$ appartient  à $\H^{q-1}(M,{\cal F}(-1))$ et où $\alpha$ est élément de $M^{*}$. Mais $(\alpha)\cup \varpi_{n}$ est égale à $\Lambda_{n,\vec{\alpha}}$ dans le groupe $\H^{2}_{c}(Y\times_{\got{c}(Y)}M,\mu_{n})$ avec les notations de la fin du~\ref{COH}.\ref{calcexpl}.\ref{ypseucocomp}. En conséquence $h_{1}\cup \varpi_{n}$ est égale à $\mbox{Cor}_{M/\got{c}(Y)}(h_{3}\cup \Lambda_{n,\vec{\alpha}})$. On déduit de ce fait et du lemme~\ref{COH}.\ref{cupclasseconst}.\ref{gagafibre} que la classe $j_{c}(h_{1}\cup \varpi_{n})_{\bullet}$ satisfait les conclusions requises, ce qui achève la démonstration du théorème.~$\Box$ 

\end{itemize}

\subsection*{Quelques commentaires à propos des classes décomposables}

\bigskip
\deux{merksus} Soit $K$ un corps et soit $n$ un entier inversible dans $K$.  Rost et Voevodsky ont annoncé avoir démontré que pour tout entier $q$ au moins égal à $1$ tout élément de $\H^{q}(K,\ZZ/n(q))$ était décomposable. Le cas où $q$ vaut 1 est tautologique, celui où il vaut 2 a été démontré il y a longtemps par Merkurjev et Suslin (et par Lichtenbaum encore avant dans le cas particulier où $K$ est $p$-adique). 

\bigskip
\deux{casvaldisc} Supposons que $K$ est un corps complet pour une valeur absolue ultramétrique $|\;|$ telle que $|K^{*}|$ soit libre de rang 1. Soit $n$ un entier inversible dans $\red{K}$ et soit $\cal F$ un faisceau étale sur $K$ donné par un module galoisien fini annulé par $n$. Soit $\pi$ une uniformisante de $K$. 

\setcounter{cptbis}{0}

\bigskip
\trois{cdoublet} Si $\red{K}$ est séparablement clos alors pour tout entier $q$ toute classe de $\H^{q}(K,{\cal F})$ est décomposable ; c'est une reformulation de la dualité pour la cohomologie de $K$ (compte-tenu du fait que les racines $n$-ièmes de l'unité sont dans $K$). 

\bigskip
\trois{caspadique} Supposons maintenant que $\red{K}$ est fini. Tout élément de $\H^{2}(K,{\cal F})$ est alors décomposable. En effet si $\cal F$ est constant et si les racines $n$-ièmes de l'unité sont dans le corps de base c'est simplement le résultat de Lichtenbaum mentionné ci-dessus. Dans le cas général on déduit l'assertion de la surjectivité de la corestriction, elle-même conséquence de la dualité de Tate et de l'injectivité des restrictions au niveau du $\H^{0}$. {\em Si de plus $\cal F$ est non ramifié} alors $\H^{1}(K,{\cal F})$ est lui aussi constitué de classes décomposables. Pour le voir on remarque tout d'abord que tout élément de $\H^{1}(K,{\cal F})$ est dans ces conditions de la forme $h+(\pi)\cup \eta$ où $h$ provient de $\H^{1}(K \zero,{\cal F})\simeq \H^{1}(\red{K},{\cal F})$, où $\eta$ appartient à $\H^{0}(K,{\cal F}(-1))$ et où $(\pi)$ désigne l'image de $\pi$ dans $\H^{1}(K,\mu_{n})$. Il suffit donc de montrer toute classe de $\H^{1}(\red{K},{\cal F})$ est décomposable. On procède dès lors comme ci-dessus en se ramenant au cas constant et en utilisant la dualité, sur les corps finis cette fois. 

\bigskip
\deux{cyclotom} On désigne toujours par $k$ un corps ultramétrique complet et l'on note $p$ l'exposant caractéristique de $\red{k}$. Soit $l$ un nombre premier différent de $p$ et soit $K$ un corps de décomposition de $X^{l}-1$ sur $k$. L'entier $[K:k]$ est premier à $l$. Soit ${\cal F}$ un faisceau sur $k$ donné par un module galoisien fini et de $l$-torsion et soit $\cal X$ une courbe algébrique lisse sur $k$. Soit $q$ un entier. Supposons que $\H^{0}(|{\cal X}_{K}|,{\rm R}^{q+1}\pi_{*}{\cal F})\to  \H^{0}(|{\cal X}_{K}\an|,{\rm R}^{q+1}\pi_{*}{\cal F})$ soit un isomorphisme ; alors par un argument immédiat de restriction-corestriction le noyau de $\H^{0}(|{\cal X}|,{\rm R}^{q+1}\pi_{*}{\cal F})\to  \H^{0}(|{\cal X}\an|,{\rm R}^{q+1}\pi_{*}{\cal F})$ est formé de classes triviales au point générique de $\cal X$, donc nulles par pureté. Cette flèche est en conséquence un isomorphisme, et il en va alors de même de $$\H^{1}(|{\cal X}|,{\rm R}^{q}\pi_{*}{\cal F})\to \H^{1}(|{\cal X}\an|,{\rm R}^{q}\pi_{*}{\cal F}).$$

\bigskip
Récapitulons ce qui précède sous forme d'un corollaire au théorème~\ref{COMP}.\ref{theocomp}. 

\bigskip
\deux{corms} {\bf Corollaire.}  {\em Soit ${\cal X}$ une courbe algébrique lisse sur un corps ultramétrique complet $k$, soit $q$ un entier et soit $n$ un entier premier à la caractéristique résiduelle de $k$. Soit $\cal F$ un faisceau étale sur $k$ donné par un module galoisien fini annulé par $n$. Les flèches $$\H^{0}(|{\cal X}|,{\rm R}^{q+1}\pi_{*}{\cal F})\to  \H^{0}(|{\cal X}\an|,{\rm R}^{q+1}\pi_{*}{\cal F})$$ et $\H^{1}(|{\cal X}|,{\rm R}^{q}\pi_{*}{\cal F})\to \H^{1}(|{\cal X}\an|,{\rm R}^{q}\pi_{*}{\cal F})$ sont alors des isomorphismes dans chacun des cas suivants : }

\bigskip
\begin{itemize}
\itb si $|k^{*}|$ est libre de rang 1 et si $\red{k}$ est séparablement clos ; 

\bigskip
\itb si $|k^{*}|$ est libre de rang 1, si $\red{k }$ est fini et si $q$ est au moins égal à $2$ ;

\bigskip
\itb si $|k^{*}|$ est libre de rang 1, si $\red{k }$ est fini, si $\cal F$ est non ramifié et si $q$ est égal à $1$ ; 

\bigskip
\itb si $q$ est égal à $1$ ou $2$ et si ${\cal F}$ est égal à $\ZZ/n(q)$ ; 

\bigskip
\itb si $q$ est supérieur ou égal à $1$ et si ${\cal F}$ est égal à $\ZZ/n(q)$, {\em modulo les résultats de Rost et Voevodsky} ; 

\bigskip
\itb si $n$ est premier, si $q$ est égal à $1$ ou $2$ et si $\cal F$ est de la forme $\ZZ/n(i)$ pour un certain $i$ ; 

\bigskip
\itb si $n$ est premier, si $q$ est supérieur ou égal à $1$ et si $\cal F$ est de la forme $\ZZ/n(i)$ pour un certain $i$, {\em modulo les résultats de Rost et Voevodsky}.~$\Box$

\end{itemize} }

\section{Exemple : le cas des courbes de Tate} \label{TATE}

\setcounter{cpt}{0}{\em On désigne toujours par $k$ un corps ultramétrique complet, et par $p$ l'exposant caractéristique de $\red{k}$. On suppose que $|k^{*}|$ n'est pas triviale.}

\bigskip
\deux{courbetate} Soit $\alpha$ un élément de $k^{*}$ de valeur absolue strictement inférieure à 1. Le quotient de $\gm\an$ par le sous-groupe $\alpha^{\ZZ}$ est un groupe $k$-analytique compact, sans bord et de dimension $1$ qui est isomorphe à l'analytification ${\cal X}\an$ d'une $k$-courbe elliptique ${\cal X}$  à réduction multiplicative. Notons $t$ la fonction coordonnée canonique sur $\gm$. Alors $r\mapsto (\sum a_{i}t^{i}\mapsto \max |a_{i}|r^{i})$ établit un homéomorphisme entre $\RR^{*}_{+}$ et $S(\gm\an)$ ; cet homéomorphisme en induit un entre le cercle $\RR^{*}_{+}/|\alpha|^{\ZZ}$ et $S({\cal X\an})$.

\bigskip
\deux{pointrattate} {\bf Remarque.} Le groupe $\H^{1}(L,\ZZ)$ est nul pour toute extension $L$ de $k$ ; on en déduit que ${\cal X}(L)$ est isomorphe à $L^{*}/(\alpha)^{\ZZ}$ pour tout corps ultramétrique complet $L$ au-dessus de $k$. 

\bigskip
\deux{produitinfini} Le produit infini $$(1-t)\prod_{n\geq 1}(1-\alpha^{n}t)\prod_{n\geq 1}(1-(\alpha^{n}/t))$$ est convergent dans l'algèbre de Fréchet des fonctions $k$-analytiques sur $\gm\an$ et définit donc une telle fonction que l'on notera $\Delta$. Elle satisfait l'équation fonctionnelle $\Delta(\alpha t)=(-1/t)\Delta(t)$.  Ses zéros sont les puissances de $\alpha$ et sont tous simples.

\bigskip
\deux{tatetri} Soit $P$ le point de $S(\gm\an)$ tel que $|t(P)|$  soit égal à 1, et soit $Q$ l'image de $P$ sur ${\cal X}\an$. Le singleton $\{P\}$ constitue une triangulation de $\gm\an$. Les couronnes de $\gm\an$ respectivement définies par les inégalités $|t|>1$ et $|t|<1$ seront notées $Y_{+}$ et $Y_{-}$. Ce sont deux arêtes de la triangulation $\{P\}$, les autres étant des pseudo-disques. Le complémentaire de $Y_{+}\cup Y_{-}$ dans $\gm\an$ est le domaine $k$-affinoïde $\gmu$ défini par l'égalité $|t|=1$. Par construction $\{P\}$ est une triangulation de $\gmu$. 

\bigskip
L'uniformisation $\gm\an\to {\cal X}\an$  établit un isomorphisme entre $\gmu$ et un domaine $k$-affinoïde $T$ de ${\cal X}\an$. Le complémentaire $Y$ de $T$ est une couronne dont les deux bouts adhèrent à $Q$. Le singleton $\{Q\}$ constitue une triangulation de ${\cal X}\an$. Sa restriction à $T$ correspond {\em via} l'isomorphisme mentionné ci-dessus à la triangulation $\{P\}$ de $\gmu$. En particulier toutes les arêtes de $\{Q\}$ incluses dans $T$ sont des pseudo-disques et $\{Q\}$ possède une et une seule arête en dehors de $T$, à savoir la couronne $Y$. 

\bigskip
\deux{orientate} Orientons $S(\gm\an)$ dans le sens où $|t|$ {\em décroît}. On en déduit une orientation de $S({\cal X}\an)$, et partant de $S(Y)$. On note $\omega$ (resp. $\varpi$) l'élément de $\H^{1}({\cal X}\an,\ZZ)$ (resp. $\H^{1}_{c}(Y,\ZZ)$) qui provient de l'image réciproque par la rétraction de la classe fondamentale topologique de la courbe lisse orientée $S({\cal X}\an)$ (resp. $S(Y)$). Pour tout entier $n$ on définit de même $\omega_{n}$ et $\varpi_{n}$ en remplaçant $\ZZ$ par $\ZZ/n$. Si $j$ désigne l'immersion $Y\hookrightarrow {\cal X}\an$ on a $j_{c}(\varpi)=\omega$ et $j_{c}(\varpi_{n})=\omega_{n}$ pour tout $n$. 

\bigskip
\deux{pictate} La flèche $\gm\an\to {\cal X}\an$ en induit une de $k^{*}$ vers ${\cal X}(k)$, dont le noyau est $\alpha^{\ZZ}$. Comme ${\cal X}$ est une courbe elliptique ${\cal X}(k)$ s'identifie naturellement à $\mbox{Pic}^{0}{\cal X}$, ce qui fournit un morphisme de groupes $$k^{*}\to \H^{1}({\cal X},\gm)\simeq \H^{1}({\cal X}\an, \gm)$$ (ce dernier isomorphisme découlant de GAGA), lequel se décrit comme suit : si $u$ est un élément de $k^{* }$ son image dans $\H^{1}({\cal X}\an,\gm)$ est égale à $u\cup \omega$, soit encore à $j_{c}(u\cup\varpi)$. 

\bigskip
Plus géométriquement le fibré correspondant $\EE$ se décrit comme suit. Soit $V$ un voisinage de $T$ coupant $Y$ selon deux composantes connexes $W_{+}$ et $W_{-}$, les notations étant choisies de sorte que l'orientation de $S(Y)$ aille du $(-)$ vers le $(+)$. On peut alors construire $\EE$ en recollant ${\cal O}_{V}$ et ${\cal O}_{Y}$ {\em via} l'identification de la section $1$ de ${\cal O}_{V}$ et de la section $(1,u)$ de ${\cal O}_{W_{-}\cup W_{+}}$. 

\bigskip
Pour justifier cette assertion on remarque que l'image réciproque de $\EE$ sur $\gm\an$ est le fibré trivial sur lequel $\ZZ$ agit par multiplication par $\alpha$ sur la base et par $u$ sur les fibres. Dès lors il existe une bijection naturelle entre l'ensemble des sections méromorphes de $\EE$ et l'ensemble des fonctions méromorphes $f$ sur $\gm\an$ satisfaisant l'équation fonctionnelle $f(\alpha t)=uf(t)$. Or $\Delta(t/u)/\Delta(t)$ est une telle fonction ; si $u$ est une puissance de $\alpha$ elle est inversible, et $\EE$ est trivial. Sinon ses zéros sont les $u\alpha^{n}$ pour $n$ parcourant $\ZZ$, et sont simples ; ses pôles sont les $\alpha^{n}$ pour $n$ parcourant $\ZZ$, et sont également simples. Le diviseur de la section correspondante de $\EE$ est donc égal à $P_{u}-O$ où $P_{u}$ est l'image de $u$ dans ${\cal X}(k)$, et $O$ le neutre (ou encore l'image de $1$). Ceci permet de conclure. 

\bigskip
\deux{cuptate} Soit $\cal F$ un faisceau sur $k$ donné par un module galoisien fini de cardinal premier à $p$. Soit $n$ un entier annulant $\cal F$. Soit $q$ un entier et soit $h$ appartenant à $\H^{q}(Y,{\cal F})$. Le~\ref{COH}.\ref{introcohsup} et les calculs explicites faits aux paragraphes~\ref{COH}.\ref{calcexpl}.\ref{ydisc}, \ref{COH}.\ref{calcexpl}.\ref{ypseucofuite} et~\ref{COH}.\ref{calcexpl}.\ref{ypseucocomp} fournissent alors une condition nécessaire et suffisante pour que la classe $j_{c}(h\cup \varpi_{n})$ soit nulle. C'est plus précisément le cas si et seulement si il existe un élément $\lambda$ appartenant à $\H^{q}(\hres(Q),{\cal F})$ (que l'on peut prolonger en une classe définie au voisinage de $Q$, classe que l'on notera encore $\lambda$) et qui vérifie : 

\bigskip
\begin{itemize} 
\item [$i)$] pour toute arête $Z$ incluse dans $T$ il existe une sous-pseudo-couronne de $Z$ sur laquelle $\lambda$ est définie et constante (c'est-à-dire provient de $\got{c}(Z)$) ; 

\bigskip
\item[$ii)$]  il existe deux sous-couronnes disjointes $W_{+}$ et $W_{-}$ de $Y$, chacune située à un bout (les signes étant toujours choisis de manière compatible avec l'orientation) et incluse dans l'ouvert de définition de $\lambda$, et telles que $\lambda_{|W_{+}}-\lambda_{|W_{-}}$ soit égale à $h$, la différence étant calculée dans $\H^{q}(Y,{\cal F})$ {\em via} les isomorphismes $\H^{q}(W_{+},{\cal F})\simeq  \H^{q}(Y,{\cal F})$ et $\H^{q}(W_{-},{\cal F})\simeq  \H^{q}(Y,{\cal F})$. 

\end{itemize} 

\bigskip
\deux{remontegm} Considérons une classe $\lambda$ comme ci-dessus. L'uniformisation $\gm\an\to{\cal X}\an$ induit un isomorphisme entre $\gmu$ et $T$ qui respecte les triangulations considérées ; la classe $\lambda$ peut donc être vue comme un élément de $\H^{q}(\hres(P),{\cal F})$, que l'on peut prolonger au voisinage de $P$ dans $\gm\an$ en une classe notée encore $\lambda$ et qui possède la propriété suivante : {\em pour toute arête $Z$ incluse dans $\gmu$ il existe une sous-pseudo-couronne de $Z$ sur laquelle $\lambda$ est définie et constante (c'est-à-dire provient de $\got{c}(Z)$)}. 

\bigskip
Les arêtes de la triangulation $\{P\}$ sont toutes incluses dans $T$ à l'exception des couronnes $Y_{+}$ et $Y_{-}$, qui sont toutes deux d'adhérence non compacte dans $\gm\an$. On déduit alors des paragraphes~\ref{COH}.\ref{calcexpl}.\ref{ydisc} et \ref{COH}.\ref{calcexpl}.\ref{ypseucofuite} que $\lambda$ s'annule dans $\H^{q+1}_{\underline{\{P\}}}(\gm\an-\{P\},{\cal F})$ ; en conséquence elle provient d'après le~\ref{COH}.\ref{introcohsup} d'une classe de $\H^{q}(\gm\an,{\cal F})$ que l'on note encore $\lambda$. 

\bigskip
Soient $V_{+}$ et $V_{-}$ deux sous-couronnes respectives de $Y_{+}$ et $Y_{-}$, aboutissant à $P$, et telles que la flèche de $\gm\an$ dans ${\cal X}\an$ induise un isomorphisme entre $V_{+}$ et une sous-couronne ${\cal V}_{+}$ de $Y$, ainsi qu'entre  $V_{-}$ et une sous-couronne ${\cal V}_{-}$ de $Y$. Les sous-couronnes ${\cal V}_{+}$ et ${\cal V}_{-}$ sont nécessairement disjointes. Notons que $t$ définit {\em via} l'isomorphisme $V_{-}\simeq {\cal V}_{-}$ une fonction coordonnée $\tau$ sur ${\cal W}_{-}$ ; celle-ci s'étend en une fonction coordonnée $\tau$ définie sur $Y$ dont la restriction à ${\cal V}_{+}$ correspond, par le biais de $V_{+}\simeq {\cal V}_{+}$, à la fonction $\alpha t$. 

\bigskip
Comme $\gm\an$ est une couronne on peut écrire $\lambda$ sous la forme $\mu +(t)\cup \eta$ où $\mu$ appartient à $\H^{q}(k,{\cal F})$ et $\eta$ à $\H^{q-1}(k,{\cal F}(-1))$. Avec les notations du~\ref{TATE}.\ref{cuptate} on a dès lors $\lambda_{|W_{-}}=\mu +(\tau)\cup \eta$ et $\lambda_{|W_{+}}=\mu +(\alpha\tau)\cup \eta$ ; par conséquent $h$ est égale à $(\alpha)\cup \eta$. 

\bigskip
Réciproquement si $\eta$ appartient à $\H^{q-1}(k,{\cal F}(-1))$ alors $j_{c}((\alpha)\cup\eta\cup \varpi_{n})$, qui n'est autre que $(\alpha)\cup \eta \cup \omega_{n}$, est triviale : en effet $(\alpha)\cup \omega_{n}$ provient par le biais du cobord de la suite de Kummer de la classe $\alpha\cup \omega$ de $\H^{q}({\cal X}\an,\gm)$ (où $\alpha$ est vu comme appartenant à $\H^{0}({\cal X}\an,\gm)$), laquelle est d'après le~\ref{TATE}.\ref{pictate} l'image de $\alpha$ par la flèche $k^{*}\to {\cal X}(k)$. Elle est donc nulle, ce qui permet de conclure. On a finalement établi la proposition suivante : 

\bigskip
\deux{noyautate} {\bf Proposition.} {\em Avec les notations ci-dessus le noyau de la flèche $$h\mapsto j_{c}(h\cup \varpi_{n})$$ qui va de $\H^{q}(Y,{\cal F})$ dans $\H^{q+1}({\cal X}\an,{\cal F})\simeq \H^{q+1}({\cal X},{\cal F})$ est le sous-groupe formé des classes de la forme $(\alpha)\cup \eta$ où $\eta$ appartient à $\H^{q-1}(k,{\cal F}(-1))$.}~$\Box$ 

\bigskip
\deux{calchzerotate} {\bf Description du groupe $\H^{0}(|{\cal X}\an|,\mbox{\rm R}^{q}\pi_{*}{\cal F})$.} Le théorème~\ref{COH}.\ref{cohocomplexe} stipule que le groupe en question est isomorphe au noyau de $$\H^{q}(\hres(Q),{\cal F})\to \H^{q+1}_{\underline{\{Q\}}}({\cal X}\an-\{Q\},{\cal F}).$$ Soit $\lambda$ appartenant à $\H^{q}(\hres(Q),{\cal F})$. On la prolonge au voisinage de $Q$ en une classe encore notée $\lambda$ ; des remarques faites au~\ref{COH}.\ref{calcexpl}.\ref{ydisc} et~\ref{COH}.\ref{calcexpl}.\ref{ypseucocomp} on déduit que $\lambda$ appartient au noyau considéré si et seulement si elle satisfait la condition $i)$ du~\ref{TATE}.\ref{cuptate} et la condition $ii)$ du même paragraphe avec $h=0$. 

\bigskip
Voyons $\lambda$ comme un élément de $\H^{q}(\hres(P),{\cal F})$ par le biais de l'isomorphisme $\hres(Q)\simeq \hres(P)$ induit par l'uniformisation. En utilisant le théorème~\ref{COH}.\ref{cohocomplexe} ainsi que le~\ref{COH}.\ref{calcexpl}.\ref{ypseucofuite} on voit exactement comme au~\ref{TATE}.\ref{remontegm} que la condition $i)$ équivaut à dire que $\lambda$ se prolonge en une classe de $\H^{q}(\gm\an,{\cal F})$ ; une telle classe est nécessairement unique d'après la proposition~\ref{GERM}.\ref{courskel}.\ref{locglob} puisque $\gm\an$ est une couronne. Pour cette même raison cet unique prolongement est alors nécessairement de la forme $\mu +(t)\cup \eta$ où $\mu$ appartient à $\H^{q}(k,{\cal F})$ et $\eta$ à $\H^{q-1}(k,{\cal F}(-1))$. Les classes $\mu$ et $\eta$ sont déterminées sans ambiguïté. 

\bigskip
Dans cette situation la différence $\lambda_{|W_{+}}-\lambda_{|W_{-}}$ (avec les notations du~\ref{TATE}.\ref{cuptate} et par le même raisonnement qu'au~\ref{TATE}.\ref{remontegm}) est égale à $(\alpha)\cup \eta$. Compte-tenu du fait que $\H^{q}(k,{\cal F})$ s'injecte dans $\H^{q}(Y,{\cal F})$ puisque $Y$ est une couronne on vient finalement d'établir la proposition suivante : 

\bigskip
\deux{deschzerotate} {\bf Proposition.} {\em Identifions $t$ à un élément de $\hres(Q)^{*}$ {\em via} l'isomorphisme $\hres(Q)\simeq \hres(P)$ induit par l'uniformisation. Notons $_{(\alpha)}\H^{q-1}(k,{\cal F}(-1))$ le sous-groupe de $\H^{q-1}(k,{\cal F}(-1))$ formé des classes annulées par le cup-produit avec $(\alpha)$. Soit $\mu$ un élément de $\H^{q}(k,{\cal F})$ et soit $\eta$ un élément de $_{(\alpha)}\H^{q-1}(k,{\cal F}(-1))$. Alors la classe $\mu+(t)\cup \eta$ de $\H^{q}(\hres(Q),{\cal F})$ s'étend d'une unique manière en un élément de $\H^{0}(|{\cal X}\an|,\mbox{\rm R}^{q}\pi_{*}{\cal F})$ et cette construction induit un isomorphisme $$\H^{q}(k,{\cal F})\oplus \;_{(\alpha)}\H^{q-1}(k,{\cal F}(-1))\simeq \H^{0}(|{\cal X}\an|,\mbox{\rm R}^{q}\pi_{*}{\cal F}).~\Box$$}

\bigskip
\deux{pictau} On conserve toujours les mêmes notations. La classe $j_{c}(\tau\cup \varpi)$ est un élément de $\mbox{Pic}\;{\cal X}\an$, qui est isomorphe à $\mbox{Pic}\;{\cal X}$ par GAGA. L'image réciproque sur $\gm\an$ du fibré correspondant $\FF$ est le fibré trivial sur lequel $\ZZ$ agit en multipliant par $\alpha$ sur la base et par $t$ sur les fibres ; l'espace des sections méromorphes de $\FF$ est donc en bijection naturelle avec l'ensemble des fonctions méromorphes $f$ sur $\gm\an$ qui satisfont l'équation fonctionnelle $f(\alpha t)=tf(t)$. Or c'est le cas de $1/\Delta$, qui n'a pas de zéros et dont les pôles sont les puissances de $\alpha$ et sont simples. La section méromorphe correspondante de $\FF$ a donc pour diviseur $(-O)$ ; en conséquence $\FF$ est de degré $(-1)$. Comme $\mbox{Pic}\;{\cal X}$ est une extension de $\ZZ$ ({\em via} le degré) par $\mbox{Pic}\zero {\cal X}$ on déduit de ce qui précède et du~\ref{TATE}.\ref{pictate} que tout élément de $\H^{1}({\cal X},\gm)$ est de la forme $j_{c}(u\tau^{m}\cup \varpi)$ pour un certain $u$ appartenant à $k^{*}$ (et uniquement déterminé modulo $\alpha^{\ZZ}$) et un certain entier relatif $m$. 

\bigskip
\deux{classezartriv} Soit $\bf x$ un point fermé de $\cal X$ et soit $\eta$ un élément de $\H^{q-1}(\kappa({\bf x}),{\cal F}(-1))$. Notons $c({\bf x})$ l'image de $\bf x$ dans $\H^{1}_{\bf x}({\cal X},\gm)$ (obtenue en trivialisant le fibré associé au diviseur $\bf x$ par la section $1$ sur l'ouvert ${\cal X}-\{\bf x\}$). Soit $\partial$ le cobord de la suite exacte de Kummer. Le cup-produit $\partial(c({\bf x}))\cup \eta$ vit dans le groupe $\H^{q+1}_{\bf x}({\cal X},{\cal F})$ ; on note $(\partial(c({\bf x}))\cup \eta)_{\bullet}$ la classe de $\H^{q+1}({\cal X},{\cal F})$ obtenue en oubliant le support. 

Supposons que $\kappa({\bf x})$ est séparable sur $k$. Le point $\bf x$ possède un antécédent $\kappa({\bf x})$-rationnel $\bf y$ sur ${\cal X}_{\kappa({\bf x})}$. On a l'égalité $\partial(c({\bf x}))\cup \eta= \mbox{Cor}_{\kappa({\bf x})/k}(\partial(c({\bf y}))\cup \eta)$. Par ailleurs $(\partial(c({\bf y}))\cup \eta) _{\bullet}$ est égal à $(\partial(c({\bf y}))_{\bullet}\cup \eta$ puisque $\eta$ provient de la cohomologie de $\kappa({\bf x})$. 

\bigskip
La classe $c({\bf y})_{\bullet}$ peut être vue comme un élément de $\H^{1}({\cal X}_{\kappa({\bf x})}\an,\gm)$ ; elle est donc de la forme $j_{c}(u\tau^{m}\cup \varpi)$ où $u$ appartient à $\kappa({\bf x})^{*}$ et où $m$ est un entier relatif. Dès lors $\partial(c({\bf y}))_{\bullet}$ s'écrit $(u)\cup \omega_{n}+mj_{c}((\tau)\cup \varpi_{n})$. On déduit de tout ceci l'égalité $$(\partial(c({\bf x}))\cup \eta)_{\bullet}=\mbox{Cor}_{\kappa({\bf x})/k}((u)\cup \eta)\cup \omega_{n}  +m j_{c}((\tau)\cup \mbox{Cor}_{\kappa({\bf x})/k}(\eta)\cup \varpi_{n})$$ dans le groupe $\H^{q+1}({\cal X},{\cal F})$.

\bigskip
\deux{theoopt} {\bf Théorème.} {\em Désignons par $S$ le sous-groupe de $\H^{q}(k,{\cal F})$ formé des classes décomposables. L'application $h\mapsto h\cup \omega_{n}$, composée avec les flèches $$\H^{q+1}({\cal X}\an,{\cal F})\simeq \H^{q+1}({\cal X},{\cal F})\to \H^{0}(|{\cal X}|,\mbox{\rm R}^{q+1}\pi_{*}{\cal F}),$$ établit un isomorphisme entre $\H^{q}(k,{\cal F})/S$ et le noyau de $$ \H^{0}(|{\cal X}|,\mbox{\rm R}^{q+1}\pi_{*}{\cal F})\to \H^{0}(|{\cal X}\an|,\mbox{\rm R}^{q+1}\pi_{*}{\cal F}).$$}

\bigskip
{\em Démonstration.} Notons pour commencer que la classe $\omega_{n}$ provient du groupe $\H^{1}(|{\cal X}\an|,\ZZ/n)$ ; en conséquence pour toute classe $h$ l'élément $h\cup \omega_{n}$ se tue dans $\H^{0}(|{\cal X}\an|,\mbox{\rm R}^{q+1}\pi_{*}{\cal F})$. De plus $h\cup \omega_{n}$ est égale à $j_{c}(h\cup \varpi_{n})$. Si $h$ est décomposable on en déduit, par la méthode suivie lors de la démonstration de l'assertion $ii)$ du théorème~\ref{COMP}.\ref{theocomp}, que $h\cup \omega_{n}$ est localement triviale pour la topologie de Zariski, et donc est nulle si on la voit dans $\H^{0}(|{\cal X}|,\mbox{\rm R}^{q+1}\pi_{*}{\cal F})$ ; on a bien ainsi défini un morphisme de $\H^{q+1}(k,{\cal F})/S$ vers le noyau de la flèche $ \H^{0}(|{\cal X}|,\mbox{\rm R}^{q+1}\pi_{*}{\cal F})\to \H^{0}(|{\cal X}\an|,\mbox{\rm R}^{q+1}\pi_{*}{\cal F})$. 

\bigskip
\setcounter{cptbis}{0}
\trois{flechesur} {\em Le morphisme est surjectif.} Considérons un élément $\lambda$ appartenant au noyau de $ \H^{0}(|{\cal X}|,\mbox{\rm R}^{q+1}\pi_{*}{\cal F})\to \H^{0}(|{\cal X}\an|,\mbox{\rm R}^{q+1}\pi_{*}{\cal F})$. La classe $\lambda$ provient d'une classe de $\H^{q+1}({\cal X},\gm)\simeq \H^{q+1}({\cal X}\an,\gm)$ que l'on notera encore $\lambda$. Par hypothèse $\lambda$ s'annule dans $\H^{0}(|{\cal X}\an|,\mbox{\rm R}^{q+1}\pi_{*}{\cal F})$. Du~\ref{COH}.\ref{introcohsup} et du lemme~\ref{COH}.\ref{nulsurs} on déduit que $\lambda$ provient du groupe $\H^{q+1}_{\underline{\{Q\}}}({\cal X}\an-\{Q\},{\cal F})$. On a vu lors de la preuve du théorème~\ref{COMP}.\ref{theocomp} que $\H^{q+1}_{\underline{\{Q\}}}({\cal X}\an-\{Q\},{\cal F})$ était isomorphe à la somme directe des $\H^{q+1}_{\underline{\{Q\}}_{|Z}}(Z,{\cal F})$ où $Z$ parcourt l'ensemble des arêtes de la triangulation $\{Q\}$. Comme elles sont toutes d'adhérence compacte on a en fait affaire à la somme directe des  $\H^{q+1}_{c}(Z,{\cal F})$. 

\bigskip
Si $Z$ est un pseudo-disque alors l'image de tout élément de $\H^{q+1}_{c}(Z,{\cal F})$ dans $\H^{l+1}({\cal X}\an,{\cal F})$ est localement triviale pour la topologie de Zariski de $\cal X$ d'après le~\ref{COMP}.\ref{theocomp}.\ref{comppdisc}, donc nulle une fois vue dans $ \H^{0}(|{\cal X}|,\mbox{\rm R}^{q+1}\pi_{*}{\cal F})$. 

\bigskip
Il reste à traiter le cas de l'arête $Y$. Tout élément de $\H^{q+1}_{c}(Y,{\cal F})$ est de la forme $(\mu +(\tau)\cup \eta)\cup \varpi_{n}$ où $\mu$ et $\eta$ appartiennent respectivement à $\H^{q}(k,{\cal F})$ et $\H^{q-1}(k,{\cal F}(-1))$. En vertu du~\ref{COMP}.\ref{theocomp}.\ref{comppcourcomp} l'image de $((\tau)\cup \eta)\cup \varpi_{n}$ dans $\H^{q+1}({\cal X}\an,{\cal F})$ est localement triviale pour la topologie de Zariski de $\cal X$ et donc est nulle dans $ \H^{0}(|{\cal X}|,\mbox{\rm R}^{q+1}\pi_{*}{\cal F})$. On peut finalement supposer que $\lambda$ s'écrit $j_{c}(\mu\cup \varpi_{n})$, c'est-à-dire encore $\mu\cup \omega_{n}$. Ceci permet de conclure à la surjectivité. 

\bigskip
\trois{flecheinj} {\em Le morphisme est injectif}. Soit $h$ appartenant à $\H^{q}(k,{\cal F})$ telle que l'élément $h\cup \omega_{n}$ de $\H^{q+1}({\cal X}\an,{\cal F})\simeq \H^{q+1}({\cal X},{\cal F})$  s'annule dans le groupe $\H^{0}(|{\cal X}|,\mbox{\rm R}^{q+1}\pi_{*}{\cal F})$. Alors cette classe est en particulier trivialisée sur un ouvert de Zariski non vide de $\cal X$ ; quitte à faire une extension finie purement inséparable de $k$ (ce qui est sans effet sur ce qu'on étudie) on peut supposer que l'ouvert en question est le complémentaire d'un ensemble fini $\{{\bf x}_{1},\ldots,{\bf x}_{r}\}$ de points fermés à corps résiduels séparables sur $k$. La suite exacte de cohomologie à support combinée avec les théorèmes de pureté en cohomologie étale assure alors que $h\cup \omega_{n}$ s'écrit $\sum (\partial(c({\bf x}_{i}))\cup \eta_{i})_{\bullet}$, avec les notations du~\ref{TATE}.\ref{classezartriv}, où $\eta_{i}$ appartient pour tout $i$ à $\H^{q-1}(\kappa({\bf x}_{i}),{\cal F}(-1))$. 

\bigskip
Des calculs faits au~\ref{TATE}.\ref{classezartriv} on déduit que $h\cup \omega_{n}$, qui n'est autre que $j_{c}(h\cup \varpi_{n})$, s'écrit $$\sum j_{c}\left( \left[\mbox{Cor}_{\kappa({\bf x}_{i})/k}((u_{i})\cup \eta_{i}) +m_{i}(\tau)\cup \mbox{Cor}_{\kappa({\bf x}_{i})/k}(\eta_{i})\right]\cup \varpi_{n}\right)$$ où, $i$ étant fixé, $u_{i}$ appartient à $\kappa({\bf x}_{i})^{*}$ et  $m_{i}$ à $\ZZ$. 

\bigskip
De la proposition~\ref{TATE}.\ref{noyautate} on déduit que l'image de $h$ dans $\H^{q}(Y,{\cal F})$ est égale à $$\sum \mbox{Cor}_{\kappa({\bf x}_{i})/k}((u_{i})\cup \eta_{i}) +m_{i}(\tau)\cup \mbox{Cor}_{\kappa({\bf x}_{i})/k}(\eta_{i})$$ modulo une classe qui est de la forme $(\alpha)\cup \beta$ où $\beta$ appartient à $\H^{q-1}(k,{\cal F}(-1))$. Quitte à retrancher cette dernière (qui est décomposable) à $h$ on peut supposer que la restriction de $h$ à $\H^{q}(Y,{\cal F})$  est exactement $$\sum \mbox{Cor}_{\kappa({\bf x}_{i})/k}((u_{i})\cup \eta_{i}) +m_{i}(\tau)\cup \mbox{Cor}_{\kappa({\bf x}_{i})/k}(\eta_{i}).$$ Comme $Y$ est une couronne $\H^{q}(k,{\cal F})$ s'injecte dans $\H^{q}(Y,{\cal F})$ et son image est en somme directe avec le groupe des classes de la forme $(\tau)\cup \eta$ ; dès lors $h$ est égale à $\sum \mbox{Cor}_{\kappa({\bf x}_{i})/k}((u_{i})\cup \eta_{i})$ et appartient en conséquence à $S$, ce qui achève la démonstration.~$\Box$

\bigskip
\deux{remopt} {\bf Remarque.} Le théorème ci-dessus montre que les conditions suffisantes $i)$ et $ii)$ du théorème~\ref{COMP}.\ref{theocomp} sont optimales. 

\section{Classes gratte-ciel et contre-exemples au principe de Hasse}\label{GRAT}

\setcounter{cpt}{0} 
{\em On désigne à nouveau par $k$ un corps ultramétrique complet {\em quelconque} (la valeur absolue peut être triviale). L'exposant caractéristique de $\red{k}$ est toujours noté $p$. La notion de {\em faisceau raisonnable} est celle de la définition~\ref{COH}.\ref{defmodfasc}.}

\bigskip
\deux{grattecour} {\bf Lemme.} {\em Soit $X$ une $k$-courbe qui est ou bien un pseudo-disque, ou bien une pseudo-couronne. Soit ${\cal F}$ un faisceau raisonnable sur $X$, soit $q$ un entier et soit $h$ appartenant à $\H^{q}(X,{\cal F})$. Supposons qu'il existe un ensemble fermé et discret $E$ de $X$ tel que $h(P)$ soit nulle pour tout $P$ n'appartenant pas à $E$. Alors $h$ est nulle.} 

\bigskip
{\em Démonstration.} Si $X$ est une pseudo-couronne alors il existe un point $P$ sur $S(X)$ qui n'est pas dans $E$, et le résultat se déduit de la proposition~\ref{GERM}.\ref{courskel}.\ref{locglob}. Si $X$ est un pseudo-disque il possède une sous-pseudo-couronne $Z$ et il existe un point de $S(Z)$ qui n'est pas dans $E$. On conclut à l'aide du~\ref{GERM}.\ref{loczeropseudisc}.~$\Box$

\bigskip
\deux{intrograt} {\bf Définition}. {\em Soit $X$ un espace $k$-analytique, soit $\cal F$ un faisceau sur $X$ et soit $q$ un entier. On dira qu'un élément $h$ de $\H^{0}(|X|,\mbox{\rm R}^{q}\pi_{*}{\cal F})$ est une {\em classe gratte-ciel} s'il existe un sous-ensemble fermé et discret $E$ de $X$ tel que la restriction de $h$ à $X-E$ soit triviale. Le sous-groupe de $\H^{0}(|X|,\mbox{\rm R}^{q}\pi_{*}{\cal F})$ formé des classes gratte-ciel sera noté $\G^{q}(X,{\cal F})$. On désignera par $\G^{q}_{c}(X,{\cal F})$ le groupe des classes gratte-ciel à support compact.}

\bigskip
\deux{grattepc} Soit $X$ une $k$-courbe qui est ou bien un pseudo-disque ou bien une pseudo-couronne. Si $\cal F$ est un faisceau raisonnable sur $X$ alors pour tout $q$ les groupes $\G^{q}(X,{\cal F})$ et $\G^{q}_{c}(X,{\cal F})$ sont nuls : en effet $\H^{0}(|X|,\mbox{\rm R}^{q}\pi_{*}{\cal F})$ s'identifie à $\H^{q}(X,{\cal F})$ d'après la proposition~\ref{GERM}.\ref{aretetriv}, et il n'y a plus alors qu'à appliquer le lemme~\ref{GRAT}.\ref{grattecour}. 

\bigskip
\deux{gratteloc} Soit $X$ un espace $k$-analytique et $\cal F$ un faisceau sur $X$. Soit $P$ un point de $X$. Soit $q$ un entier. On notera $\G_{P}^{q}(X,{\cal F})$ le sous-groupe de $\H^{q}(\hres(P),{\cal F}_{P})$ formé des classes qui s'annulent sur un voisinage épointé de $P$. Si $E$ est un sous-ensemble fermé et discret de $X$ le sous-groupe de $\G^{q}(X,{\cal F})$ (resp. de $\G^{q}_{c}(X,{\cal F})$) formé des classes à support inclus dans $E$ est naturellement isomorphe à $\prod\limits_{P\in E} \G_{P}^{q}(X,{\cal F})$ (resp.$\bigoplus\limits_{P\in E} \G_{P}^{q}(X,{\cal F})$). 

\setcounter{cptbis}{0} 
\bigskip
\trois{equivgratzero} On a donc équivalence entre les trois propositions : 

\bigskip
\begin{itemize}
\itb $\G^{q}(X,{\cal F})$ est nul. 
\itb $\G^{q}_{c}(X,{\cal F})$ est nul. 
\itb Pour tout point $P$ de $X$ le groupe $\G^{q}_{P}(X,{\cal F})$ est trivial. 
\end{itemize} 

\bigskip
\trois{gpcour} En particulier si $X$ est un pseudo-disque ou une pseudo-couronne et si $\cal F$ est raisonnable alors $\G^{q}_{P}(X,{\cal F})$ est nul pour tout point $P$ de $X$. 

\bigskip
\deux{trigratte} {\bf Lemme.} {\em Soit $X$ une $k$-courbe quasi-lisse et soit $P$ un point de $X$ ou bien de type $(2)$, ou bien de type $(3)$. Soit $\cal F$ un faisceau raisonnable sur $X$, soit $q$ un entier, et soit $h$ appartenant au groupe $\H^{q}(\hres(P),{\cal F}_{P})$. Alors $h$ appartient à $\G^{q}_{P}(X,{\cal F})$ si et seulement si la restriction de $h$ à $\H^{q}({\cal O}({\cal E})_{\rm alg},{\cal F})$ (restriction définie  {\em via} l'isomorphisme entre la cohomologie de ${\cal O}_{X,P}$ et celle de $\hres(P)$ ) est triviale pour tout système de composantes $\cal E$ adhérent à $P$.}

\bigskip
{\em Démonstration.} La classe $h$ s'étend en un élément de $\H^{q}(V,{\cal F})$ pour un certain voisinage $V$ de $P$. Supposons que $h$ appartienne à $\G^{q}_{P}(X,{\cal F})$. On peut alors restreindre $V$ de sorte que $h_{|V-\{P\}}$ soit nulle. Sous cette hypothèse  la classe $h_{|{\cal E}_{V}}$ est {\em a fortiori} triviale, et comme $\H^{q}({\cal O}({\cal E})_{\rm alg},{\cal F})$ est naturellement isomorphe à la limite inductive des $\H^{q}({\cal E}_{V},{\cal F})$ d'après le lemme~\ref{GERM}.\ref{equivetale} le résultat cherché s'ensuit aussitôt. 

\bigskip
Pour la réciproque la proposition~\ref{TRI}.\ref{triexist} permet de supposer que $P$ est un sommet d'une triangulation de $V$, et donc de se ramener au cas où $V-\{P\}$ est une somme disjointe de pseudo-disques et de pseudo-couronnes aboutissant à $P$. Soit $Y$ une composante connexe de $V-\{P\}$ et soit $\cal E$ le système de composantes correspondant (ou l'un des deux si $Y$ est une pseudo-couronne telle que $Y\cup\{P\}$ soit un cercle). Comme la restriction de $h$ à ${\cal O}({\cal E})_{\rm alg}$ est triviale le lemme~\ref{GERM}.\ref{equivetale} assure qu'il existe une sous-pseudo-couronne $Z$ de $Y$ aboutissant à $P$ telle que $h_{|Y}$ soit nulle. Or l'application de restriction $$\H^{q}(Y,{\cal F})\to \H^{q}(Z,{\cal F})$$ est injective : si $Y$ est-elle même une pseudo-couronne c'est la remarque~\ref{PSEU}.\ref{remsouspseu}. Si $Y$ est un pseudo-disque alors la restriction de $\cal F$ à $Y$ provient d'un faisceau raisonnable sur $\got{c}(Y)$ par le~\ref{PSEU}.\ref{cohogen}, faisceau que l'on note encore $\cal F$. D'autre part $\got{c}(Z)$ est égal à $\got{c}(Y)$ et $\H^{q}(\got{c}(Z),{\cal F})\to \H^{q}(Z,{\cal F})$ est injective d'après le~\ref{GERM}.\ref{constcour}. Finalement $h_{|Y}$ est nulle, et ce pour toute composante connexe $Y$ de $V-\{P\}$. On en déduit que $h$ appartient à $\G^{q}(X,P)$.~$\Box$ 

\bigskip
\deux{gratterig}  {\bf Lemme.} {\em Soit $X$ une $k$-courbe quasi-lisse. Soit $P$ un point de $X$ qui est ou bien rigide, ou bien de type $(3)$ et non isolé. Soit $\cal F$ un faisceau raisonnable sur $X$ et soit $q$ un entier. Le groupe $\G^{q}_{P}(X,{\cal F})$ est nul.} 

\bigskip
{\em Démonstration.} Si $P$ est rigide il possède alors un voisinage qui est un pseudo-disque, et donc $\G^{q}_{P}(X,{\cal F})$ est nul d'après le~\ref{GRAT}.\ref{gratteloc}.\ref{gpcour} ci-dessus. 

\bigskip
Soit $P$ un point de type $(3)$ non isolé dans $X$ et soit $h$ un élément du groupe $\H^{q}({\cal O}_{X,P},{\cal F})$ qui appartient à $\G^{q}_{P}(X,{\cal F})$. Comme $P$ est non isolé il y a au moins un système de composantes $\cal E$ de $X$ adhérent à $P$. Le lemme~\ref{GRAT}.\ref{trigratte} assure que la restriction de $h$ à ${\cal O}({\cal E})_{\rm alg}$ est triviale. Comme ${\cal O}({\cal E})_{\rm alg}$ s'identifie à ${\cal O}_{X,P}$ en vertu du~\ref{GERM}.\ref{henslocal}.\ref{oealgtype3} la classe $h$ est elle-même nulle, ce qui termine la démonstration.~$\Box$ 

\bigskip
\deux{ouvtri} {\bf Proposition.} {\em Soit $X$ une $k$-courbe quasi-lisse sans point de type $(3)$ isolé. Soit $\bf S$ une triangulation de $X$ et soit ${\bf S}_{(2)}$ l'ensemble des sommets de type $(2)$. Soit $\cal F$ un faisceau raisonnable sur $X$ et soit $U$ un ouvert de $X$. Soit $q$ un entier. Le support de toute classe gratte-ciel de $\H^{0}(|U|,\mbox{\rm R}^{q}\pi_{*}{\cal F})$ est alors contenu dans ${\bf S}_{(2)}\cap U$ ; autrement dit  $\G^{q}(U,{\cal F})$ (resp. $\G_{c}^{q}(U,{\cal F})$ ) est naturellement isomorphe à $$\prod_{P\in U\cap {\bf S}_{(2)}} \G^{q}_{P}(X,{\cal F})\;(\mbox{resp.}\;\bigoplus_{P\in U\cap {\bf S}_{(2)}} \G^{q}_{P}(X,{\cal F})\;).$$}

\bigskip
{\em Démonstration.} Soit $P$ un point de $U$. Notons que comme $U$ est un ouvert de $X$ les groupes $\G^{q}_{P}(U,{\cal F})$ et $\G^{q}_{P}(X,{\cal F})$ coïncident. Si $P$ n'est pas dans $\bf S$ il possède dans $X$ un voisinage qui est un pseudo-disque ou une pseudo-couronne, et en conséquence $\G^{q}_{P}(U,{\cal F})$ est nul d'après le~\ref{GRAT}.\ref{gratteloc}.\ref{gpcour}. Si $P$ est dans ${\bf S}-{\bf S}_{(2)}$ alors $\G^{q}_{P}(U,{\cal F})$ est trivial en vertu du lemme~\ref{GRAT}.\ref{gratterig}. Le résultat cherché découle de ces faits et du~\ref{GRAT}.\ref{gratteloc}.~$\Box$ 

\bigskip
{\em A partir de maintenant on fixe un nombre premier $l$ différent de $p$. On désigne par $\cal F$ un $\red{k}$-module galoisien fini annulé par $l^{n}$ pour un certain entier $n$. On peut voir $\cal F$ comme un faisceau étale sur $\red{k}$ aussi bien que sur $k$. On se donne une $k$-courbe quasi-lisse $X$.} 

\bigskip
\deux{rgdiv} Soit $I$ un ensemble tel que la partie $l$-primaire de $\sqrt{|k^{*}|}/|k^{*}|$ soit isomorphe à $(\QQ_{l}/\ZZ_{l})^{(I)}$. Si $P$ est un point de type $(2)$ de $X$ alors $|k^{*}|$ est d'indice fini  dans $|\hres(P)^{*}|$. La partie $l$-primaire de $\sqrt{|\hres(P)^{*}|}/|\hres(P)^{*}|$ est en conséquence elle aussi isomorphe à $(\QQ_{l}/\ZZ_{l})^{(I)}$. Si $r$ est un entier on notera ${I\choose r}$ l'ensemble des parties de $I$ à $r$ éléments. 

\bigskip
\deux{rgdivindep} {\bf Remarque.} Seul importe le {\em cardinal} de $I$. Il dépend {\em a priori} du nombre premier $l$ fixé ; toutefois si le groupe abélien $|k^{*}|$ est abstraitement isomorphe à une somme directe $\ZZ^{(B)}\oplus D$ où $B$ est un ensemble et $D$ un $\QQ$-espace vectoriel alors $|I|$ est égal à  $|B|$ et ce quelque soit $l$. 

\bigskip
\deux{courberes} Soit $P$ un point de type $(2)$ de $X$. Le corps résiduel $\red{\hres(P)}$ est le corps des fonctions d'une courbe projective, normale et géométriquement intègre sur une certaine extension finie de $\red{k}$, et que l'on appellera la {\em courbe résiduelle}. Si $\cal Q$ est un point fermé de ladite courbe (ou encore une $\red{k}$-valuation discrète de $\red{\hres(P)}$) on désignera par $\red{\hres(P)}_{\cal Q}$ le hensélisé correspondant. Rappelons que l'ensemble des systèmes de composantes adhérents à $P$ est en bijection naturelle avec l'ensemble des points fermés d'un ouvert de Zariski non vide de la courbe résiduelle, ouvert que l'on notera ${\cal U}_{X,P}$ et qui est égal à la courbe toute entière dans le cas où $P$ est un point intérieur de $X$. La proposition ci-dessous est une conséquence immédiate du~\ref{GRAT}.\ref{rgdiv} et  des lemmes~\ref{HENS}.\ref{decompocup} et~\ref{HENS}.\ref{hasseloc}. 

\bigskip
\deux{hqtype2} {\bf Proposition.} {\em On conserve les notations introduites ci-dessus. Soit $q$ un entier et soit $P$ un point de type $(2)$ de $X$. Il existe un isomorphisme $$\H^{q}(\hres(P),{\cal F})\simeq \bigoplus_{r\leq \min(q,|I|)} \left[\;\H^{q-r}(\red{\hres(P)},{\cal F}(-r))\;\right]^{\left({I\choose r}\right)}$$ qui possède la propriété suivante : {\em soit $h$ appartenant à $\H^{q}(\hres(P),{\cal F})$. Ecrivons $h=\sum h_{i}$ au moyen de l'isomorphisme en question. Si $\cal E$ est un système de composantes adhérent à $P$ et si $\cal Q$ désigne le point fermé correspondant de la courbe résiduelle, alors l'image de $h$ dans $\H^{q}({\cal O}({\cal E})_{\rm alg},{\cal F})$ ({\em via} l'isomorphisme entre la cohomologie de $\hres(P)$ et celle de ${\cal O}_{X,P}$) est nulle si et seulement si la restriction à la cohomologie de $\red{\hres(P)}_{\cal Q}$ de chacune des $h_{i}$ est triviale.}~$\Box$} 

\bigskip
\deux{excohot2} {\bf Exemples.} Compte-tenu de la remarque~\ref{GRAT}.\ref{rgdivindep} la proposition ci-dessus admet les déclinaisons suivantes : si $|k^{*}|$ est un $\QQ$-espace vectoriel alors le groupe $\H^{q}(\hres(P),{\cal F})$ est isomorphe à $\H^{q}(\red{\hres(P)},{\cal F})$. Si $|k^{*}|$ est libre de rang 1 (autrement dit si $k$ est complet pour une valuation discrète) alors $\H^{q}(\hres(P),{\cal F})$ est isomorphe à $\H^{q}(\red{\hres(P)},{\cal F})\oplus \H^{q-1}(\red{\hres(P)},{\cal F}(-1))$. Dans le cas où $|k^{*}|$ est libre de rang 2 on trouve $$\H^{q}(\red{\hres(P)},{\cal F})\oplus (\H^{q-1}(\red{\hres(P)},{\cal F}(-1))\;)^{2}\;\oplus \H^{q-2}(\red{\hres(P)},{\cal F}(-2)),\;\mbox{\em etc.}$$

\bigskip
\deux{ouvuxp} Soit $P$ un point de type $(2)$ de $X$. Si $\cal U$ est un ouvert de Zariski non vide de la courbe résiduelle correspondante, si $\cal G$ est un $\red{k}$-faisceau étale et si $q$ est un entier on notera $\cha^{q}_{\cal U}(\red{\hres(P)},{\cal G})$ le noyau de $$\H^{q}(\red{\hres(P)},{\cal G})\to  \prod_{\tiny {\cal Q}\; \mbox{pt fermé de}\; {\cal U}} \H^{q}(\red{\hres(P)}_{\cal Q},{\cal G}).$$ Lorsque $\cal U$ est la courbe résiduelle toute entière on omettra de le rappeler en indice et l'on écrira simplement $\cha^{q}(\red{\hres(P)},{\cal G})$. Notons que $\cha^{0}_{\cal U}(\red{\hres(P)},{\cal G})$ est nul. 

\bigskip
Le lemme~\ref{GRAT}.\ref{trigratte} et la proposition~\ref{GRAT}.\ref{hqtype2} permettent aussitôt d'établir le résultat suivant : 

\bigskip
\deux{gptype2} {\bf Proposition.} {\em Soit $P$ un point de type $(2)$ de $X$. Le groupe $\G^{q}_{P}(X,{\cal F})$ est alors isomorphe à $$\bigoplus_{r\leq \min(q-1,|I|)} \left[\;\cha_{{\cal U}_{X,P}}^{q-r}(\red{\hres(P)},{\cal F}(-r))\;\right]^{\left({I\choose r}\right)}.$$ Si $P$ est intérieur (ce qui est par exemple toujours le cas si $X$ est lisse) on peut donc écrire $$\G^{q}_{P}(X,{\cal F})\simeq \bigoplus_{r\leq \min(q-1,|I|) }\left[\;\cha^{q-r}(\red{\hres(P)},{\cal F}(-r))\;\right]^{\left({I\choose r}\right)}.\;\Box$$}

\bigskip
\deux{exemplecha} {\bf Exemples.} Décrivons plus précisément $\G^{q}_{P}(X,{\cal F})$ {\em lorsque $P$ est intérieur} et sous certaines hypothèses particulières. 

\setcounter{cptbis}{0}
\bigskip
\trois{chadiv} Si $|k^{*}|$ est divisible alors $\G^{q}_{P}(X,{\cal F})$ est isomorphe à $\cha^{q}(\red{\hres(P)},{\cal F})$. 

\bigskip
\trois{chavaldisc} Si $|k^{*}|$ est libre de rang 1 alors $\G^{q}_{P}(X,{\cal F})$ est isomorphe à $$\cha^{q}(\red{\hres(P)},{\cal F})\oplus \cha^{q-1}(\red{\hres(P)},{\cal F}(-1)).$$ 

\bigskip
\trois{chaalgclos} Si $|k^{*}|$ est libre de rang 1 et si $\red{k}$ est séparablement clos alors ${\cal F}$ est constant. Le groupe $\G_{P}^{2}(X,{\cal F})$ est alors isomorphe à $\cha^{1}(\red{\hres(P)},{\cal F})$ et  partant à ${\cal F}^{2g}$, où $g$ est le genre de la courbe résiduelle. 

\bigskip
\deux{gpcourbelisse} {\bf Application au cas où $X$ est lisse.} Supposons que $X$ est lisse et soit $\bf S$ une triangulation de $X$. Soit ${\bf S}_{(2)}$ l'ensemble des points de type $(2)$ de $\bf S$.  La proposition~\ref{GRAT}.\ref{ouvtri} assure que pour tout ouvert $U$ de $X$ et tout entier $q$ les groupe $\G^{q}(U,{\cal F})$ et $\G^{q}_{c}(U,{\cal F})$ sont respectivement isomorphes à $\prod\limits_{P\in U\cap {\bf S}_{(2)}} \G^{q}_{P}(X,{\cal F})$ et à $\bigoplus\limits_{P\in U\cap {\bf S}_{(2)}} \G^{q}_{P}(X,{\cal F})$. Les groupes $\G^{q}_{P}(X,{\cal F})$ peuvent être décrits à l'aide de la proposition~\ref{GRAT}.\ref{gptype2} ci-dessus. 

\bigskip
\deux{coropasgratte} {\bf Corollaire.} {\em Soit $q$ un entier. Supposons que $X$ est lisse et possède une triangulation $\bf S$ telle que pour tout sommet $P$ de type $(2)$ et tout entier $r$ inférieur ou égal à $\min(q-1,|I|)$ le groupe $\cha^{q-r}(\red{\hres(P)},{\cal F}(-r))$ soit trivial. Alors $\G^{q}(U,{\cal F})$ et $\G^{q}_{c}(U,{\cal F})$ sont nuls pour tout ouvert $U$ de $X$.~$\Box$}
 
\bigskip
\deux{trpurtriv} Si $F$ est un corps quelconque et si $\cal G$ est un $F$-faisceau étale donné par un module galoisien fini de cardinal inversible dans $F$ alors $\cha^{q}(F(T),{\cal G})$ est trivial pour tout entier $q$ (\cf. \cite{jps}, \S 4 p. 122). Par ailleurs si $F$ est un corps fini, si  $n$ est un entier inversible dans $F$ et si $\cal C$ est une $F$-courbe projective, lisse et géométriquement intègre alors $\cha^{1}(F({\cal C}),\ZZ/n)$ et $\cha^{2}(F({\cal C}),\mu_{n})$ sont nuls. On en déduit les corollaires suivants :

\setcounter{cptbis}{0}

\bigskip
\trois{ratpasgratte} {\bf Corollaire.} {\em Supposons que $X$ est lisse et possède une triangulation $\bf S$ telle que pour tout sommet $P$ de type $(2)$ le corps $\red{\hres(P)}$ soit transcendant pur sur une extension finie de $\red{k}$. Alors $\G^{q}(U,{\cal F})$ et $\G^{q}_{c}(U,{\cal F})$ sont nuls pour tout ouvert $U$ de $X$ et tout entier $q$.~$\Box$}

\bigskip
\trois{mumpasgratte} {\bf Corollaire.} {\em Soit $\cal X$ une $k$-courbe algébrique de Mumford. Soit $U$ un ouvert de ${\cal X}\an$. Alors pour tout entier $q$ les groupes $\G^{q}(U,{\cal F})$ et $\G^{q}_{c}(U,{\cal F})$ sont nuls.}~$\Box$

\bigskip
\trois{ppasgratte}  {\bf Corollaire.} {\em Supposons que $|k^{*}|$ est libre de rang 1, que $\red{k}$ est fini et que $X$ est lisse. Soit $n$ un entier inversible dans $\red{k}$. Alors $\G^{3}(X,\ZZ/n(2))$ et $\G^{2}(X,\ZZ/n(1))$ sont nuls, ainsi que $\G^{3}_{c}(X,\ZZ/n(2))$ et $\G^{2}_{c}(X,\ZZ/n(1))$.~$\Box$}

\section{Deux théorèmes de dualité}\label{DUAL}

\subsection*{Un résultat de profinitude} 

\setcounter{cpt}{0}
\deux{notationdual} {\bf Le contexte.} On travaille toujours avec un corps ultramétrique complet $k$. On note $p$ l'exposant caractéristique de $\red{k}$. On choisit une clôture algébrique $k^{a}$ de $k$, on note $G$ le groupe de Galois de $k^{a}$ sur $k$ et $\ka$ le complété de $k^{a}$. On pourra se reporter à la définition~\ref{COH}.\ref{defmodfasc} concernant la notion de {\em faisceau raisonnable.}

\bigskip
On aura besoin de plusieurs assertions intermédiaires concernant la finitude de certains groupes de cohomologie. 

\bigskip
\deux{kacourbeepointe} {\bf Lemme.} {\em Soit $\cal Y$ une $\ka$-courbe algébrique projective et lisse. Donnons-nous une famille finie et non vide $(Y_{i})$ de domaines $k$-analytiques de ${\cal Y}\an$ satisfaisant les conditions suivantes : 

\bigskip
\begin{itemize}
\itb pour tout $i$ le domaine $Y_{i}$ est isomorphe ou bien à un disque ouvert, ou bien à un disque fermé contenu dans un domaine de ${\cal Y}\an$ lui-même isomorphe à un disque ouvert et ne rencontrant pas les autres $Y_{i}$. 

\bigskip
\itb les $Y_{i}$ sont deux à deux disjoints. 

\end{itemize}

\bigskip
Soit $Y$ le domaine complémentaire dans ${\cal Y}\an$ de la réunion des $Y_{i}$. Alors pour tout faisceau raisonnable $\cal F$ sur $Y$ et tout entier $q$ le groupe $\H^{q}(Y,{\cal F})$ est fini.}

\bigskip
{\em Démonstration.} Si $q$ est nul le résultat est évident puisque $Y$ est connexe. Si $q$ est au moins égal à $2$ le groupe $\H^{q}(Y,{\cal F})$ est nul : en effet la proposition 4.3.4 de \cite{brk2} permet de se ramener au cas où tous les $Y_{i}$ sont des disques fermés "prolongeables" en des disques ouverts deux à deux disjoints, auquel cas $Y$ est lisse et non compacte. Il n'y a plus alors qu'à appliquer la dualité de Poincaré pour les espaces $\ka$-analytiques lisses (\cite{brk2}, th. 7.4.3). Il reste le cas où $q$ vaut $1$ : le lemme 6.3.12 de \cite{brk2} permet cette fois de supposer que les $Y_{i}$ sont tous isomorphes à des disques {\em ouverts}. En vertu du théorème 6.3.9 et  du corollaire 6.3.11 de {\em loc. cit.} le lemme se déduit alors de la finitude de la cohomologie des $\ka$-courbes {\em algébriques.}~$\Box$

\bigskip
\deux{finicohotririg} {\bf Lemme.} {\em Soit $Y$ une $\ka$-courbe possédant une triangulation dont les sommets sont tous rigides et en nombre fini. Alors pour tout faisceau raisonnable $\cal F$ sur $Y$ et tout entier $q$ les groupes $\H^{q}(Y,{\cal F})$ sont finis.}

\bigskip
{\em Démonstration.} Si $P$ est un sommet alors l'ensemble des arêtes y aboutissant est fini puisque $P$ est rigide. En conséquence le nombre total d'arêtes est fini. Le lemme se déduit alors aussitôt de la première suite exacte du~\ref{COH}.\ref{introcohsup}, du lemme~\ref{PSEU}.\ref{cohogalois} et du fait que pour tout sommet $P$ le corps $\hres(P)$ n'est autre que $\ka$·~$\Box$

\bigskip
\deux{tribonrecouv} {\bf Proposition.} {\em Soit $X$ une $k$-courbe possédant une triangulation $\bf S$. Il existe alors un recouvrement de $X$ par une famille $(V_{i})$ d'ouverts, que l'on peut prendre finie si $\bf S$ est finie, et qui possède les propriétés suivantes : 

\bigskip
\begin{itemize}

\item[$i)$] pour tout couple $(i,j)$ d'indices (non nécessairement distincts), pour tout entier $q$ et pour tout faisceau raisonnable $\cal F$ sur $(V_{i}\cap V_{j})\times_{k}\ka$ le groupe $\H^{q}((V_{i}\cap V_{j})\times_{k}\ka,{\cal F})$ est fini ; 

\bigskip
\item [$ii)$] pour tout couple d'indices $(i,j)$ {\em distincts} l'intersection $V_{i}\cap V_{j}$ est une somme disjointe finie de pseudo-couronnes ; 

\bigskip
\item[$iii)$] si $(i,j,r)$ sont trois indices deux à deux distincts alors $V_{i}\cap V_{j}\cap V_{r}$ est vide ; 

\bigskip
\item[$iv)$] pour tout ensemble fini d'indices $I$ la réunion des $V_{i}$ pour $i$ appartenant à $I$ possède une triangulation finie. 

\end{itemize}
} 

\bigskip
{\em Démonstration.} Les arêtes qui aboutissent à un sommet donné sont toutes des pseudo-disques à l'exception d'un nombre fini d'entre elles. Puisque $X$ possède une triangulation son lieu de quasi-lissité est dense en vertu de la proposition~\ref{TRI}.\ref{triexistsing} ; en particulier si $P$ est un sommet non rigide de $\bf S$ alors $X$ est quasi-lisse en $P$, et ce dernier possède de ce fait un voisinage dans $X$  isomorphe à un domaine $k$-affinoïde d'une courbe algébrique lisse.  Ces remarques entraînent l'existence d'une famille $(V_{P})_{P\in {\bf S}}$ d'ouverts de $X$ possédant les propriétés suivantes : 

\bigskip
\begin{itemize}

\itb Si $P$ est de type $(2)$ ou $(3)$ alors $V_{P}$ contient toutes les arêtes de sommet $P$ qui sont des pseudo-disques, et $V_{P}\times_{k}\ka$ vérifie les conditions imposées à  $Y$ dans l'énoncé du lemme~\ref{DUAL}.\ref{kacourbeepointe}. 

\bigskip
\itb Pour tout sommet $P$ et toute arête $Y$ aboutissant à $P$ qui est une pseudo-couronne l'intersection $V_{P}\cap Y$ est une sous-pseudo-couronne stricte de $Y$. 

\bigskip
\itb Si $P$ et $Q$ sont deux points distincts de $\bf S$ et si $Y$ est une arête qui les relie (arête qui est alors nécessairement une pseudo-couronne) l'intersection $V_{P}\cap V_{Q}\cap Y$ est une sous-pseudo-couronne de $Y$.

\end{itemize} 

\bigskip
La famille formée des  $V_{P}$ et des arêtes qui sont des pseudo-couronnes dont l'adhérence n'est pas compacte satisfait alors les conditions requises grâce aux lemmes~\ref{PSEU}.\ref{cohogalois},~\ref{DUAL}.\ref{kacourbeepointe} et~\ref{DUAL}.\ref{finicohotririg}.$\Box$ 

\bigskip
\deux{invrecouvchange} {\bf Remarque.} Conservons les notations de la proposition ci-dessus ; alors si $L$ est une extension finie de $k$ le recouvrement de $X\times_{k}L$ par les $V_{i}\times_{k}L$ satisfait également les assertions $i)$, $ii)$, $iii)$ et $iv)$.  

\bigskip
\deux{fixen} On fixe à partir de maintenant un entier $n$ premier à $p$. Un faisceau étale en $\ZZ/n$-modules sera simplement appelé un {\em $\ZZ/n$-module.} Si l'on est sur un corps il sera dit {\em fini} s'il provient d'un module galoisien fini.

\bigskip
\deux{cohoprofin} {\bf Proposition.} {\em Supposons que pour tout $\ZZ/n$-module fini $\cal G$ sur $k$ et pour tout entier $q$ les groupes $\H^{q}(k,{\cal G})$ soient finis. Soit $X$ une $k$-courbe possédant un recouvrement par une famille d'ouverts $(V_{i})$ satisfaisant les conditions $i)$, $ii)$ et $iii)$ de la proposition~{\rm \ref{DUAL}.\ref{tribonrecouv}}. Pour tout ensemble fini $J$ d'indices notons $X_{J}$ la réunion des $V_{i}$ pour $i$ appartenant à $J$.  Soit $\cal F$ un $\ZZ/n$-module raisonnable sur $X$. Les propriétés suivantes sont alors vérifiées : 

\bigskip
\begin{itemize}
\item [$1)$] pour tout $J$ et tout $q$ les groupes $$\H^{q}(X_{J},{\cal F}),\;\H^{0}(|X_{J}|, \mbox{\rm R}^{q}\pi_{*}{\cal F})\;\mbox{et}\;\H^{1}(|X_{J}|,\mbox{\rm R}^{q}\pi_{*}{\cal F})$$ sont finis ;

\bigskip
\item [$2)$]  pour tout $q$ le groupe $$\H^{q}(X,{\cal F}), \; \mbox{(resp.}\;\H^{0}(|X|,\mbox{\rm R}^{q}\pi_{*}{\cal F}),\;\mbox{resp.}\;  \H^{1}(|X|,\mbox{\rm R}^{q}\pi_{*}{\cal F})\;\mbox{)}$$ s'identifie à $$\lim_{\leftarrow} \H^{q}(X_{J},{\cal F}),\;\mbox{(resp.}\; \lim_{\leftarrow}\H^{0}(|X_{J}|,\mbox{\rm R}^{q}\pi_{*}{\cal F}),\;\mbox{resp.}\; \lim_{\leftarrow}\H^{1}(|X_{J}|,\mbox{\rm R}^{q}\pi_{*}{\cal F})\;\mbox{)}$$ et est en particulier profini, et fini si la famille $(V_{i})$ est finie ;

\bigskip
\item [$3)$]  pour tout $J$ et tout $q$ la restriction $\H^{1}(|X|,\mbox{\rm R}^{q}\pi_{*}{\cal F})\to \H^{1}(|X_{J}|,\mbox{\rm R}^{q}\pi_{*}{\cal F}) $ est surjective.
\end{itemize}}

\bigskip
{\em Démonstration.} Commençons par une remarque. Donnons-nous un espace $k$-analytique $Y$ et notons $Y_{\ka}$ le $\ka$-espace $Y\times_{k}\ka$. Soit $\cal G$ un $\ZZ/n$-module raisonnable sur $Y$. De la suite spectrale $$\H^{i}(G,\H^{j}(Y_{\ka},{\cal G}))\Rightarrow \H^{i+1}(Y,{\cal G})$$ et de l'hypothèse faite sur $k$ on déduit aussitôt que si $\H^{i}(Y_{\ka},{\cal G})$ est fini pour tout $i$ alors $\H^{q}(Y,{\cal G})$ est fini pour tout $q$. Par ailleurs si $Y$ est une courbe on dispose quel que soit $q$ d'une suite exacte $$0\to \H^ {1}(|Y|,\mbox{R}^{q-1}\pi_{*}{\cal G})\to \H^{q}(Y,{\cal G})\to \H^{0}(|Y|,\mbox{R}^{q}\pi_{*}{\cal G})\to 0.$$  Dans ce cas si $\H^{q}(Y,{\cal G})$ est fini pour tout $q$ il en va de même de $\H^ {1}(|Y|,\mbox{R}^{q}\pi_{*}{\cal G})$ et $\H^{0}(|Y|,\mbox{R}^{q}\pi_{*}{\cal G})$.

\setcounter{cptbis}{0}

\bigskip
\trois{mayviet} Démonstration des assertions $1)$ et $2)$. En vertu de la remarque ci-dessus et des hypothèses faites sur le recouvrement $(V_{i})$ les groupes $\H^{q}(V_{i}\cap V_{j},{\cal F})$ et $\H^{0}(|V_{i}\cap V_{j}|,\mbox{R}^{q}\pi_{*}{\cal F})$ sont finis pour tout entier $q$ et tout couple $(i,j)$ d'indices (non nécessairement distincts). La propriété $iii)$ du recouvrement $(V_{i})$ entraîne l'existence d'une suite exacte de Mayer-Vietoris {\small $$\ldots\to \prod_{i} \H^{q-1}(V_{i},{\cal F})\to \prod_{i\neq j} \H^{q-1}(V_{i}\cap V_{j},{\cal F})\to \H^{q}(X,{\cal F})\to\prod_{i} \H^{q}(V_{i},{\cal F})\to \ldots$$} que l'on notera $(M)$. Soit $J$ un ensemble fini d'indices et soit  $(M_{J})$ la suite de Mayer-Vietoris associée au recouvrement $(V_{i})_{i\in J}$ de la réunion $X_{J}$ des $V_{i}$ pour $i$ appartenant à $J$. D'après ce qui précède le terme $\H^{q}(X_{J},{\cal F})$ de la suite $(M_{J})$ est encadré pour tout entier $q$ par deux groupes finis ; il est donc lui-même fini. On en déduit l'assertion $1)$ à l'aide grâce à la remarque ci-dessus. 

\bigskip
La limite projective $(\widehat{M})$ des $(M_{J})$ est  une limite projective de suites exactes de groupes finis, elle est en conséquence encore exacte. Or $(\widehat{M})$ se décrit simplement : elle est obtenue à partir de $(M)$ en remplaçant chacun des termes $\H^{q}(X,{\cal F})$ par $\lim\limits_{\leftarrow} \H^{q}(X_{J},{\cal F})$. On en déduit par chasse au diagramme que $$\H^{q}(X,{\cal F})\to \lim_{\leftarrow}\H^{q}(X_{J},{\cal F})$$ est un isomorphisme pour tout entier $q$, ce qui achève la démonstration de l'assertion $2)$ relative aux groupes $\H^{q}(X,{\cal F})$. 

\bigskip
Comme $\mbox{R}^{q}\pi_{*}{\cal F}$ est un faisceau le groupe $\H^{0}(|X|,\mbox{\rm R}^{q}\pi_{*}{\cal F})$ est isomophe pour tout $q$ à la limite projective des $\H^{0}(|X_{J}|,\mbox{\rm R}^{q}\pi_{*}{\cal F})$. Pour conclure prenons un entier $q$ et considérons le diagramme commutatif {\small $$\diagram 0\rto& \H^{1}(|X|,\mbox{R}^{q}\pi_{*}{\cal F})\rto \dto&\H^{q+1}(X,{\cal F})\rto\dto&\H^{0}(|X|,\mbox{R}^{q+1}\pi_{*}{\cal F})\dto \rto& 0\\ 0\rto& \lim\limits_{\leftarrow} \H^{1}(|X_{J}|,\mbox{R}^{q}\pi_{*}{\cal F})\rto &\lim\limits_{\leftarrow} \H^{q+1}(X_{J},{\cal F})\rto&\lim\limits_{\leftarrow}\H^{0}(|X_{J}|,\mbox{R}^{q+1}\pi_{*}{\cal F})\rto & 0\enddiagram$$} dans lequel les deux lignes sont exactes : pour celle du haut c'est évident, quant à celle du bas elle est la limite projective d'une famille de suites exactes courtes de groupes abéliens {\em finis} d'après ce qui précède ; elle est donc exacte. Par ailleurs on a vu ci-dessus que la flèche verticale du milieu et celle de droite sont des isomorphismes ; c'est en conséquence également le cas de celle de gauche, ce qui achève de démontrer l'assertion $2)$. 

\bigskip
\trois{hxhxjsurj} Etablissons la validité de l'assertion $3)$. Soit $q$ un entier. D'après les résultats de finitude que l'on vient d'établir il suffit de démontrer que pour tout ensemble fini d'indices $J$ et tout $i$ n'appartenant pas à $J$ la restriction $$\H^{1}(|X_{J}\cup V_{i}|,\mbox{\rm R}^{q}\pi_{*}{\cal F})\to \H^{1}(|X_{J}|,\mbox{\rm R}^{q}\pi_{*}{\cal F}) $$ est surjective. On se donne donc de tels $J$ et $i$. Soit $h$ appartenant au groupe $ \H^{1}(|X_{J}|,\mbox{\rm R}^{q}\pi_{*}{\cal F})$. Notons $\tilde{h}$ son image dans $\H^{q+1}(X_{J},{\cal F})$. Si $P$ est un point de $X_{J}$ alors $\tilde{h}(P)$ est nulle. Des propriétés $ii)$ et $iii)$ de la proposition~\ref{DUAL}.\ref{tribonrecouv} on déduit que $X_{J}\cap V_{i}$ est une somme disjointe et finie de pseudo-couronnes. Le lemme~\ref{GERM}.\ref{aretetriv} assure alors que $\tilde{h}_{|X_{J}\cap V_{i}}$ est nulle. La suite exacte de Mayer-Vietoris relative au recouvrement de $X_{J}\cup V_{i}$ par $X_{J}$ et $V_{i}$ fournit en particulier une classe $\eta$ dans $\H^{q+1}(X_{J}\cup V_{i},{\cal F})$  dont la restriction à $V_{i}$ (resp. à $X_{J}$) est nulle (resp. égale à $\tilde{h}$). L'image de $\eta$ dans $\H^{0}(|X_{J}\cup V_{i}|,\mbox{R}^{q+1}\pi_{*}{\cal F})$ est triviale ; en conséquence $\eta$ provient d'un unique élément $\omega$ de $\H^{1}(|X_{J}\cup V_{i}|,\mbox{\rm R}^{q}\pi_{*}{\cal F})$ dont la restriction à $X_{J}$ est par construction égale à $h$.~$\Box$ 

\bigskip
\deux{conclucohoprofin} {\bf Corollaire.} {\em Supposons que pour tout $\ZZ/n$-module fini $\cal G$ sur $k$ et pour tout entier $q$ les groupes $\H^{q}(k,{\cal G})$ soient finis. Soit $X$ une $k$-courbe, soit $\cal F$ un $\ZZ/n$-module raisonnable sur $X$ et soit $q$ un entier. Alors les groupes $\H^{q}(X,{\cal F})$, $\H^{0}(|X|,\mbox{\rm R}^{q}\pi_{*}{\cal F})$ et $\H^{1}(|X|,\mbox{\rm R}^{q}\pi_{*}{\cal F})$ sont profinis. Si $(X\times_{k}\ra{k})_{red}$ possède une triangulation finie, ce qui est en particulier le cas si $X$ est compacte ou si $X$ elle-même possède une triangulation finie, les groupes en question sont finis.}

\bigskip
{\em Démonstration.} On peut remplacer $k$ par $\ra{k}$ et donc le supposer parfait. On peut ensuite réduire $X$, et ainsi se ramener au cas où son lieu de quasi-lissité est dense. Elle possède alors une triangulation $\bf S$ par la proposition~\ref{TRI}.\ref{triexistsing}. Il n'y a plus qu'à appliquer les propositions~\ref{DUAL}.\ref{tribonrecouv} et~\ref{DUAL}.\ref{cohoprofin}.~$\Box$ 

\subsection*{Le cas où existe un module dualisant sur $k$} 

\bigskip
\deux{contextedual} {\bf A partir de maintenant et jusqu'à la fin du chapitre on fait les suppositions suivantes sur la cohomologie du corps $k$ :} 

\bigskip
\begin{itemize}
\item [$i)$] Pour toute extension finie séparable $L$ de $k$, pour tout $\ZZ/n$-module fini $\cal F$ sur $L$ et pour tout entier $q$ le groupe $\H^{q}(L,{\cal F})$ est fini.

\bigskip
\item[$ii)$] Il existe un $\ZZ/n$-module fini $\DD$ sur $k$, un entier $d$ et pour toute extension finie séparable $L$ un isomorphisme $\mbox{\rm Tr}_{L}:\H^{d}(L,\DD)\simeq \ZZ/n$, le système des $\mbox{\rm Tr}_{L}$ étant compatible aux corestrictions. 

\bigskip
\item[$iii)$] Pour tout entier $q$, pour toute extension finie séparable $L$ de $k$ et pour tout $\ZZ/n$-module fini $\cal F$ sur $L$ fini l'application $$\H^{q}(L,{\cal F})\times\H^{d-q}(L,{\cal F}^{\vee}\otimes\DD)\to \ZZ/n$$ induite par le cup-produit et la trace établit une dualité parfaite entre les groupes finis $\H^{q}(L,{\cal F})$ et $\H^{d-q}(L,{\cal F}^{\vee}\otimes\DD)$.

\end{itemize}

\bigskip
\deux{remdimcoho} {\bf Remarque.} Sous ces hypothèses pour toute extension finie séparable $L$ de $k$, pour tout $\ZZ/n$-module fini $\cal F$ sur $L$ et tout entier $q$ {\em strictement supérieur à $d$} le groupe $\H^{q}(L,{\cal F})$ est nul. 

\bigskip
\deux{dualpad} {\bf Exemples.} Si $|k^{*}|$ est libre de rang 1 et si $\red{k}$ est séparablement clos, ou bien si $k$ est fini (avec la valeur absolue triviale) alors on est dans cette situation avec $d$ égal à $1$ et en prenant pour $\DD$ le faisceau constant $\ZZ/n$. Si $|k^{*}|$ est libre de rang 1 et si $\red{k}$ est fini on est également dans ce cadre, avec cette fois-ci $d$ égal à $2$ et en prenant pour $\DD$ le faisceau $\mu_{n}$. 

\bigskip
\deux{intropoincare} Soit $X$ une $k$-courbe lisse. De la suite spectrale $$\H^{p}(G,\H^{q}_{c}(X_{\ka},\DD(1))\;)\Rightarrow \H^{p+q}_{c}(X,\DD(1))$$ et des hypothèses faites sur $k$ découle l'existence d'un isomorphisme entre les groupes $\H^{d+2}_{c}(X,\DD(1))$ et $\H^{d}(G,\H^{2}_ {c}(X_{\ka},\DD(1))\;)$. La trace définie dans le cas algébriquement clos (\cite{brk2}, th. 6.2.1) induit une flèche de $\H^{2}_ {c}(X_{\ka},\DD(1))$ vers $\DD$ puis un morphisme $$\mbox{\rm Tr}_{X} : \H^{d+2}_{c}(X,\DD(1))\simeq \H^{d}(G,\H^{2}_ {c}(X_{\ka},\DD(1))\;)\to \H^{d}(G,\DD)\simeq \ZZ/n.$$ Soit $X$ une $k$-courbe lisse. Soit $\cal F$ un $\ZZ/n$-module sur $X$. Notons $\breve{\cal F}$ le faisceau ${\cal F}^{\vee}\otimes \DD(1)$. La composée du cup-produit et de l'application $\mbox{Tr}_{X}$ définit pour tout entier $q$ un accouplement $$\H^{q}(X,{\cal F})\times \H^{d+2-q}_{c}(X,\breve{\cal F})\to \ZZ/n.$$

\bigskip
\deux{duapoincare}  {\bf Proposition (dualité de Poincaré).} {\em Le système de morphismes $(\mbox{\rm Tr}_{X})$ possède les propriétés suivantes : 

\bigskip
\begin{itemize}
\item[$i)$] Si $\phi : Y\to X$ est un morphisme plat et quasi-fini alors $\mbox{\rm Tr}_{Y}=\mbox{\rm Tr}_{X}\circ \mbox{\rm Tr}_{\phi}$ où $\mbox{\rm Tr}_{\phi}$ a été défini par Berkovich ({\rm\cite{brk2}, \S. 5.4}). En particulier si $\phi$ est une immersion ouverte alors $\mbox{\rm Tr}_{Y}=\mbox{\rm Tr}_{X}\circ \phi_{c}$.

\bigskip
\item[$ii)$] Pour toute $k$-courbe lisse $X$, pour tout $\ZZ/n$-module raisonnable $\cal F$ sur $X$ et tout entier $q$ l'accouplement $\H^{q}(X,{\cal F})\times \H^{d+2-q}_{c}(X,\breve{\cal F})\to \ZZ/n$ est continu par rapport à première variable pour la topologie profinie. Il induit une dualité parfaite entre le groupe profini $\H^{q}(X,{\cal F})$ et le groupe discret $\H^{d+2-q}_{c}(X,\breve{\cal F})$.

\end{itemize}}
 
\bigskip
{\em Démonstration.} Rappelons que le groupe $\H^{q}(X,{\cal F})$ est profini par le corollaire~\ref{DUAL}.\ref{conclucohoprofin}. L'assertion $i)$, ainsi que l'assertion $ii)$ lorsque les groupes $\H^{i}(X_{\ka},{\cal F})$ sont finis pour tout $i$, découlent formellement (comme dans le cas algébrique, \cf. \cite{sai}, lemme 2.9) des propriétés de la cohomologie de $k$ et de la dualité de Poincaré pour les courbes lisses sur un corps algébriquement clos (\cite{brk2}, th. 7.4.3) Traitons maintenant le cas général. Comme $X$ est lisse elle possède une triangulation $\bf S$. On peut donc lui appliquer la proposition~\ref{DUAL}.\ref{tribonrecouv} et obtenir ainsi une famille $(V_{i})$ d'ouverts recouvrant $X$ et assujettis à certaines conditions. Pour toute famille finie $J$ d'indices notons $(X_{J})$ la réunion des $V_{i}$ pour $i$ appartenant à $J$. La proposition~\ref{DUAL}.\ref{cohoprofin} assure alors (à l'aide de la remarque~\ref{DUAL}.\ref{invrecouvchange}) que pour tout $J$ et pour tout entier $i$ les groupes $\H^{i}(X_{J,\ka},{\cal F})$ sont finis ; la proposition est donc vraie en ce qui concerne chacun des $X_{J}$. Soit $q$ un entier. 

\bigskip
Soit $h$ une classe appartenant à $\H^{d+2-q}_{c}(X,\breve{\cal F})$. Elle provient d'une classe $\eta$ de $\H^{d+2-q}_{c}(X_{I},\breve{\cal F})$ pour un certain $I$. En vertu de l'assertion $i)$ une classe $\lambda$ de  $\H^{q}(X,{\cal F})$ appartient à l'orthogonal de $h$ si et seulement si $\lambda_{|X_{I}}$ appartient à l'orthogonal de $\eta$. Comme $\H^{q}(X,{\cal F})$ s'identifie à la limite projective des $\H^{q}(X_{J},{\cal F})$ d'après la proposition~\ref{DUAL}.\ref{cohoprofin} l'orthogonal de $h$ est un sous-groupe ouvert de $\H^{q}(X,{\cal F})$, ce qui montre la continuité par rapport à la première variable. 

\bigskip
Si $\lambda$ est un élément de $\H^{q}(X,{\cal F})$ appartenant à l'orthogonal de $\H^{d+2-q}_{c}(X,\breve{\cal F})$ alors $\lambda_{|X_{J}}$ appartient pour tout $J$ à l'orthogonal de $\H^{d+2-q}_{c}(X_{J},\breve{\cal F})$, toujours d'après l'assertion $i)$. L'assertion $ii)$ étant valable pour chacun des $X_{J}$ la classe $\lambda_{|X_{J}}$ est triviale quelque soit $J$. Comme $\H^{q}(X,{\cal F})$ est naturellement isomorphe à la limite projective des $\H^{q}(X_{J},{\cal F})$ la classe $\lambda$ est nulle. 

\bigskip
Soit $\chi$ un élément de $\H^{d+2-q}_{c}(X,\breve{\cal F})^{\vee}$. Pour tout $J$ la forme $\ZZ/n$-linéaire $\chi$ en induit une sur $\H^{d+2-q}_{c}(X_{J},\breve{\cal F})$, qui provient d'un {\em unique} $h_{J}$ appartenant à $\H^{q}(X_{J},{\cal F})$ puisque l'assertion $ii)$ est vraie pour $X_{J}$. En vertu de ladite unicité les $h_{J}$ définissent un élément de $\lim\limits_{\leftarrow} \H^{q}(X_{J},{\cal F})$ ; il existe donc une classe $h$ appartenant à $\H^{q}(X,{\cal F})$ tel que $h_{|X_{J}}$ soit égal à $h_ {J}$ pour tout $J$. Par construction la forme $\chi$ est égale à l'accouplement avec $h$.~$\Box$ 

\bigskip

{\em On fixe une $k$-courbe lisse $X$ ainsi qu'un $\ZZ/n$-module raisonnable $\cal F$ sur icelle.}

\bigskip
\deux{dualrqp} Le groupe $\H^{1}_{c}(|X|, \mbox{R}^{d+1}\pi_{*}\DD(1))$ est isomorphe à $\H^{d+2}_{c}(X,\DD(1))$, comme on le voit à l'aide des suites spectrales usuelles, et en se fondant sur des arguments de dimension cohomologique ; quant à $\H^{d+2}_{c}(X,\DD(1))$ il s'envoie lui-même dans $\ZZ/n$ par la trace. Soit $i$ un élément de $\{0,1\}$ et soit $q$ un entier. On peut définir grâce à ce qui précède et à l'aide du cup-produit un accouplement $$\H^{i}(|X|,\mbox{R}^{q}\pi_{*}{\cal F})\times \H^{1-i}_{c}(|X|,\mbox{R}^{d+1-q}\pi_{*}\breve{\cal F})\to \H^{1}_{c}(|X|, \mbox{R}^{d+1}\pi_{*}\DD(1))\to \ZZ/n.$$ Dans les énoncés qui suivent on considèrera $\H^{i}(|X|,\mbox{R}^{q}\pi_{*}{\cal F})$ comme muni de sa topologie profinie (voir le corollaire~\ref{DUAL}.\ref{conclucohoprofin}) ; rappelons par ailleurs que les groupes $\G^{q}(X,{\cal F})$ et $\G^{q}_{c}(X,{\cal F})$ ont été définis et étudié au chapitre~\ref{GRAT} (voir notamment la définition~\ref{GRAT}.\ref{intrograt}). 

\bigskip
\deux{maintheodua1} {\bf Théorème.} {\em L'accouplement $$\H^{0}(|X|,\mbox{\rm R}^{q}\pi_{*}{\cal F})\times\H^{1}_{c}(|X|,\mbox{\rm R}^{d+1-q}\pi_{*}\breve{\cal F})\to \ZZ/n$$ est continu par rapport à la première variable. Son noyau à droite est trivial  ; son noyau à gauche est le groupe $\G^{q}(X,{\cal F})$ qui est de ce fait profini. On a ainsi construit une dualité parfaite entre le groupe profini $$\H^{0}(|X|,\mbox{\rm R}^{q}\pi_{*}{\cal F})/\G^{q}(X,{\cal F})$$ et le groupe discret $\H^{1}_{c}(|X|,\mbox{\rm R}^{d+1-q}\pi_{*}\breve{\cal F})$.}

\bigskip
\deux{maintheodua2} {\bf Théorème.} {\em L'accouplement $$\H^{0}_{c}(|X|,\mbox{\rm R}^{q}\pi_{*}{\cal F})\times\H^{1}(|X|,\mbox{\rm R}^{d+1-q}\pi_{*}\breve{\cal F})\to \ZZ/n$$ est continu par rapport à la seconde variable. Son noyau à gauche est le groupe $\G^{q}_{c}(X,{\cal F})$. On a ainsi construit une dualité parfaite entre le groupe discret $$\H^{0}_{c}(|X|,\mbox{\rm R}^{q}\pi_{*}{\cal F})/\G^{q}_{c}(X,{\cal F})$$ et le groupe profini $\H^{1}(|X|,\mbox{\rm R}^{d+1-q}\pi_{*}\breve{\cal F})$.}

\bigskip
{\em Démonstration des deux théorèmes.} La courbe $X$ étant lisse elle possède une triangulation $\bf S$. La proposition~\ref{DUAL}.\ref{tribonrecouv} permet de construire une famille $(V_{i})$ d'ouverts recouvrant $X$ et vérifiant quatre conditions qu'on ne rappelle pas ici. Pour toute famille finie $J$ d'indices on va noter $X_{J}$ la réunion des $V_{i}$ pour $i$ appartenant à $J$. Les $X_{J}$ satisfont un certain nombre de propriétés données par la proposition~\ref{DUAL}.\ref{cohoprofin}. Pour chacun des deux théorèmes et pour tout ouvert $U$ de $X$ on notera $[.,.]_{U}$ (resp.  $<.,.>_{U}$) l'accouplement étudié (resp. l'accouplement fourni par la dualité de Poincaré, \cf. prop.~\ref{DUAL}.\ref{duapoincare}) sur la $k$-courbe lisse $U$. 

\bigskip
\setcounter{cptbis}{0}
\noindent
{\bf Démonstration du théorème~\ref{DUAL}.\ref{maintheodua1}}. Soient $h$ et $\eta$ appartenant respectivement à $\H^{0}(|X|,\mbox{\rm R}^{q}\pi_{*}{\cal F})$ et $\H^{1}_{c}(|X|,\mbox{\rm R}^{d+1-q}\pi_{*}\breve{\cal F})$. Soit $\lambda$ un antécédent de $h$ dans $\H^{q}(X,{\cal F})$ et soit $\tilde{\eta}$ l'image de $\eta$ dans $\H^{d+2-q}_{c}(X,\breve{\cal F})$. D'après le~\ref{GERM}.\ref{introaltercup} les éléments $[h,\eta]_{X}$ et $<\lambda,\tilde{\eta}>_{X}$ coïncident. En conséquence $h$ appartient à l'orthogonal de la classe $\eta$ si et seulement si $\lambda$ appartient à celui de $\tilde{\eta}$ ; la continuité de $[.,.]_{X}$ par rapport à la première variable suit aussitôt compte-tenu de l'assertion $ii)$ de la proposition~\ref{DUAL}.\ref{duapoincare} et du caractère ouvert de la surjection $\H^{q}(X,{\cal F})\to \H^{0}(|X|,\mbox{\rm R}^{q}\pi_{*}{\cal F})$. 

\bigskip
On voit de même que $\eta$ appartient au noyau à droite de $[.,.]_{X}$ si et seulement si $\tilde{\eta}$ appartient à celui de $<.,.>_{X}$. Or ce dernier est nul par la proposition~\ref{DUAL}.\ref{duapoincare}. Comme $\H^{1}_{c}(|X|,\mbox{\rm R}^{d+1-q}\pi_{*}\breve{\cal F})$ s'injecte dans $\H^{d+2-q}_{c}(X,\breve{\cal F})$ le noyau à droite de $[.,.]_{X}$ est bien trivial. 

\bigskip
Traitons maintenant le cas du noyau à gauche. On conserve les notations $h$ et $\lambda$ introduites ci-dessus. Pour toute arête $Y$ notons $j_{Y}$ l'immersion ouverte correspondante dans $X$. L'image de l'application $$\bigoplus_{\tiny \mbox{arêtes}} \H^{d+2-q}_{c}(Y,\breve{\cal F})\to \H^{d+2-q}_{c}(X,\breve{\cal F})$$ qui envoie $\sum \rho_{Y}$ sur $\sum j_{Y,c}(\rho_{Y})$ est égale à celle de l'injection $$\H^{1}_{c}(|X|,\mbox{\rm R}^{d+1-q}\pi_{*}\breve{\cal F})\hookrightarrow \H^{d+2-q}_{c}(X,\breve{\cal F})$$ d'après le lemme~\ref{COH}.\ref{nulsurs}. Donnons-nous un élément $\sum \rho_{Y}$ de $\bigoplus \H^{d+2-q}_{c}(Y,\breve{\cal F})$. De l'assertion $i)$ de la proposition~\ref{DUAL}.\ref{duapoincare} découle l'égalité $$\left<\lambda, \sum j_{Y,c}(\rho_{Y})\right>_{X}=\sum <\lambda_{|Y},\rho_{Y}>_{Y}.$$ En conséquence $h$ appartient au noyau à gauche de $[.,.]_{X}$ si et seulement si $\lambda_{|Y}$ appartient au noyau à gauche de $<.,.>_{Y}$ pour toute arête $Y$. Mais par la proposition~\ref{DUAL}.\ref{duapoincare} l'accouplement  $<.,.>_{Y}$ est non dégénéré quelle que soit $Y$ ; il s'ensuit que $h$ appartient à l'orthogonal de $\H^{1}_{c}(|X|,\mbox{\rm R}^{d+1-q}\pi_{*}\breve{\cal F})$ si et seulement si $\lambda_{|Y}$ est triviale pour tout $Y$. La nullité de $\lambda_{|Y}$ équivaut à celle de $h_{|Y}$ par la proposition~\ref{GERM}.\ref{aretetriv} ; et en vertu de la proposition~\ref{GRAT}.\ref{ouvtri} la classe $h$ est une classe gratte-ciel si et seulement si $h_{|Y}$ est nulle pour toute arête $Y$. Le théorème~\ref{DUAL}.\ref{maintheodua1} est donc démontré. 

\bigskip
\noindent
{\bf Démonstration du théorème~\ref{DUAL}.\ref{maintheodua2}.} Soient $h$ et $\eta$ appartenant respectivement à $\H^{0}_{c}(|X|,\mbox{\rm R}^{q}\pi_{*}{\cal F})$ et $\H^{1}(|X|,\mbox{\rm R}^{d+1-q}\pi_{*}\breve{\cal F})$. Soit $\lambda$ un antécédent de $h$ dans $\H^{q}_{c}(X,{\cal F})$ et soit $\tilde{\eta}$ l'image de $\eta$ dans $\H^{d+2-q}(X,\breve{\cal F})$. D'après le~\ref{GERM}.\ref{introaltercup} les éléments $[h,\eta]_{X}$ et $<\lambda,\tilde{\eta}>_{X}$ coïncident. En conséquence $\eta$ appartient à l'orthogonal de $h$ si et seulement si $\tilde{\eta}$ appartient à celui de $\lambda$ ; l'assertion $ii)$ de la proposition~\ref{DUAL}.\ref{duapoincare} permet de conclure aussitôt à la continuité de $[.,.]_{X}$ par rapport à la seconde variable. 

\bigskip
Il est clair d'après ce qui précède que $\eta$ appartient au noyau à droite de $[.,.]_{X}$ si et seulement si $\tilde{\eta}$ appartient à celui de $<.,.>_{X}$. Or ce dernier est nul par la proposition~\ref{DUAL}.\ref{duapoincare}. Comme $\H^{1}(|X|,\mbox{\rm R}^{d+1-q}\pi_{*}\breve{\cal F})$ s'injecte dans $\H^{d+2-q}(X,\breve{\cal F})$ le noyau à droite de $[.,.]_{X}$ est bien trivial. 

\bigskip
Traitons maintenant le cas du noyau à gauche. Soit $h$ appartenant au groupe $\H^{0}_{c}(|X|,\mbox{\rm R}^{q}\pi_{*}{\cal F})$. Elle provient d'une classe $\tilde{h}$ de la cohomologie à support compact de $X_{J}$ pour un certain $J$. La flèche  $$\H^{1}(|X|,\mbox{\rm R}^{d+1-q}\pi_{*}\breve{\cal F})\to \H^{1}(|X_{J}|,\mbox{\rm R}^{d+1-q}\pi_{*}\breve{\cal F})$$ est surjective par la proposition~\ref{DUAL}.\ref{cohoprofin} ;  on en déduit aussitôt que $h$ appartient au noyau à gauche de $[.,.]_{X}$ si et seulement $\tilde{h}$ appartient à celui de $[.,.]_{X_{J}}$. Comme il est par ailleurs clair que $h$ est une classe gratte-ciel si et seulement si c'est le cas pour $\tilde{h}$ on peut finalement supposer que $X$ est égal à $X_{J}$ et donc que la triangulation $\bf S$ est {\em finie}. Le groupe $\H^{0}_{c}(|X|,\mbox{\rm R}^{q}\pi_{*}{\cal F})$ est un sous-groupe de $\H^{0}(|X|,\mbox{\rm R}^{q}\pi_{*}{\cal F})$. Compte-tenu de la formule d'adjonction, du théorème~\ref{DUAL}.\ref{maintheodua1} déjà établi et du fait que $\G^{q}_{c}(X,{\cal F})$ est l'intersection de $\G^{q}(X,{\cal F})$ avec $\H^{0}_{c}(|X|,\mbox{\rm R}^{q}\pi_{*}{\cal F})$ il suffit de démontrer que la flèche canonique $$\H^{1}_{c}(|X|,\mbox{\rm R}^{d+1-q}\pi_{*}\breve{\cal F})\to \H^{1}(|X|,\mbox{\rm R}^{d+1-q}\pi_{*}\breve{\cal F})$$ est surjective. 

\bigskip
Soit $\eta$ appartenant à  $\H^{1}(|X|,\mbox{\rm R}^{d+1-q}\pi_{*}\breve{\cal F})$ et soit $\tilde{\eta}$ son image dans le groupe $\H^{d+2-q}(X,\breve{\cal F})$. Pour tout point de $P$ de $X$ la classe $\tilde{\eta}(P)$ est nulle. Dès lors la restriction de $\tilde{\eta}$ à toute arête est triviale par le lemme~\ref{GERM}.\ref{aretetriv}. Autrement dit $\tilde{\eta}_{|X-{\bf S}}$ est égale à zéro. Par injectivité de $$\H^{1}_{c}(|X-{\bf S}|,\mbox{\rm R}^{d+1-q}\pi_{*}\breve{\cal F})\to \H^{d+2-q}(X-{\bf S},\breve{\cal F})$$ la restriction de $\eta$ à $X-{\bf S}$ est nulle. Or $\bf S$ est fini, donc compact. Ceci achève la démonstration.~$\Box$ 

\bigskip
\deux{mumperfdual} {\bf Remarque.} En conservant les mêmes notations si les groupes $\G^{q}(X,{\cal F})$ et $\G^{q}_{c}(X,{\cal F})$ sont nuls alors les accouplements construits sont des dualités parfaites ; c'est par exemple toujours le cas d'après le corollaire~\ref{GRAT}.\ref{trpurtriv}.\ref{mumpasgratte}  si $\cal F$ provient d'un $\ZZ/n$-module galoisien fini sur $\red{k}$ et si $X$ est un ouvert de l'analytification d'une courbe de Mumford. 

\section{Quelques résultats revisités}\label{RE}

\setcounter{cpt}{0}
On se propose d'indiquer comment différents résultats déjà démontrés par d'autres méthodes se trouvent être des cas particuliers de ce qui précède. Ainsi le corollaire ci-dessous reprend les théorèmes 4.2 et 5.2 de \cite{duc}. 

\bigskip
\deux{h3mun2} {\bf Corollaire.} {\em Soit $k$ un corps ultramétrique complet. On suppose que $k$ est local, c'est-à-dire que $|k^{*}|$ est libre de rang $1$ et que $\red{k}$ est fini. Soit $n$ un entier inversible dans $\red{k}$. 

\bigskip

\begin{itemize} 
\item[$i)$] Soit $X$ une courbe $k$-analytique lisse et $\Delta$ son squelette ; orientons-le arbitrairement. Le groupe $\H^{0}(|X|,\mbox{\rm R}^{3}\pi_{*}\ZZ/n(2))$ s'identifie alors au groupe $\mbox{\rm Harm}(\Delta,\ZZ/n)$ des cochaînes harmoniques sur $\Delta$ à coefficients dans $\ZZ/n$. 

\bigskip
\item[$ii)$] Si $\cal X$ est une $k$-courbe algébrique lisse alors $$\H^{0}(|{\cal X}|,\mbox{\rm R}^{3}\pi_{*}\ZZ/n(2))\to \H^{0}(|{\cal X}\an|,\mbox{\rm R}^{3}\pi_{*}\ZZ/n(2))$$ est un isomorphisme.
\end{itemize}}

\bigskip
{\em Démonstration.} Établissons l'assertion $i)$. Le corollaire~\ref{GRAT}.\ref{trpurtriv}.\ref{ppasgratte} assure que $\G^{3}(X,\ZZ/n(2))$ est nul. Par ailleurs les hypothèses du~\ref{DUAL}.\ref{contextedual} sont vérifiées avec $d$ égal à $2$ et en prenant pour $\DD$ le module $\ZZ/n(1)$. En vertu du théorème~\ref{DUAL}.\ref{maintheodua1} le groupe profini $\H^{0}(X,\mbox{\rm R}^{3}\pi_{*}\ZZ/n(2))$ est alors isomorphe au dual de $\H^{1}_{c}(|X|,\ZZ/n)$, ce qui permet de conclure. 

\bigskip
Quant à l'assertion $ii)$ elle découle directement du corollaire~\ref{COMP}.\ref{corms}.~$\Box$ 

\bigskip
On peut également retrouver au niveau fini et $p$-primaire la dualité de Lichtenbaum entre le groupe de Picard et le groupe de Brauer d'une courbe $p$-adique projective et lisse (\cite{lic}). 

\bigskip
\deux{brauer} {\bf Corollaire.} {\em Soit $k$ un corps ultramétrique complet. On suppose que $k$ est local, c'est-à-dire que $|k^{*}|$ est libre de rang $1$ et que $\red{k}$ est fini. Soit $n$ un entier inversible dans $\red{k}$. 

\bigskip

\begin{itemize} 
\item[$i)$] Soit $X$ une courbe $k$-analytique lisse. Le groupe profini $\H^{0}(|X|,\mbox{\rm R}^{2}\pi_{*}\ZZ/n(1))$ est naturellement isomorphe au dual du groupe discret $\H^{1}_{c}(|X|,\mbox{\rm R}^{1}\pi_{*}\ZZ/n(1))$. 

\bigskip
\item[$ii)$] Soit $\cal X$ une $k$-courbe algébrique lisse. Les flèches $$\H^{0}(|{\cal X}|,\mbox{\rm R}^{2}\pi_{*}\ZZ/n(1))\to \H^{0}(|{\cal X}\an|,\mbox{\rm R}^{2}\pi_{*}\ZZ/n(1))$$ et $$\H^{1}(|{\cal X}|,\mbox{\rm R}^{1}\pi_{*}\ZZ/n(1))\to \H^{1}(|{\cal X}\an|,\mbox{\rm R}^{1}\pi_{*}\ZZ/n(1))$$ sont des isomorphismes. Si $\cal X$ est projective le groupe $_{n}\mbox{\rm Br}\;{\cal X}$  s'identifie naturellement au dual de $(\mbox{\rm Pic}\;{\cal X})/n$. 
\end{itemize}}

\bigskip
{\em Démonstration.} Établissons l'assertion $i)$. Le corollaire~\ref{GRAT}.\ref{trpurtriv}.\ref{ppasgratte} assure que $\G^{2}(X,\ZZ/n(1))$ est nul. Par ailleurs les hypothèses du~\ref{DUAL}.\ref{contextedual} sont vérifiées avec $d$ égal à $2$ et en prenant pour $\DD$ le module $\ZZ/n(1)$. Le théorème~\ref{DUAL}.\ref{maintheodua1} donne alors l'assertion $i)$. 

\bigskip
En ce qui concerne l'assertion $ii)$ le fait que les deux flèches soient des isomorphismes se déduit aussitôt du corollaire~\ref{COMP}.\ref{corms}. Considérons la suite exacte $$0\to \H^{1}(|{\cal X}|,\mbox{\rm R}^{1}\pi_{*}\ZZ/n(1))\to \H^{2}({\cal X},\ZZ/n(1))\to  \H^{0}(|{\cal X}|,\mbox{\rm R}^{2}\pi_{*}\ZZ/n(1)).$$ Par pureté un élément de $\H^{2}({\cal X},\ZZ/n(1))$ s'annule dans $\H^{0}(|{\cal X}|,\mbox{\rm R}^{2}\pi_{*}\ZZ/n(1))$ si et seulement si il s'annule aux points génériques de $\cal X$, ce qui équivaut à dire que son image dans $\mbox{\rm Br}\;{\cal X}$ est triviale. La suite exacte $$0\to (\mbox{\rm Pic}\;{\cal X})/n\to \H^{2}({\cal X},\ZZ/n(1)) \to _{n}\mbox{\rm Br}\;{\cal X}\to 0$$ fournit donc un isomorphisme entre $(\mbox{\rm Pic}\;{\cal X})/n$ et $\H^{1}(|{\cal X}|,\mbox{\rm R}^{1}\pi_{*}\ZZ/n(1))$ ainsi qu'un second entre $_{n} \mbox{\rm Br}\;{\cal X}$ et $\H^{0}(|{\cal X}|,\mbox{\rm R}^{2}\pi_{*}\ZZ/n(1))$. Lorsque $\cal X$ est projective le groupe $\H^{1}_{c}(|{\cal X}\an|,\mbox{\rm R}^{1}\pi_{*}\ZZ/n(1))$ est égal à $\H^{1}(|{\cal X}\an|,\mbox{\rm R}^{1}\pi_{*}\ZZ/n(1))$ et le résultat souhaité s'en déduit immédiatement.~$\Box$ 

\bigskip
Le corollaire ci-dessous figurait sans démonstration détaillée à la fin de \cite{duc} (pour $\cal F$ cyclique, ce qui est une restriction inoffensive) ; dans le cas projectif on peut le voir comme la reformulation à propos de la jacobienne de la courbe d'un théorème établi par Ogg dans \cite{ogg} pour toutes les variétés abéliennes (voir les explications de ce fait à la fin de \cite{duc}). 

\bigskip
\deux{ogg} {\bf Corollaire.} {\em Soit $k$ un corps ultramétrique complet tel que $|k^{*}|$ soit libre de rang 1 et tel que $\red{k}$ soit séparablement clos. 

\bigskip
\begin{itemize}

\item[$i)$] Soit $X$ une courbe $k$-analytique lisse, soit $E$ l'ensemble fermé et discret des points $P$ de type $(2)$ de $X$ tels que le genre $g(P)$ de la courbe résiduelle en $P$ soit strictement positif. Soit $\cal F$ un groupe abélien fini de cardinal inversible dans $\red{k}$. Soit $\Delta$ le squelette de $X$, arbitrairement orienté. On dispose alors d'une suite exacte de groupes profinis $$0\to \prod_{P \in E}{\cal F}^{2g(P)}\to \H^{0}(|X|, \mbox{\rm R}^{2}\pi_{*}{\cal F})\to \mbox{\rm Harm}(\Delta,{\cal F}) \to 0.$$

\bigskip
\item[$ii)$] Soit $\cal X$ une $k$-courbe algébrique lisse. La flèche $$\H^{0}(|{\cal X}|,\mbox{\rm R}^{2}\pi_{*}{\cal F})\to \H^{0}(|{\cal X}\an|,\mbox{\rm R}^{2}\pi_{*}{\cal F})$$ est un isomorphisme.
\end{itemize}}

\bigskip
{\em Démonstration.} Les hypothèses du~\ref{DUAL}.\ref{contextedual} sont vérifiées avec $d$ égal à $1$ et en prenant pour $\DD$ le module $\ZZ/n$. On déduit de la proposition~\ref{GRAT}.\ref{ouvtri} et du~\ref{GRAT}.\ref{exemplecha}.\ref{chaalgclos} que le groupe $\G^{2}(X,{\cal F})$ est naturellement isomorphe à $\prod \nolimits_{P\in E} {\cal F}^{2g(P)}$. Par le théorème~\ref{DUAL}.\ref{maintheodua1} le quotient $\H^{0}(|{\cal X}|,\mbox{\rm R}^{2}\pi_{*}{\cal F})/\G^{2}(X,{\cal F})$ s'identifie au dual de $\H^{1}_{c}(|X|,{\cal F}^{\vee})$ ; l'assertion $i)$ est ainsi établie. L'assertion $ii)$ se déduit du corollaire~\ref{COMP}.\ref{corms}.~$\Box$

\end{document}